\documentclass[12pt]{amsart}

\usepackage{dsfont}
\usepackage{amsxtra}
\usepackage{amsmath}
\usepackage{amscd}
\usepackage{amssymb}
\usepackage{amsfonts}
\usepackage[all]{xy}
\usepackage{mathrsfs}
\usepackage{amsthm} 
\usepackage{color}
\usepackage{upgreek }
\usepackage{tikz}
\usepackage{enumerate}
\usepackage{geometry}
\usepackage{ae}
\usepackage{enumitem}
\usepackage{blkarray}
\usepackage{multirow}
\usepackage{mathtools}
\usepackage{graphicx}

\numberwithin{itemcounter}{subsection}

\theoremstyle{plain}

\newtheorem{theorem}{Theorem}[section]
\newtheorem{lemma}[theorem]{Lemma}
\newtheorem{definition-lemma}[theorem]{Definition-Lemma}
\newtheorem{claim}[theorem]{Claim}
\newtheorem{proposition}[theorem]{Proposition}
\newtheorem{conjecture}[theorem]{Conjecture}
\newtheorem{corollary}[theorem]{Corollary}
\theoremstyle{definition}
\newtheorem{definition}[theorem]{Definition}
\theoremstyle{remark}
\newtheorem{remark}[theorem]{Remark}
\newtheorem{example}[theorem]{Example}
\numberwithin{equation}{section}

\def\bbA{\mathbb{A}}

\def\bbC{\mathbb{C}}

\def\bbF{\mathbb{F}}
\def\bbG{\mathbb{G}}

\def\bbN{\mathbb{N}}

\def\bbQ{\mathbb{Q}}

\def\bbZ{\mathbb{Z}}

\def\scrA{\mathscr{A}}

\def\scrC{\mathscr{C}}

\def\scrF{\mathscr{F}}

\def\scrH{\mathscr{H}}
\def\scrI{{I}}

\def\scrK{\mathscr{K}}
\def\scrL{\mathscr{L}}

\def\scrO{\mathscr{O}}
\def\scrP{\mathscr{P}}

\def\scrS{\mathscr{S}}

\def\scrU{\mathscr{U}}
\def\scrV{\mathscr{V}}

\def\scrX{\mathscr{X}}

\def\frakg{\mathfrak{g}}

\def\fraks{\mathfrak{s}}
\def\frakl{\mathfrak{l}}

\def\frakF{\mathfrak{F}}

\def\frakH{\mathfrak{H}}

\def\frakS{\mathfrak{S}}

\def\calL{\mathcal{L}}

\def\calO{\mathcal{O}}

\def\frakg{\mathfrak{g}}
\def\frakh{\mathfrak{h}}
\def\frakl{\mathfrak{l}}

\def\fraksl{\mathfrak{sl}}

\def\bfb{\mathbf{b}}

\def\bfe{\mathbf{e}}

\def\bfL{\mathbf{L}}

\def\bfB{\mathbf{B}}
\def\bfC{\mathbf{C}}

\def\bfF{\mathbf{F}}
\def\bfG{\mathbf{G}}
\def\bfH{\mathbf{H}}

\def\bfM{\mathbf{M}}
\def\bfN{\mathbf{N}}

\def\bfP{\mathbf{P}}

\def\bfS{\mathbf{S}}
\def\bfT{\mathbf{T}}
\def\bfU{\mathbf{U}}
\def\bfV{\mathbf{V}}
\def\bfW{\mathbf{W}}

\def\bfP{\mathbf{P}}
\def\bfT{\mathbf{T}}
\def\bfV{\mathbf{V}}

\def\geqs{\geqslant}

\def\simto{\overset{\sim}\to}

\def\bfSp{{\mathbf{Sp}}}

\def\bfSO{{\mathbf{SO}}}
\def\SO{\operatorname{SO}\nolimits}
\def\Sp{\operatorname{Sp}\nolimits}

\def\Unip{{\operatorname{Unip}\nolimits}}
\def\Irr{{\operatorname{Irr}\nolimits}}
\def\wIrr{{\operatorname{WIrr}\nolimits}}
\def\sp{{\operatorname{sp}\nolimits}}

\def\e{{\operatorname{e}\nolimits}}
\def\op{{\operatorname{op}\nolimits}}

\def\GL{\operatorname{GL}\nolimits}
\def\GU{\operatorname{GU}\nolimits}
\def\diag{\operatorname{diag}\nolimits}

\def\Frac{{\operatorname{Frac}\nolimits}}
\def\k{{\operatorname{k}\nolimits}}

\def\P{\operatorname{P}\nolimits}
\def\Q{\operatorname{Q}\nolimits}

\def\X{\operatorname{X}\nolimits}

\def\head{\operatorname{hd}\nolimits}
\def\cl{{\operatorname{cl}\nolimits}}
\def\lw{{\operatorname{\,lw}\nolimits}}
\def\up{{\operatorname{\,up}\nolimits}}

\def\soc{\operatorname{soc}\nolimits}
\def\rk{{\operatorname{rk}\nolimits}}

\def\min{{\operatorname{min}\nolimits}}

\def\Hom{\operatorname{Hom}\nolimits}

\def\End{\operatorname{End}\nolimits}
\def\Res{\operatorname{Res}\nolimits}
\def\res{\operatorname{res}\nolimits}
\def\Ind{\operatorname{Ind}\nolimits}

\def\Id{\operatorname{Id}\nolimits}

\def\wt{\operatorname{wt}\nolimits}
\def\hw{{\operatorname{hw}\nolimits}}

\def\height{\operatorname{ht}\nolimits}

\def\umod{\operatorname{-umod}\nolimits}
\def\mod{\operatorname{-mod}\nolimits}

\def\height{\operatorname{ht}\nolimits}

\newcommand{\bT}{{\mathbf{T}}}

\author{O. Dudas, M. Varagnolo, E. Vasserot}

\title[]{Categorical actions on unipotent representations\\
of finite classical groups}

\begin{document}
\begin{abstract}
We review the categorical representation of a Kac-Moody algebra
on unipotent representations of finite unitary groups in non-defining characteristic
given in \cite{DVV}, using Harish-Chandra induction and restriction.
Then, we extend this construction to finite reductive groups of types B or C in
non-defining characteristic. We show that the decategorified
representation is isomorphic to a direct sum of level 2 Fock spaces. 
We deduce that the Harish-Chandra branching graph coincides with the crystal graph of
these Fock spaces. We also obtain derived equivalences between blocks, yielding
Brou\'e's abelian defect group conjecture for unipotent $\ell$-blocks at linear primes
$\ell$. 
\end{abstract}

\thanks{This research was partially supported by the ANR grant number 
ANR-13-BS01-0001-01.}

\maketitle

\setcounter{tocdepth}{3}

\tableofcontents

\section*{Introduction}

Let $\bfG$ be a connected algebraic group defined over a finite field $\bbF_q$. The finite group $G =\bfG(\bbF_q)$ of its $\bbF_q$-rational points is a \emph{finite reductive group}. The irreducible representations
of $G$ over fields of characteristic $\ell$ prime to $q$ fall into Harish-Chandra series, which 
are defined in terms of Harish-Chandra induction $R_L^G$ and restriction ${}^*R_L^G$
from proper Levi subgroups $L$ of $G$. The isomorphism classes in each series are parametrized
by the simple modules of the \emph{ramified Hecke algebras}, which are realized as the
endomorphism algebras of the Harish-Chandra induction of cuspidal representations.
Therefore, the classification of isomorphism classes of irreducible representations of $G$
can be reduced to the following two problems :
\begin{itemize}
  \item[(a)] classification of the cuspidal irreducible representations,
  \item[(b)] determination of the ramified Hecke algebras.
\end{itemize}

\smallskip

This was achieved by Lusztig in \cite{Lu84} when $\ell = 0$ but it remains open for
representations in positive characteristic for most of the finite reductive
groups. By results of Geck-Hiss-Malle \cite{GHM96}, we know however that the ramified Hecke algebras 
are indeed Hecke algebras of finite type, only the parameters of the deformation are unknown in general. 

\smallskip

When $G$ is a classical group, it turns out that most of the structure of the ramified Hecke algebras does not
depend on the corresponding unipotent cuspidal representations. This suggests to rather study the endomorphism
algebra of the Harish-Chandra induction functor $R_L^G$ rather
than the endomorphism algebra of the induced representation. This was already
achieved in \cite{CR} for $G= \mathrm{GL}_n(q)$. Our goal is to extend 
 Chuang-Rouquier's approach to other classical groups. 

\smallskip

In the first part of this paper we focus on the case of finite unitary
groups $\mathrm{GU}_n(q)$, reviewing our previous work \cite{DVV}, whereas
in the last part we deal with the case of groups of types B or C.
We will work with both ordinary representations (characteristic zero) and
modular representations in non-defining characteristic (characteristic $\ell$ prime to $q$). More precisely, the field of coefficients $R$ of the representations will
be an extension of either $\bbQ_\ell$ or $\bbF_\ell$. 

\smallskip

Let $\mathrm{G}_n(q)$ be one of the families of finite classical
groups among $\mathrm{GU}_{2n}(q)$, $\mathrm{GU}_{2n+1}(q)$, 
$\mathrm{Sp}_{2n}(q)$ and $\mathrm{SO}_{2n+1}(q)$.
Using the tower of inclusion of groups
$\cdots \subset \mathrm{G}_n(q) \subset  \mathrm{G}_{n+1}(q) \subset \cdots$
one can form the abelian category of unipotent representations
$$ \scrU_R := \bigoplus_{n \geqslant 0} \mathrm{G}_n(q)\umod.$$
Furthermore, under mild assumption on $\ell$, we can modify the Harish-Chandra
induction and restriction functors to obtain an adjoint pair $(E,F)$ of functors
on $\scrU_R$. The functor $F$ corresponds to a Harish-Chandra induction from
$\mathrm{G}_n(q)$ to $\mathrm{G}_{n+1}(q)$ whereas $E$ corresponds to the restriction.
Note that only specific Levi subgroups are considered, and we must work with a
variation of the usual Harish-Chandra theory (the \emph{weak Harish-Chandra theory})
introduced in \cite{GHJ}. 

\smallskip

In this framework, problem (a) amounts to finding the modules $V$ such that
$EV=0$ and problem (b) is about the structure of $\mathrm{End}_G(F^m V)$
for such cuspidal modules $V$. As mentioned before, most of the structure
of this endomorphism algebra is already contained in $\mathrm{End}(F^m)$. 
In \S\ref{sec:rep-datum} and \S\ref{sec:repdatumB}, we construct natural
transformations $X$ of $F$ and $T$ of $F^2$ which
endow $\mathrm{End}(F^m)$ with a morphism from an affine Hecke algebra
$\bfH_{m}^{q^\delta}$ of type $A_{m-1}^{(1)}$ with parameter $q^\delta$ where the integer $\delta$ is determined as follows
\begin{itemize}
\item $\delta=2$ if $\mathrm{G}_n(q)$ is $\mathrm{GU}_{2n}(q)$ or $\mathrm{GU}_{2n+1}(q),$
\item $\delta=1$ if $\mathrm{G}_n(q)$ is $\mathrm{Sp}_{2n}(q)$ or $\mathrm{SO}_{2n+1}(q).$
\end{itemize}
Back to our original problem, the evaluation at a cuspidal module $V$ provides
a natural map $\bfH_{m}^{q^\delta} \to \mathrm{End}_G(F^mV)$.
Then, we prove that this map induces a natural isomorphism between
$\mathrm{End}_G(F^m V)$ and a level 2 cyclotomic quotient of $\bfH_{m}^{q^\delta}$ whose parameters are naturally attached to $V$, see Theorems \ref{thm:cat} and
\ref{thm:HL-BC}.

\smallskip

Next, we prove that the eigenvalues of $X$ belong to
\begin{itemize}
\item $I_2=(-q)^\bbZ$ if $\mathrm{G}_n(q)$ is $\mathrm{GU}_{2n}(q)$ or $\mathrm{GU}_{2n+1}(q),$
\item $I_1=q^\bbZ\sqcup(-q^\bbZ)$ if $\mathrm{G}_n(q)$ is $\mathrm{Sp}_{2n}(q)$ or $\mathrm{SO}_{2n+1}(q).$
\end{itemize}
Then, we can form a Lie algebra $\frakg$ corresponding to the quiver
with vertices $I_\delta$ and arrows given by multiplication by $q^\delta$.
When working in characteristic zero, the Lie algebra $\frakg$ is isomorphic to
two copies of $\fraks\frakl_\bbZ$. In positive characteristic
$\ell$ prime to $q$ it will depend on whether $-q$ is a power of $q^\delta$. This reflects
the difference of behaviour for unipotent representations when working at
linear or unitary primes. To explain this we denote by $d,e,f$ the order of
$q^2$, $-q$ and $q$ modulo $\ell$. Two situations occur:
\begin{itemize}
 \item $\ell$ is a \emph{linear prime} if $\mathrm{G}_n(q)$ is $\mathrm{GU}_{2n}(q)$
 or $\mathrm{GU}_{2n+1}(q)$ and $e$ is even (then $e=2d$ and $-q\notin q^{2\bbZ}$ mod $\ell$), or if
 $\mathrm{G}_n(q)$ is $\mathrm{Sp}_{2n}(q)$ or $\mathrm{SO}_{2n+1}(q)$
 and $f$ is odd (then $f=d$ and $-q\notin q^\bbZ$ mod $\ell$). 
 In that case $\frakg$ is a subalgebra of  $(\widehat{\fraks\frakl}_{d})^{\oplus2}$;

 \item $\ell$ is a \emph{unitary prime} if $\mathrm{G}_n(q)$ is $\mathrm{GU}_{2n}(q)$
 or $\mathrm{GU}_{2n+1}(q)$ and $e$ is odd (then $e=d$ and $-q\equiv(q^2)^{(e+1)/2}$ mod $\ell$), or if
 $\mathrm{G}_n(q)$ is $\mathrm{Sp}_{2n}(q)$ or $\mathrm{SO}_{2n+1}(q)$
 and $f$ is even (then $f=2d$ and $-q\equiv q^{d+1}$ mod $\ell$). In that case $\frakg$ is isomorphic to $\widehat{\fraks\frakl}_{f/\delta}$.
\end{itemize}
In each case we prove that the representation datum $(E,F,X,T)$
induces a categorical action of $\frakg$ on $\scrU_R$.
See Theorems \ref{thm:char0} and \ref{thm:charl} for finite unitary groups, and Theorems
\ref{thm:B}, \ref{thm:BB} and \ref{thm:BBB} for groups of types B and C. 

\smallskip

In particular, let $E = \bigoplus E_i$ and $F = \bigoplus F_i$ be the decomposition of
the functors into generalized $i$-eigenspaces for $X$. Then $[E_i],$
$[F_i]$ act as the Chevalley generators of $\frakg$
on the Grothendieck group $[\scrU_R]$ of $\scrU_R$ and many
problems on $\scrU_R$ have a Lie-theoretic counterpart. For example,
\begin{itemize}[leftmargin=8mm]
  \item weakly cuspidal modules correspond to highest weight vectors,
  \item the decomposition of $\scrU_R$ into weak Harish-Chandra series corresponds
  to the decomposition of the $\frakg$-module $[\scrU_R]$ into a direct sum
  of irreducible highest weight modules,
  \item the parameters of the ramified Hecke algebra attached to a weakly cuspidal unipotent module $V$
  are given by the weight of $[V]$,
  \item the blocks of $\scrU_R$, or equivalently the unipotent $\ell$-blocks, correspond
  to the weight spaces for the action of $\frakg$ (inside a Harish-Chandra series
  if $\ell$ is a linear prime). 
\end{itemize}
Such observations were already used in other situations, e.g., 
for cyclotomic rational double affine Hecke algebras.

\smallskip

For this dictionary to be efficient one needs to determine the $\frakg$-module
structure on $[\scrU_R]$. This is done in \S \ref{sec:g-e} and \S \ref{sec:uk-typeBC}
by looking at the action
of $[E_i]$ and $[F_i]$ on the basis of $[\scrU_R]$ formed by unipotent characters
and their $\ell$-reduction, which play a role similar to the role of the  \emph{standard  modules} in the categorifications of
 cyclotomic rational double affine Hecke algebras mentioned above.
On this basis the action can be made explicit, and we prove that there is a natural
$\frakg$-module isomorphism
$$  [\scrU_R] \, \mathop{\longrightarrow}\limits^\sim \, \bigoplus_{t \geqslant 0} \bfF(Q_t)$$
between the Grothendieck group of $\scrU_R$ and a direct sum of
level 2 Fock spaces $\bfF(Q_t)$, each of which corresponds to an ordinary
Harish-Chandra series. Through this
isomorphism, the basis of unipotent characters (or their $\ell$-reduction) is
sent to the standard monomial basis. 

\smallskip

Our original motivation for constructing a categorical action of $\frakg$
on $\scrU_R$ comes from a conjecture of Gerber-Hiss-Jacon \cite{GHJ} for
finite unitary groups, which
predicts an explicit relation between the Harish-Chandra branching graph
and the crystal graph of the Fock spaces $\bigoplus_{t \geqslant 0} \bfF(Q_t)$ when $e$ is odd.
See also  \cite{GH}.
Using our categorical methods and the unitriangularity of the 
decomposition matrix we obtained in \cite{DVV} a complete proof of the conjecture,
see Theorem \ref{thm:wHC}. We extend here this result to groups of type B
and C (see Theorem \ref{thm:crystalBC}). Note however that in order to do so,
we must rely on a conjecture of Geck concerning the unitriangular shape of the decomposition 
matrices that we have recalled in Conjecture \ref{conj:unitriangular}.
\smallskip

A similar result can be deduced when $\ell$ is a linear prime. However, in 
that case, the situation is already well-understood by the work 
of Gruber-Hiss \cite{GrH} on classical groups. The case where 
$\ell$ is unitary is considered as more challenging
and our categorification techniques give the first major result in that direction since
the case of $\mathrm{GL}_n(q)$ was solved by Dipper-Du
\cite{DiDu}. This solves completely the problem of classification of
irreducible unipotent modules for unitary groups and groups of types B or C
mentioned at the beginning of the introduction, and yields a combinatorial 
description of the weak Harish-Chandra series.

\smallskip

For finite unitary groups, another categorical construction can be used in order to get the usual (non weak) Harish-Chandra series, 
by adapting some 
techniques from \cite{SV}. It relies on a categorification of the Heinsenberg representation of level 2 on the Fock spaces
$\bfF(Q_t)$ mentioned above. We explain this in \S \ref{sec:heisenberg}.
This yields a complete classification of the cuspidal unipotent modules, as explained in Theorem 
\ref{thm:cuspidal} and \S \ref{sec:FLOTW} (and a computation of the parameters of all ramified Hecke algebras). 

\smallskip

By the work of Chuang-Rouquier, categorical actions also provide
derived equivalences between weight spaces. In our situation, these
weight spaces are exactly the unipotent $\ell$-blocks and we obtain
many derived equivalences between blocks with the same local structure.
Together with Livesey's construction of good blocks in the linear
prime case, we deduce a proof of Brou\'e's abelian defect group
conjecture see Theorems  \ref{thm:broue} (or \cite{DVV}) and \ref{thm:derivedBC}
(which does not rely on Conjecture  \ref{conj:unitriangular}).

\smallskip

The paper is organized as follows. In Section \ref{part:categoricalactions}
we set our notations and recall basic facts on categorical
actions, perfect bases and derived 
equivalences. In Section \ref{part:fock} we introduce the Fock spaces,
which are certain integrable representations of Kac-Moody algebras.
They have a crystal graph which can be defined combinatorially. In Section \ref{part:modular}
we recall standard results on unipotent representations of finite
reductive groups in non-defining characteristic. In Section \ref{part:unitarygroups}
we recall the categorical representation on the unipotent modules of unitary groups given in \cite{DVV},
and the main applications proved there.
In Section \ref{sec:heisenberg} we explain the role of the Heisenberg categorical action and its relation with
Harish-Chandra series. The results here are new and were announced in \cite{DVV}.
Section \ref{sec:BC} deals with the categorical representation on the unipotent modules of groups of type B,C
and its applications to weak Harish-Chandra series and derived equivalences.
The results in this final section are new.

\medskip

\section{Categorical representations}\label{part:categoricalactions}

Throughout this section, $R$ will denote a noetherian commutative domain (with unit).

\subsection{Rings and categories\label{ssec:rings-cat}}\hfill\\

An \emph{$R$-category} $\scrC$ is an additive category enriched over the tensor
category of $R$-modules. All the functors $F$ on $\scrC$ will be assumed to be $R$-linear.
Given such a functor, we denote by $1_F$ or sometimes $F$ the identity element in the endomorphism
ring $\End(F)$. The identity functor on $\scrC$ will be denoted by by $1_\scrC$.
A composition of functors $E$ and $F$ is written as $EF$, while a composition of
morphisms of functors (or natural transformations) $\psi$ and $\phi$ is written as
$\psi\circ\phi$. We say that $\scrC$ is \emph{Hom-finite} if the Hom spaces are finitely
generated over $R$.  If the category $\scrC$ is abelian or exact, we denote by $[\scrC]$ the
complexified Grothendieck group and by $\Irr(\scrC)$ the set of isomorphism classes of
simple objects of $\scrC$. The class of an object $M$ of $\scrC$ in the Grothendieck group
is denoted by $[M]$. An exact endofunctor $F$ of $\scrC$ induces a linear map
on $[\scrC]$ which we will denote by $[F]$.

\smallskip

Assume that $\scrC$ is Hom-finite. Given an object $M\in\scrC$ we set $\scrH(M)=
\End_\scrC(M)^\text{op}$. It is an $R$-algebra which is finitely generated as an $R$-module.
We denote by $\frakF_M$ the functor
$\frakF_M=\Hom_\scrC(M,-):\scrC\longrightarrow\scrH(M)\mod.$

\smallskip

Assume now that $\scrC=H\mod$, where $H$ is an $R$-algebra with 1 which is finitely
generated and free over $R$. We abbreviate $\Irr(H)=\Irr(\scrC)$. Given an homomorphism
$R\to S$, we can form the $S$-category $S\scrC=SH\mod$ where $SH=S\otimes_R H$. 
Given another $R$-category $\scrC'$ as above and an exact functor
$F:\scrC\to\scrC'$, then $F$ is represented by a projective object $P\in\scrC$. 
We set $SF=\Hom_{S\scrC}(SP,-):S\scrC\to S\scrC'$. 

\smallskip

Let $K$ be the field of fractions
of $R$, $A\subset R$ be a subring which is integrally closed in $K$ and $\theta:R\to \k$
be a ring homomorphism into a field $\k$ such that $\k$ is the field of fractions of
$\theta(A)$. If $\k H$ is split, then there is a \emph{decomposition map} $d_\theta\,:
\,[KH\mod]\longrightarrow[\k H\mod]$, see e.g. \cite[sec.~3.1]{GJ} for more details.

\medskip

\subsection{Kac-Moody algebras of type $A$ and their representations}\label{sec:quivers}\hfill\\

The Lie algebras which will act on the categories we will study will always be finite
sums of Kac-Moody algebras of type $A_\infty$ or $A_{e-1}^{(1)}$. They will arise
from quivers of the same type. 

\smallskip

\subsubsection{Lie algebra associated with a quiver}\label{subsec:quivers}
Let $v \in R^\times$ and $\scrI \subset R^\times$. We assume that $v \neq 1$ and that
$\scrI$ is stable by multiplication by $v$ and $v^{-1}$ with finitely many orbits.
To the pair $(\scrI,v)$ we associate a quiver $\scrI(v)$ (also denoted by $\scrI$)
as follows:
\begin{itemize}[leftmargin=8mm]
  \item the vertices of $\scrI(v)$ are the elements of $\scrI$;
  \item the arrows of $\scrI(v)$ are $i \to iv$ for $i \in \scrI$. 
\end{itemize}
Since $\scrI$ is assumed to be stable by multiplication by $v$ and $v^{-1}$, such a
quiver is the disjoint union of quivers of type $A_\infty$ if $v$ is not a root of unity,
or of cyclic quivers of type $A_{e-1}^{(1)}$ if $v$ is a primitive $e$-th root of $1$.

\smallskip

The quiver $\scrI(v)$ defines a symmetric generalized Cartan matrix
$A = (a_{ij})_{i,j\in \scrI}$ with $a_{ii}= 2$, $a_{ij} =-1$ when
$i \rightarrow j$ or $j \rightarrow i$ and $a_{ij}=0$ otherwise. To this Cartan matrix one can
associate the (derived) Kac-Moody algebra $\frakg_\scrI'$ over $\bbC$, which
has Chevalley generators $e_i,f_i$ for $i\in \scrI$, subject to the usual
relations. 

\smallskip

More generally, let $(\X_\scrI,\X_\scrI^\vee,\langle\bullet,\bullet\rangle_\scrI,
\{\alpha_i\}_{i \in \scrI}, \{\alpha_i^\vee\}_{i \in \scrI})$ be a \emph{Cartan datum} associated with $A$, \emph{i.e.}, 
we assume that
\begin{itemize}[leftmargin=8mm]
 \item $\X_\scrI$ and $\X^\vee_\scrI$ are free abelian groups,
 \item the simple coroots $\{\alpha_i^\vee\}$ are linearly independent in $\X_\scrI^\vee$, 
 \item for each $i\in \scrI$ there exists a fundamental weight $\Lambda_i\in\X_\scrI$
 satisfying $\langle \alpha_j^\vee,\Lambda_i \rangle_\scrI = \delta_{ij}$ for all $j\in\scrI$,
 \item $\langle \bullet , \bullet \rangle_\scrI : \X_\scrI^\vee\times \X_\scrI 
 \longrightarrow \bbZ$ is a perfect pairing such that $\langle \alpha_j^\vee,\alpha_i\rangle_\scrI =~a_{ij}$. 
\end{itemize}
Let $\Q_\scrI^\vee = \bigoplus \bbZ \alpha_i^\vee$ be the coroot lattice and $\P_\scrI = \bigoplus \bbZ \Lambda_i$ be the weight lattice.
Then, the Kac-Moody algebra $\frakg_\scrI$ corresponding to this datum is the Lie algebra generated by
the Chevalley generators $e_i, f_i$ for $i \in \scrI$ and the Cartan algebra
$\frakh = \bbC\otimes \X_\scrI^{\vee}$. An element $h \in \frakh$ acts
by $[h,e_i] = \langle h,\alpha_i\rangle e_i$.
The Lie algebra $\frakg_\scrI'$ is the derived 
subalgebra $[\frakg_\scrI,\frakg_\scrI]$.  

\begin{example}\label{ex:onegenerator}
When $\scrI = v^{\bbZ}$ two cases arise.
 \begin{itemize}[leftmargin=8mm]
  \item[(a)] If $\scrI$ is infinite, then $\frakg_\scrI'$ is isomorphic to 
  $\mathfrak{sl}_\bbZ$, the Lie algebra of traceless matrices with finitely many
  non-zero entries.
  \item[(b)] If $v$ has finite order $e$, then $\scrI$ is isomorphic
  to a cyclic quiver of type $A_{e-1}^{(1)}$. We can form $\X^\vee = \Q^\vee 
  \oplus \bbZ \partial$ and $\X = \P \oplus \bbZ \delta$ with $\langle \partial,
  \Lambda_i\rangle = 0$, $\langle \partial,\alpha_i\rangle = \delta_{i1}$ and
  $\delta = \sum_{i \in \scrI} \alpha_i$. The pairing is non-degenerate, and
  $\frakg_\scrI$ is isomorphic to the Kac-Moody algebra
  $$\widehat{\fraks\frakl}_e = \fraks\frakl_e(\bbC) \otimes \bbC[t,t^{-1}] \oplus
  \bbC c \oplus \bbC \partial.$$
  An explicit isomorphism sends $e_{v^i}$ (resp. $f_{v^i}$) to the matrix
  $E_{i,i+1} \otimes 1$ (resp. $E_{i+1,i} \otimes 1$) if $i \neq e$ and $e_1$ 
  (resp. $f_1$) to $E_{e,1}\otimes t$ (resp. $E_{1,e}\otimes t^{-1}$).
  Via this isomophism the central element $c$ corresponds to $\sum_{i \in \scrI}
  \alpha_i^\vee$, and the derived algebra $\frakg_\scrI'$
  to $\widetilde{\fraks\frakl}_e=\fraks\frakl_e(\bbC) \otimes \bbC[t,t^{-1}] \oplus \bbC c$.
 \end{itemize}
\end{example}

To avoid cumbersome notation, we may write $\frakg=\frakg_\scrI$, $\P=\P_{\scrI}$,
$\Q^\vee=\Q_{\scrI}^\vee$, etc. when there is no risk of confusion.

\smallskip

\subsubsection{Integrable representations}\label{subsec:oint}
Let $V$ be a $\frakg$-module. Given $\omega \in \X$, the $\omega$-\emph{weight space}
of $V$ is 
  $V_\omega= \{v\in V\,\mid\,\alpha^\vee\cdot v=\langle\alpha^\vee,\omega\rangle\,v,\,
  \forall \alpha^\vee\in\Q^\vee\}.$
We denote by $\calO^\text{int}$ the category of \emph{integrable highest weight} modules,
\emph{i.e.} $\frakg$-modules $V$ satisfying
\begin{itemize}[leftmargin=8mm]
  \item $V = \bigoplus_{\omega \in \X} V_\omega$ and $\dim V_\omega < \infty$ for all
  $\omega \in \X$,
  \item the action of $e_i$ and $f_i$ is locally nilpotent for all $i \in \scrI$,
  \item there exists a finite set $F \subset \X$ such that 
  $\mathrm{wt}(V) \subset F + \sum_{i \in \scrI} \bbZ_{\leqslant 0} \alpha_i$.
\end{itemize}
Let $\X^+ = \{\omega \in \X\, \mid\, \langle \alpha_i^\vee,\omega\rangle \in \bbN \text{ for all } i \in \scrI\}$ be the set of \emph{integral dominant weights}. 
Given $\Lambda \in  \X^+$, we denote
by $\bfL(\Lambda)$ the simple integrable highest weight module with highest weight $\Lambda$.

\smallskip

\subsubsection{Quantized enveloping algebras}

Let $u$ be a formal variable and $A=\bbC[u,u^{-1}]$.
Let $U_u(\frakg)$ be the quantized enveloping algebra over $\bbC(u)$.
Let $U_A(\frakg)\subset U_u(\frakg)$ be
Lusztig's divided power version of $U_u(\frakg)$.
For each integral weight $\Lambda$ the module $\bfL(\Lambda)$ admits a deformed version $\bfL_u(\Lambda)$ over $U_u(\frakg)$
and an integral form $\bfL_A(\Lambda)$ which is the $U_A(\frakg)$-submodule of
$\bfL_u(\Lambda)$ generated by the highest weight vector $|\Lambda\rangle$.
Let $\calO^\text{int}_u$ be the category consisting of the $U_u(\frakg)$-modules which are (possibly infinite) direct sums of $\bfL_u(\Lambda)$'s.
If $V_u\in\calO^\text{int}_u,$ then its integral form $V_A$ is the corresponding sum of the modules $\bfL_A(\Lambda)$.
It depends on the choice of a family of highest weight vectors of the constituents of $V_u$.

\medskip

\subsection{Categorical representations on abelian categories}
\label{sec:categoricalactions}\hfill\\

In this section we recall from \cite{CR,R08} the notion of a categorical action of $\frakg$.
It consists of the data of functors $E_i$, $F_i$ lifting the Chevalley
generators $e_i$, $f_i$ of $\mathfrak{g}$, together with an action of an affine
Hecke algebra on $(\bigoplus_{i\in I} F_i)^{m}$.

\smallskip

\subsubsection{Affine Hecke algebras and representation data}\label{subsec:hecke}
Let $\scrC$ be an abelian $R$-category and $v \in R^\times$.

\begin{definition} 
A \emph{representation datum} on $\scrC$ with parameter $v$ is a tuple $(E,F,X,T)$ where
$E$, $F$ are bi-adjoint functors $\scrC\to\scrC$ and $X\in\End(F)^\times$, $T\in\End(F^2)$
are endomorphisms of functors satisfying the following conditions: 
\begin{itemize}[leftmargin=8mm]
  \item[(a)] $1_FT\circ T1_F\circ 1_FT=T1_F\circ 1_FT\circ T1_F$,

  \item[(b)] $(T+1_{F^2})\circ(T-v1_{F^2})=0$,

  \item[(c)] $T\circ(1_FX)\circ T=vX1_F$.
\end{itemize}
\end{definition}

This definition can also be formulated in terms of actions of affine Hecke algebras.
For $m\geqs 1,$ let $\bfH_{R,m}^v$ be the \emph{affine Hecke algebra} of type $A_{m-1}$ over $R$.
It is
generated by $T_1,\ldots,T_{m-1}$, $X^{\pm 1}_1,\ldots,X^{\pm 1}_m$
subject to the well-known relations.
We will also set $\bfH_{R,0}^v=R$. 

\smallskip

Given $(E,F)$ a pair of biadjoint functors, and $X\in\End(F)$, $T\in\End(F^2)$,
the tuple $(E,F,X,T)$ is a representation datum if and only if 
for each $m\in\bbN$, the map
$$ \begin{array}{rcl} \phi_{F^m}\, :\, \bfH_{R,m}^v & \longrightarrow & \End(F^m) \\[4pt]
  X_k & \longmapsto &  1_{F^{m-k}} X 1_{F^{k-1}}\\
  T_l & \longmapsto & 1_{F^{m-l-1}}T1_{F^{l-1}}\\
  \end{array}$$
is a well-defined $R$-algebra homomorphism.

\smallskip

\subsubsection{Categorical representations}\label{subsec:categoricalaction}
We assume now that $R$ is a field and that
$\scrC$ is Hom-finite. We fix a pair $(\scrI,v)$ as in \S \ref{sec:quivers} and
we denote by $\frakg = \frakg_\scrI$ the Lie algebra associated to that pair.

\begin{definition}[\cite{R08}]\label{df:cat1}
A \emph{$\frakg$-representation} on $\scrC$ consists of a representation datum
$(E,F,X,T)$ on $\scrC$ and of a decomposition $\scrC=\bigoplus_{\omega\in\X}\scrC_\omega$.
For each $i\in \scrI,$ let $F_i,$ $E_i$ be the generalized $i$-eigenspaces
of $X$ acting on $F,$ $E$ respectively. We assume in addition that
\begin{itemize}[leftmargin=8mm]

\item[(a)] $F=\bigoplus_{i\in \scrI} F_i$ and $E=\bigoplus_{i\in \scrI} E_i$,

\item[(b)] the action of $[E_i]$ and $[F_i]$ for $i\in \scrI$ endow $[\scrC]$ with 
a structure of integrable $\frakg$-module such that $[\scrC]_\omega=[\scrC_\omega]$,

\item[(c)] $E_i(\scrC_\omega)\subset\scrC_{\omega+\alpha_i}$ and
$F_i(\scrC_\omega)\subset\scrC_{\omega-\alpha_i}$.
\end{itemize}
\end{definition}

We say that the tuple $(E, F, X, T)$ and the decomposition $\scrC=\bigoplus_{\omega\in \X}
\scrC_{\omega}$ is a \emph{$\frakg$-categorification} of the integrable $\frakg$-module
$[\scrC]$.

\medskip

\subsection{Minimal categorical representations} \label{subsec:minimalcat}\hfill\\

For most of the results in the rest of  Section \ref{part:categoricalactions} we will assume that $R$ is a field
and that $\scrI$ is finite. In particular $v \in R^\times$ will be a root of unity.

\smallskip

Let $m \geqslant 0$, $v \in R^\times$ and $\bfH_{R,m}^v$ be the affine Hecke algebra as
defined in \S \ref{subsec:hecke}. We fix a tuple $Q=(Q_1,\ldots,Q_l)$ in $(R^\times)^l$. The
\emph{cyclotomic Hecke algebra} $\bfH^{Q;\,v}_{R,m}$ is the quotient of $\bfH_{R,m}^v$ 
by the two-sided ideal generated by $\prod_{i=1}^l(X_1-Q_{i})$.

\smallskip

Assume now that $R$ is a field. Any finite dimensional $\bfH^{Q;\,v}_{R,m}$-module $M$
is the direct sum of the weight subspaces
  $$M_\nu=\{v\in M\,\mid\,(X_r-i_r)^dv=0,\
  r\in[1,m],\, d\gg 0\},\quad\nu=(i_1,\dots,i_m)\in R^m.$$
Decomposing the regular module, we get a system of orthogonal idempotents $\{e_\nu\,;\,\nu
\in R^m\}$ in $\bfH^{Q;\,v}_{R,m}$ such that $e_\nu M=M_\nu$ for each $M$. 
The eigenvalues of $X_r$ are always of the form $Q_i v^j$ for some $i \in \{1,\ldots,l\}$
and $j \in \bbZ$. As a consequence, if we set $\scrI = \bigcup Q_i v^\bbZ$, then
$e_\nu = 0$ unless $\nu \in \scrI$. The pair $(\scrI, v)$ satisfies the assumptions of
\S \ref{sec:quivers} and we can consider a corresponding Kac-Moody algebra $\frakg_\scrI$ and
its root lattice $\Q_\scrI$. Given $\alpha\in \Q_\scrI^+$ of height $m$, let
$e_\alpha=\sum_{\nu} e_\nu$ where the sum runs over the set of all tuples such that
$\sum_{r=1}^m \alpha_{i_r}=\alpha$. The nonzero $e_\alpha$'s are the primitive central
idempotents in $\bfH^{Q;\,v}_{R,m}$. 

\smallskip

To the dominant weight $\Lambda_Q=\sum_{i=1}^l\Lambda_{Q_{i}}$ of $\frakg_\scrI$ and to
any $\alpha \in \Q_\scrI^+$ we associate the following abelian categories:
  $$\scrL(\Lambda_Q) = \bigoplus_{m\in\bbN}\bfH^{Q;\, v}_{R,m}\mod \quad \text{and} \quad
  \scrL(\Lambda_Q)_{\Lambda_Q-\alpha}=e_\alpha\bfH^{Q;\,v}_{R,m}\mod.$$
For any $m<n$, the $R$-algebra embedding of the affine Hecke algebras $\bfH_{R,m}^v
\hookrightarrow \bfH_{R,n}^v$ given by $T_i\mapsto T_i$ and $X_j\mapsto X_j$ 
induces an embedding $\bfH^{Q;\, v}_{R,m}\hookrightarrow \bfH^{Q;\,v}_{R,n}$.
The $R$-algebra $\bfH^{Q;\,v}_{R,n}$ is free as a left and as a right
$\bfH^{Q;\,v}_{R,m}$-module. This yields a pair of exact adjoint functors 
$(\Ind^{n}_m,\,\Res^{n}_m)$ between $\bfH^{Q;\,v}_{R,n}\mod$ and $\bfH^{Q;\,v}_{R,m}\mod$.
They induce endofunctors $E$ and $F$ of $\scrL(\Lambda_Q)$ by $E=\bigoplus_{m\in\bbN}
\Res_m^{m+1}$ and $F=\bigoplus_{m\in\bbN}\Ind_m^{m+1}$. The right multiplication on 
$\bfH^{Q;\, v}_{R,m+1}$ by $X_{m+1}$ yields an endomorphism of the functor $\Ind_m^{m+1}$.
The right multiplication by $T_{m+1}$ yields an endomorphism of $\Ind_m^{m+2}.$
We define $X\in\End(F)$ and $T\in\End(F^2)$ by $X=\bigoplus_{m}X_{m+1}$ and 
$T=\bigoplus_{m}T_{m+1}$.

\smallskip

This construction yields a categorification of the simple
highest module $\bfL(\Lambda_Q)$ of $\frakg_\scrI$. 
Indeed, a theorem of Kang and Kashiwara implies that this holds in the more
general setting of cyclotomic quiver Hecke algebras of arbitrary type. 

\begin{theorem}[\cite{KK12}, \cite{K12}]\label{thm:minimalcat}  \hfill  
  \begin{itemize}[leftmargin=8mm]
  \item[$\mathrm{(a)}$] The endofunctors $E$ and $F$ of $\scrL(\Lambda_Q)$ are biadjoint. 
  \item[$\mathrm{(b)}$] The tuple $(E,F,X,T)$ and the decomposition $\scrL(\Lambda_Q)=
  \bigoplus_{\omega\in \X} \scrL(\Lambda_Q)_{\,\omega}$ is a $\frakg_\scrI$-categorification
  of $\bfL(\Lambda_Q)$.
  \qed
  \end{itemize}
\end{theorem}

This categorical representation is called the \emph{minimal categorical
$\frakg_\scrI$-representation} of highest weight $\Lambda_Q$. 

\smallskip

The $\frakg_\scrI$-modules we are interested in are direct sums of various irreducible
highest weight modules $\bfL(\Lambda_Q)$.
Let $(\scrI,v)$ as in \S \ref{sec:quivers},
and $\frakg = \frakg_\scrI$ be a corresponding Kac-Moody algebra. 
Let $(E,F,X,T)$ be a $\frakg$-representation on an abelian $R$-category $\scrC$.
We want to relate $\scrC$ to minimal categorical $\frakg_\scrI$-representations.
To do that, recall that for any $m \geqslant 0$ we have an $R$-algebra homomorphism 
$\phi_{F^m} \, : \, \bfH^{Q;\, v}_{R,m}\ \longrightarrow \End(F^m)^\op$. Given an
object $M$ in $\scrC$, it specializes to an $R$-algebra homomorphism 
  $$ \bfH^{Q;\, v}_{R,m} \ \longrightarrow \End(F^m M)^\op =: \scrH(F^mM).$$

\begin{proposition}[\cite{R08}]\label{prop:unicite}
Assume that the simple
roots are linearly independent in $\X$. 
Let $(E,F,X,T)$ be a  representation of $\frakg$ in a abelian $R$-category $\scrC$,
and $M\in\scrC_{\omega}$. Assume that $EM=0$ and $\End_\scrC(M)=R$. Then
  \begin{itemize}[leftmargin=8mm]
    \item[$\mathrm{(a)}$] $\omega \in \X^+$ is an integral dominant weight,    
    \item[$\mathrm{(b)}$] if we write $\Lambda_Q = \sum_{i \in \scrI} \langle 
    \alpha_i^\vee,\omega\rangle \Lambda_i = 
    \sum_{p=1}^l \Lambda_{Q_p}$ for some $Q=(Q_1,\ldots,Q_l) \in \scrI^l$ and $l \geqslant 1$,
    then for all $m \geq 0$ the map $\phi_{F^m}$ factors to an $R$-algebra isomorphism 
    $$\bfH_{R,m}^{Q;\,v} \ \mathop{\longrightarrow}\limits^\sim \ \scrH(F^mM).$$
    \qed
  \end{itemize}
\end{proposition}

\medskip

\subsection{Crystals}\hfill\\

We start by a review of Kashiwara's theory of perfect bases and crystals.
We will be working with the Kac-Moody algebra $\frakg$ coming from a pair $(\scrI,v)$ as
in \S \ref{sec:quivers}. 

\begin{definition}
An {\em abstract crystal} is a set $B$ together with maps
$\mathrm{wt} : B \to \P$, $\varepsilon_i,$ $\varphi_i : B \to {\bbZ}
\sqcup \{-\infty \}$ and $\widetilde{e_i},\widetilde{f_i}$ : $B \to
B\sqcup\{0\}$ for all $i \in \scrI$ satisfying the following properties:
 \begin{itemize}[leftmargin=8mm]
  \item[(a)] $\varphi_i(b)=\varepsilon_i(b)+\left<\alpha^\vee_i,\mathrm{wt}(b)\right>,$
    
  \item[(b)] $\mathrm{wt}(\widetilde{e_i}b)=\mathrm{wt}(b)+\alpha_i$ and $\mathrm{wt}
  (\widetilde{f_i}b)=\mathrm{wt}(b)-\alpha_i,$
  
  \item[(c)] $b=\widetilde{e_i}b'$ if and only if $ \widetilde{f_i}b=b'$, where
  $b, b'  \in B$, $i \in \scrI$,

  \item[(d)] if $\varphi_i(b)=-\infty$, then $\widetilde{e_i}b=\widetilde{f_i}b=0$,
  
  \item[(e)] if $b\in B$ and $\widetilde{e_i}b \in B$,  then
  $\varepsilon_i(\widetilde{e_i}b)= \varepsilon_i(b)-1$ and $\varphi_i(\widetilde{e_i}b)=
  \varphi_i(b)+1$,

  \item[(f)] if $b\in B$ and  $\widetilde{f_i}b \in B$, then 
  $\varepsilon_i(\widetilde{f_i}b)= \varepsilon_i(b)+1$ and $\varphi_i(\widetilde{f_i}b)=
  \varphi_i(b)-1$.
\end{itemize}
\end{definition}

Note that by (a), the map $\varphi_i$ is entirely determined by $\varepsilon_i$ and
$\mathrm{wt}$. We may therefore omit $\varphi_i$ in the data of an abstract crystal
and denote it by $(B,\widetilde e_i,\widetilde f_i, \varepsilon_i, \mathrm{wt})$.   

\smallskip

An isomorphism between crystals $B_1$, $B_2$ is a bijection ${\psi}:B_1\sqcup\{0\} \longrightarrow B_2\sqcup \{0\}$ such that $\psi(0)=0$ which commutes with $\mathrm{wt}$,
$\varepsilon_i$, $\varphi_i$ $\widetilde{f}_i,$, $\widetilde{e}_i$.

\smallskip

Let $V_u$ be an integrable $U_u(\frakg)$-module in $\calO^\text{int}_u$.
Let $V_A$ be an integral form of $V_u$.
A \emph{lower crystal lattice} in $V_u$ is a free $\bbC[u]$-submodule $\calL$ of $V_A$ such that
$V_A=A\calL$, $\calL=\bigoplus_{\lambda\in \X}\calL_\lambda$ with $\calL_\lambda=\calL\cap(V_A)_\lambda$ 
and $\calL$ is preserved by the \emph{lower Kashiwara crystal operators}  $\widetilde e_i^\lw$, $\widetilde f_i^\lw$ on $V_u$.
A \emph{lower crystal basis} of $V_u$ is a pair $(\calL,B)$ where $\calL$ is a lower crystal lattice of $V_u$ and $B$ is a basis of $\calL/u\calL$ 
such that we have $B=\bigsqcup_{\lambda\in\X} B_\lambda$ where $B_\lambda=B\cap(\calL_\lambda/u\,\calL_\lambda)$,
$\widetilde e_i^\lw(B),\,\widetilde f_i^\lw(B)\subset B\sqcup\{0\}$ and $b'=\widetilde f_i^\lw b$ if and only if $b=\widetilde e_i^\lw b'$ for each
$b,b'\in B$.
A \emph{lower global basis} of $V_u$ (or a \emph{canonical basis}) is an $A$-basis $\bfB$ of $V_A$ such that the lattice 
$\calL=\bigoplus_{b\in\bfB}\bbC[u]\,b$ and the basis $B=\{b\,\text{mod}\,u\,\calL\,|\,b\in\bfB\}$ of $\calL/u\,\calL$ form a lower crystal basis.
One defines in a similar way an \emph{upper crystal lattice}, an \emph{upper crystal basis} and an 
\emph{upper global basis} (or a \emph{dual canonical basis}) using the
\emph{upper Kashiwara crystal operators} $\widetilde e_i^\up$, $\widetilde f_i^\up$ on $V_u$, see, e.g., \cite[def.~4.1,4.2]{KOP}.
Any $U_u(\frakg)$-module in $\calO^\text{int}_u$ admits a lower crystal,  an upper crystal and a global basis.

\smallskip

If $(\calL,B)$, $(\calL^\vee,B^\vee)$ are lower, upper crystal bases, then $(B,\widetilde e_i^\lw,\widetilde f_i^\lw)$, $(B^\vee,\widetilde e_i^\up,\widetilde f_i^\up)$
are abstract crystals. 

\smallskip

Let $E_i,$ $F_i,$ $u^h$ with $i\in I$, $h\in\X^\vee,$ be the standard generators of $U_u(\frakg)$.
There exists a unique non-degenerate symmetric bilinear form $(\bullet,\bullet)$ on the module $\bfL_u(\Lambda)$ with highest 
weight vector $|\Lambda\rangle$ satisfying
\begin{itemize}[leftmargin=8mm]
\item $(|\Lambda\rangle,|\Lambda\rangle)=1$,
\item $(E_ix,y)=(x,F_iy)$, $(F_ix,y)=(x,E_iy)$, $(u^{h}x,y)=(x,u^{h}y)$, 
\item $(\bfL_u(\Lambda)_\lambda,\bfL_u(\Lambda)_\mu)=0$ if $\lambda\neq\mu$.
\end{itemize}
If $(\calL, B)$ is a lower crystal basis of $\bfL_u(\Lambda)$ then the pair $(\calL^\vee,B^\vee)$ such that
$\calL^\vee=\{x\in\bfL_u(\Lambda)\,|\,(x,\calL)\subset \bbC[u]\}$ and 
$B^\vee$ is the basis of $\calL^\vee/u\,\calL^\vee$ which is dual to $B$ with respect to the non-degenerate bilinear form
$\calL^\vee/u\,\calL^\vee\times\calL/u\,\calL\to\bbC$ induced by $(\bullet,\bullet),$ is an upper crystal basis.
Further, taking a basis element in $B$ to the dual basis element in $B^\vee$ is a crystal isomorphism 
$(B,\widetilde e_i^\lw,\widetilde f_i^\lw)\to(B^\vee,\widetilde e_i^\up,\widetilde f_i^\up)$.
Therefore, if $\bfB$ is a lower global basis of $\bfL_u(\Lambda)$ then the dual basis $\bfB^\vee$ with respect to the non-degenerate bilinear form
$(\bullet,\bullet)$ is an upper global basis and the corresponding abstract crystals 
$(B,\widetilde e_i^\lw,\widetilde f_i^\lw)$ and $(B^\vee,\widetilde e_i^\up,\widetilde f_i^\up)$ are canonically isomorphic.

\medskip

\subsection{Perfect bases}\hfill\\

The crystals that we will consider in this paper all come from particular bases of $\frakg_\scrI$-modules called \emph{perfect bases}. 
Let us define them.
Let $V \in \calO^\mathrm{int}$ be an integrable highest weight $\frakg_\scrI$-module. 
For each $i \in \scrI$ and $x \in V$ we define
  $$\ell_i(x)= \max\{k\in \bbN\,\mid\, e_i^{k}x\not=0\} = 
  \min\{k\in \bbN\,\mid\, e_i^{k+1}x=0\}$$
with the convention that $\ell_i(0) = -\infty$.
For each integer $k$, we also consider the vector spaces
  $$ V_i^{\leqslant\, k}=\{x\in V\,\mid\,\ell_i(x)\leqslant k\}, \quad 
   V^{\leqslant\, k}=\bigcap_{i\in I}V_i^{\leqslant\, k},\quad V_i^{k}=V_i^{\leqslant\, k}/\,V_i^{<\, k}.
  $$

\begin{definition} \label{def:perfect}
A basis $B$ of $V$ is {\em perfect} if
 \begin{itemize}[leftmargin=8mm]
  \item[(a)] $B=\bigsqcup_{\mu \in \X_\scrI} B_\mu$ where $B_\mu= B \cap V_\mu,$
  \item[(b)]  for any $i\in \scrI$, there is a map $\bfe_i : B\to B\sqcup\{0\}$
  such that for any $b \in B$, we have
  \begin{itemize}[leftmargin=8mm]
    \item[(i)] if $\ell_i(b)=0$, then $\bfe_ib=0$,
    \item[(ii)] if $\ell_i(b)>0$, then $\bfe_i b\in B$ and
    $e_ib\in\bbC^\times\,\bfe_ib+V_i^{<\,\ell_i(b)-1}$,
  \end{itemize}
  \item[(c)] if $\bfe_i b= \bfe_i b' \neq 0$ for $b, b' \in B$, then $b=b'.$
\end{itemize}
\end{definition}

\smallskip

Any $\frakg$-module in $\calO^{\text{int}}$ admits a perfect basis.
More precisely, we have the following.

\begin{proposition}\label{prop:upper-perfect}
 If $V$ is an integrable $\frakg$-module in $\calO^\text{int}$ with a quantum deformation $V_u$, then
the specialization at $u=1$ of an upper global basis of $V_u$ is a perfect basis of $V$. 
\qed
\end{proposition}

\smallskip

To any categorical representation we associate a perfect basis as in
\cite[prop. 6.2]{S}. More precisely, let $R$ be a field (of any characteristic) and
consider a $\frakg$-representation on an abelian artinian $R$-category $\scrC$.
Then, for each $i\in\scrI$ we define the maps
\begin{align*}
&\widetilde E_i\,:\,\Irr(\scrC)\to\Irr(\scrC)\sqcup\{0\},\quad [L]\mapsto [\soc(E_i(L))],\\
&\widetilde F_i\,:\,\Irr(\scrC)\to\Irr(\scrC)\sqcup\{0\},\quad [L]\mapsto [\head(F_i(L))].
\end{align*}

\begin{proposition}\label{prop:PBfromcategorification}
The tuple $\big(\Irr(\scrC),\widetilde E_i,\widetilde F_i\big)$
defines a perfect basis of $[\scrC]$.
\qed
\end{proposition}

We now recall how to construct an abstract crystal from a perfect basis $B$. We 
set $\widetilde e_i=\bfe_i$. For all $b\in B$ we set $\widetilde{f_i}b=b'$ if $\bfe_ib'=b$
for some $b'\in B$, and $0$ otherwise. Then it follows easily from the definition that
$(B,\widetilde{e_i}, \widetilde{f_i},\ell_i, \mathrm{wt})$ is an abstract crystal.

\smallskip

We finish this section with a result which will be important to identify the crystal 
obtained by the categorification with the crystal of some Fock space.
For each $i\in \scrI$ and $k\in\bbN$, we set
$B^{\leqslant\,k}=V^{\leqslant k}\cap B$ and $B^{\leqslant\,k}_i=V^{\leqslant k}_i\cap B$.

\begin{proposition} \label{prop:perfect-bases}
Let $B$ and $B'$ be perfect bases of $V \in \calO^{\mathrm{int}}$. 
\begin{itemize}[leftmargin=8mm]
  \item[$\mathrm{(a)}$]  $B^{\leqslant\,k}$ and $B^{\leqslant\,k}_i$ are bases of $V^{\leqslant\,k}$
  and $V^{\leqslant\,k}_i.$ 
  \item[$\mathrm{(b)}$] 
Assume that there is a bijection 
$\varphi\,:\,B\to B'$ and a partial order $\leqslant$ on $B$ such that
$\varphi( b)\in b + \sum_{c > b}\bbC\,c$ for each $b\in B$. Then the map
$\varphi$ is a crystal isomorphism $B\mathop{\to}\limits^\sim B'$.
\qed
\end{itemize}
\end{proposition}

\medskip

\subsection{Derived equivalences} \hfill\\

Given $V$ an integrable $\frakg$-module and $i \in \scrI$, one can consider
the action of the simple reflection $s_i = \exp(-f_i) \exp(e_i) \exp(-f_i)$
on $V$. 
For each weight $\omega\in \X$, this action maps a weight space $V_\omega$ to $V_{s_i(\omega)}$ with
$s_i(\omega) = \omega - \langle \alpha_i^\vee,\omega \rangle \alpha_i$. 
If $\scrC$ is a categorification of $V$, then it restricts to an $\fraks\frakl_2(\bbC)$-categorification in the sense 
of Chuang-Rouquier. In particular, the simple objects are weight vectors for the categorical $\fraks\frakl_2(\bbC)$-action.
Thus, the theory of Chuang-Rouquier can be applied and
\cite[thm. 6.6]{CR} implies that $s_i$ 
can be lifted to a derived equivalence $\Theta_i$ of $\scrC$.

\begin{theorem}\label{thm:reflection}
Assume that $R$ is a field. 
Let $(E,F,X,T)$
be a representation of $\frakg$ in an abelian $R$-category $\scrC$, and
$i \in \scrI$. Then there exists a derived self-equivalence $\Theta_i$ of
$\scrC$ which restricts to derived equivalences
 $ \Theta_i \, : \, D^b(\scrC_\omega) \mathop{\longrightarrow}\limits^\sim
 D^b(\scrC_{s_i(\omega)})$
for all weightn $\omega \in \X$. Furthermore, $[\Theta_i] = s_i$ as a linear map
of $[\scrC]$.
\qed
\end{theorem}

\medskip

\section{Representations on Fock spaces}\label{part:fock}

Let $R$ be a noetherian commutative domain with unit. As in \S\ref{sec:quivers}, we fix an element $v \in R^\times$
and a subset $\scrI$ of $R^\times$ which is stable by multiplication by $v$ and $v^{-1}$.
We explained in \S\ref{subsec:quivers} how one can associate a Lie algebra $\frakg=\frakg_\scrI$ to this data. 
In this section we recall the construction of (charged) Fock spaces which are particular
integrable representations of $\frakg$. 

\medskip

\subsection{Combinatorics of $l$-partitions}\label{sec:combinatorics}

\subsubsection{Partitions and $l$-partitions}\label{subsec:partitions}
A \emph{partition} of $n$ is a non-increasing sequence of non-negative integers
$\lambda = (\lambda_1 \geqslant \lambda_2 \geqslant \cdots)$ whose terms add up to $n$. We denote by 
$\scrP_n$ the set of partitions of $n$ and by $\scrP=\bigsqcup_n\scrP_n$ be the set
of all partitions. Given a partition $\lambda$, we write $|\lambda|$ for the \emph{weight}
of $\lambda$.
We associate to $\lambda=(\lambda_1,\lambda_2,\dots)$ the \emph{Young diagram} $Y(\lambda)$
defined by $Y(\lambda)=\{(x,y)\in\bbZ_{>0}\times\bbZ_{>0}\,\mid\,y\leqslant \lambda_x\}.$
It may be visualised by an array of boxes in left justified rows with $\lambda_x$ boxes in
the $x$-th row. If $\lambda$, $\mu$ are partitions of $n$ then we write $\lambda\geqslant\mu$
if for all $n\geqslant i\geqslant 1$ we have $\sum_{j=1}^i\lambda_j\geqslant\sum_{j=1}^i\mu_j$.
This relation defines a partial order on $\scrP$ called the \emph{dominance order}.
Let $\lambda^*$ denote the partition dual (or conjugate) to $\lambda$.

\smallskip

An \emph{$l$-partition} of $n$ is an $l$-tuple of partitions whose weights add up to $n$.
We denote by $\scrP^l_n$ be the set of $l$-partitions of $n$ and by $\scrP^l=\bigsqcup_n\scrP^l_n$
the set of all $l$-partitions. The Young diagram of the $l$-partition $\lambda=(\lambda^{1},\ldots,\lambda^{l})$ is the set $Y(\lambda)=\bigsqcup_{p=1}^lY(\lambda^p)\times\{p\}.$
Its weight is the integer $|\lambda|=\sum_p|\lambda^{p}|.$

\smallskip 

\subsubsection{Residues and content}\label{subsec:contents}
We fix $Q = (Q_1,\ldots,Q_l) \in \scrI^l$.
Let $\lambda$ be an $l$-partition and $A=(x,y,p)$ be a node in $Y(\lambda)$. 
The \emph{$(Q,v)$-shifted residue} of the node $A$ is the element of $\scrI$ given by
$\res(A,Q)_\scrI=v^{y-x}Q_{p}.$  Let $n_i(\lambda,Q)_\scrI$ be the number of nodes of
$(Q,v)$-shifted residue $i$ in $Y(\lambda)$. If $\lambda$, $\mu$ are $l$-partitions such
that $|\mu|=|\lambda|+1$ we write $\res(\mu-\lambda,Q)_\scrI=i$ if $Y(\mu)$ is obtained by
adding a node of $(Q,v)$-shifted residue $i$ to $Y(\lambda)$.
A \emph{charge} of the tuple $Q = (Q_1,\ldots,Q_l)$ is an $l$-tuple of integers
$s = (s_1,\ldots,s_l)$ such that $Q_p = v^{s_p}$ for all $p = 1,\ldots,l$.
Conversely, given $\scrI \subset R^\times$ and $v \in R^\times$ as in \S\ref{sec:quivers}, any
$\ell$-tuple of integers $s \in \bbZ^l$ defines a tuple $Q = (v^{s_1},\ldots,v^{s_l})$ with
charge $s$. 
The \emph{$s$-shifted content} of the box $A=(x,y,p)$ is the integer $\text{ct}^s(A)=s_p+y-x$.
It is related to the residue of $A$ by the formula $\res(A,Q)_\scrI=v^{\text{ct}^s(A)}$.
We will also write $p(A)=p$.
We will call \emph{charged $l$-partition} a pair $(\mu,s)$ in $\scrP^l\times\bbZ^l$.

\smallskip

\subsubsection{$l$-cores and $l$-quotients}\label{subsec:l-core}
We start with the case $l=1$. The set of \emph{$\beta$-numbers} of a charged partition
$(\lambda,d) \in \scrP\times\bbZ$ is the set given by $\beta_d(\lambda)=\{\lambda_u+d+1-u\,\mid
\,u\geqslant 1\}.$ 
The charged partition $(\lambda,d)$ is uniquely determined by the set $\beta_d(\lambda)$.

\smallskip

For any positive integer $e $, an \emph{$e$-hook} of $(\lambda,d)$ is a pair $(x,x+e)$ such
that $x+e\in \beta_d(\lambda)$ and $x\not\in \beta_d(\lambda)$. Removing the $e$-hook
$(x,x+e)$ corresponds to replacing $x+e$ with $x$ in $\beta_d(\lambda)$. We say
that the charged partition $(\lambda,d)$ is an \emph{$e$-core} if it does not have any
$e$-hook. This does not depend on $d$. 

\smallskip

Next, we construct a bijection $\tau_l:\scrP\times\bbZ\to\scrP^l\times\bbZ^l$.
It takes the pair $(\lambda,d)$ to $(\mu,s),$ where $\mu=(\mu^1,\dots,\mu^l)$ is an
$l$-partition and $s=(s_1,\dots,s_l)$ is a $l$-tuple in
$\bbZ^l(d)=\{s\in\bbZ^l\,\mid\,s_1+\cdots+s_l=d\}.$
The bijection is uniquely determined by the relation
$  \beta_d(\lambda)=\bigsqcup_{p=1}^l\big(p-l+l\beta_{s_p}(\mu^p)\big).$

\smallskip

The bijection $\tau_l$ takes the pair $(\lambda,0)$ to $(\lambda^{[l]},\lambda_{[l]}),$ 
where $\lambda^{[l]}$ is the \emph{$l$-quotient} of $\lambda$ and $\lambda_{[l]}$ lies
in $\bbZ^l(0)$. Since $\lambda$ is an $l$-core if and only if $\lambda^{[l]}=\emptyset$,
this bijection identifies the set of $l$-cores and $\bbZ^l(0)$. We define the \emph{$l$-weight} $w_l(\lambda) := |\lambda^{[l]}|$ of the partition $\lambda$ to be the weight of its $l$-quotient. 

\smallskip

We will mostly consider the bijection $\tau_l$ for $l=2$.
In particular, a 2-core is either $\Delta_0=\emptyset$ or a triangular partition $\Delta_t=(t,t-1,\dots,1)$ with $t\in\bbN$.
We abbreviate $\sigma_t=(\Delta_t)_{[2]}$, and we write $\sigma_t=(\sigma_1,\sigma_2)$. We have
\begin{align}\label{sigma}
\sigma_t=\begin{cases}\big(-t/2,\,t/2\big)&\ \text{if}\ t\ \text{is\ even},\\
\big((1+t)/2,-(1+t)/2\big)&\ \text{if}\ t\ \text{is\ odd}.
\end{cases}
\end{align}
For each bipartition $\mu$, let $\varpi_t(\mu)$ denote the unique partition 
with 2-quotient $\mu$ and 2-core $\Delta_t.$ 
Thus, the bijection $\tau_2$ maps $(\varpi_t(\mu),0)$ to the pair $(\mu,\sigma_t)$.

\medskip

\subsection{Fock spaces}
\label{subsec:Fock}\hfill\\

For a reference for the results presented in this section, see for example \cite{U}, \cite{Y05}.
Let $Q = (Q_1,\ldots,Q_l) \in \scrI^l$. It defines an integral dominant weight $\Lambda_Q=\sum_{p=1}^l\Lambda_{Q_{p}} \in \P^+$. The \emph{Fock space} $\bfF(Q)_{\scrI}$ is the $\bbC$-vector
space with basis $\{|\lambda,Q\rangle_{\scrI}\,|\,\lambda\in\scrP^l\}$ called  the 
\emph{standard monomial basis}, and action of $e_i, f_i$ for all $i \in \scrI$ given by
\begin{equation}\label{eq:EF}
  f_i(|\lambda,Q\rangle_{\scrI})=\sum_{\mu}|\mu,Q\rangle_{\scrI},\qquad
  e_i(|\mu,Q\rangle_{\scrI})=\sum_{\lambda}|\lambda,Q\rangle_{\scrI},
\end{equation}
where the sums run over all partitions such that $\res(\mu-\lambda,Q)_\scrI=i$.
This endows $\bfF(Q)_{\scrI}$ with a structure of $\frakg'$-module.
The Fock space $\bfF(Q)_{\scrI}$ can also be equipped with a symmetric non-degenerate bilinear
form $\langle\bullet,\bullet\rangle_{\scrI}$ for which the standard monomial basis is orthonormal.
To avoid cumbersome notation, we shall omit the subscript $\scrI$ when not necessary.
It is easy to see that every element of the standard monomial basis is a weight vector.

\smallskip

 The $\frakg'$-submodule of $\bfF(Q)$ generated by $|\emptyset,Q\rangle$ is
 isomorphic to $\bfL(\Lambda_Q)$. Furthermore, if $\scrI= A_\infty$, then
 $\bfF(Q)=\bfL(\Lambda_Q)$. 
Using the minimal categorification $\scrL(\Lambda_Q)$ of $\bfL(\Lambda_Q),$ the map $\bfL(\Lambda_Q)\to\bfF(Q)$
can be made more explicit. 

\smallskip

To explain this, let us first recall briefly the definition of the \emph{Specht modules}.
Assume that $R$ has characteristic 0 and contains a primitive $l$-th root $\zeta$ of 1, so
$R$ is a splitting field of the complex reflection group $G(l,1,m)$.
Let $\Irr(R\frakS_m)=\{\phi_\lambda\,|\,\lambda\in\scrP_m\}$
be the standard labelling of the characters of the symmetric group.
Then $$\Irr(RG(l,1,m))=\{\scrX_{\lambda}\,|\,\lambda\in\scrP^l_m\}$$
is the labelling of the simple modules such that $\scrX_\lambda$ is induced from the 
$G(l,1,|\lambda_1|)\times\ldots\times G(l,1,|\lambda_l|)$-module
$$\phi_{\lambda^{(1)}}\chi^0\otimes
\phi_{\lambda^{(2)}}\chi^1\otimes\cdots\otimes\phi_{\lambda^{(l)}}\chi^{l-1}.$$
Here, we denote by $\chi^p$ the one-dimensional module of the $|\lambda_p|$-th cartesian power of 
the cyclic group $G(l,1,1)$ given by the $p$-th 
power of the determinant, see, e.g., \cite[sec.~5.1.3]{GJ}.
Recall that (for every field $R$) the $R$-algebra $\bfH^{Q;\,v}_{R,m}$ is split and
that it is semisimple if and only if we have, see, e.g., \cite[sec.~3.2]{Ma},
\begin{equation*}\label{(A)}
\prod_{i=1}^m(1+v+\cdots+v^{i-1})\,\prod_{a<b}\prod_{-m<r<m}(v^r\,Q_{a}-Q_{b})\neq 0.
\end{equation*}
Thus, by Tits' deformation theorem, under the specialization $v\mapsto 1$ and $Q_p\mapsto\zeta^{p-1}$, the labelling of $\Irr(R G(l,1,m))$ yields a 
canonical labelling  
$$\Irr(\bfH^{Q;\,v}_{R,m})=\{S(\lambda)_R^{Q;\,v}\,|\,\lambda\in\scrP^l_m\}.$$
Now, if $R$ is a commutative domain with fraction field $K$ of characteristic 0 as above, we define 
the $\bfH^{Q;\,v}_{R,m}$-module $S(\lambda)^{Q,\,v}_R$ as in \cite[sec.~2.4.3]{RSVV} or \cite[sec.~5.3]{GJ},
using $S(\lambda)^{Q;\,v}_K$ and the dominance order on $\scrP^l_m$,
and if $\theta\,:\,R\to\k$ 
is a ring homomorphism such that $\k$ is the fraction field of $\theta(R)$
we set $S(\lambda)^{Q;\,v}_\k=\k S(\lambda)^{Q,\,v}_R$.

\begin{proposition}\label{prop:explicitiso}
Let $R$ be a field of characteristic 0 which contains a primitive $l$-th root of 1.
The composition $[\scrL(\Lambda_Q)] \simto \bfL(\Lambda_Q) \to \bfF(Q)$ obtained from
Theorem \ref{thm:minimalcat}  sends the class
of $S(\lambda)^{Q,\,v}_R$ to the standard monomial $|\lambda,Q\rangle$.
\qed
\end{proposition}

For each $p = 1,\ldots, l$, let $\scrI_p$ be the subquiver of $\scrI$ corresponding
to the subset $v^\bbZ Q_p$ of $\scrI$. We define a relation on $\{1,\ldots,l\}$ by 
$i \sim j \iff \scrI_i = \scrI_j$. Let  $\Omega = \{1,\ldots,l\}/\sim$ be the set of 
equivalence classes for this relation. Given $p \in \Omega$, we denote by $Q_p$ the tuple
of $(Q_{i_1},\ldots,Q_{i_r})$ where $(i_1,\ldots,i_r)$ is the ordered set of elements
in $p$. The decomposition $\scrI = \bigsqcup_{p \in \Omega} \scrI_p$ yields a canonical
decomposition of Lie algebras
 $\frakg'_\scrI  = \bigoplus_{p \in \Omega} \frakg'_{\scrI_p}$.
The corresponding decomposition of Fock spaces is given in the following proposition.

\begin{proposition}\label{prop:tensorfock}
The map $|\lambda,Q\rangle_\scrI \longmapsto \otimes_{p \in \Omega} |\lambda^p,
Q_p\rangle_{\scrI_p}$ yields an isomorphism of $\frakg'_I$-modules
$ \bfF(Q)_\scrI \ \mathop{\longrightarrow}\limits^\sim\ \bigotimes_{p\in \Omega}\bfF({Q_p})_{\scrI_p}.$
 \qed
\end{proposition}

\medskip

\subsection{Charged Fock spaces}\label{sec:chargedfock}\hfill\\

A \emph{charged Fock space} is a pair $\bfF(s)=(\bfF(Q),s)$ such that $s\in\bbZ^l$ is a charge
of $Q$, that is  $Q = (v^{s_1},\ldots,v^{s_l})$. Throughout this section, we will always
assume that $\scrI$ is either of type $A_\infty$ or a cyclic quiver. For more general quivers we
can invoke Proposition \ref{prop:tensorfock} to reduce to that case.

\smallskip

\subsubsection{The $\frakg$-action on the Fock space}\label{subsec:xgrading}
The action of $\frakg'$ on $\bfF(Q)$ can be extended to an action of $\frakg$ 
when $Q$ admits a charge $s$. We describe this action in the case where $v$ has finite
order $e$, and $l=1$. In that case $\scrI = v^{\bbZ}$ is isomorphic to the cyclic quiver  $A^{(1)}_{e-1}$
and the charge $s$ is just an integer $d \in \bbZ$ such that $Q = v^d$. If we fix the affine
simple root to be $\alpha_1$, then $\X = \P \oplus \bbZ \delta$ and
$\X^\vee = \Q^\vee \oplus \bbZ \partial$ with $\delta = \sum_{i \in \scrI} \alpha_i$
and $\partial = \Lambda_1^\vee$ (see Example \ref{ex:onegenerator} for more details). 

\smallskip

We define the integer
\begin{align}\label{Delta}\Delta(d,e)=\big(\bar d(1-\bar d/e)+d(d/e-1)\big)/2,\end{align}
where $\bar d$ is the residue of $d$ modulo $e$ in $[0,e-1]$. 
Then, the action of the derivation $\partial$ on $\bfF(d)=(\bfF(Q),d)$ is
 $$ \partial(|\lambda,Q\rangle) = -\big(n_1(\lambda,Q)+\Delta(d,e)\big)|\lambda,Q\rangle.$$
For this action the weight of a standard basis element is 
\begin{equation}\label{eq:weightlevel1}
 \mathrm{wt}(|\lambda,Q\rangle) = \Lambda_Q-\sum_{i\in \scrI} n_i(\lambda,Q)\,\alpha_i-
 \Delta(d,e)\,\delta.
\end{equation}

\smallskip

We now describe the action of the affine Weyl group of $\frakg$ on $\bfF(d)$. 
For $i \in \scrI \smallsetminus\{1\}$, we denote by $\alpha_i^\cl = 2 \Lambda_i - 
\Lambda_{iv}-\Lambda_{iv^{-1}}$ and $\Lambda_i^\cl=\Lambda_i-\Lambda_1$ the $i$-th simple root and fundamental weight of
$\fraks \frakl_e$. These (classical) simple roots span the lattice of classical roots $\Q^\cl$.
The affine Weyl group of $\frakg$ is $W = 
\mathfrak{S}_\scrI \ltimes \Q^\cl$. It acts linearly on $\X$. We will denote by
$t_\gamma$ the action of an element $\gamma \in \Q^\cl$, i.e., for each $\alpha\in\X$ we set
\begin{align}\label{transl}
t_\gamma(\alpha)=\alpha+(\alpha:\delta)\,\gamma-(\alpha:\gamma)\,\delta-\frac{1}{2}(\alpha:\delta)(\gamma:\gamma)\,\delta
\end{align}
where $(\bullet:\bullet)$ is the standard symmetric non-degenerate bilinear form on $\X\times \X$.
Now, we consider the element
$\pi_s=\sum_{i\in\scrI}(s_i-s_{iv})\,\Lambda_{i}$. If $s \in \bbZ^{\scrI}(d)$,
then $\pi_s - \Lambda_Q^\cl \in \Q^\cl$ and we can consider the 
corresponding operator $t_{\pi_s- \Lambda_Q^\cl}$ on $\X$.

\begin{proposition}\label{prop:Delta} 
Let $\lambda$ and $\nu$ be two partitions. 
 Let $s \in \bbZ^\scrI(d)$ be such that $(\emptyset,s) = \tau_e(\lambda_{[e]},d)$, where
$\lambda_{[e]}$ is the $e$-core of $\lambda$.
\begin{itemize}[leftmargin=8mm]
  \item[$\mathrm{(a)}$] The weight of $|\lambda,Q\rangle$ equals
   $ \mathrm{wt}(|\lambda,Q\rangle)  = t_{\pi_s-\Lambda_Q^\cl}(\Lambda_Q)-w_e(\lambda)
   \,\delta.$
  \item[$\mathrm{(b)}$] The weights of $|\lambda,Q\rangle$ and $|\nu,Q\rangle$
  are $W$-conjugate if and only if $w_e(\lambda)=w_e(\nu)$. 
  \qed
\end{itemize}
\end{proposition}

In the particular case where the charge is zero, then
$s = \lambda_{[e]}$, and 
\begin{align}\label{formB} \mathrm{wt}(|\lambda,1\rangle)  = t_{\pi_{\lambda_{[e]}}}(\Lambda_1)-w_e(\lambda)
   \,\delta.\end{align}
Therefore weight spaces are parametrized by pairs $(\nu,w)$ where $\nu$ is an $e$-core
and $w$ is a non-negative integer. The basis element $|\lambda,1\rangle$ is in the weight
space corresponding to $(\lambda_{[e]},w_e(\lambda))$. 

\smallskip

\subsubsection{The crystal of the Fock space}\label{subsec:crystalfock}
We explain here how to associate an abstract crystal to a charged Fock space $\bfF(s)$. 
By Proposition \ref{prop:tensorfock}, we can assume that $\scrI$ is either cyclic or of
type $A_\infty$. We assume that the reader is familiar with \cite{U}.
\smallskip

When $\scrI$ has type $A_\infty$, Uglov's bases coincide with the standard monomial
basis and the discussion is trivial in that case. We will therefore assume that
$v$ has finite order $e$ and $\scrI = v^\bbZ$, so that $\scrI$ has type $A_{e-1}^{(1)}$.

\smallskip

The $\frakg$-module $\bfF(Q)$ admits a quantum deformation $\bfF_{\!u}(s)$ with an
$A$-lattice $\bfF_{\!A}(s)$,  which is a free $A$-module with basis
$\{|\mu,s\rangle\,|\,\mu\in\scrP^l\}$.
It is equipped with an integrable representation of $U_A(\frakg)$.
Note that the action of the Chevalley generators $e_i$ and $f_i$ 
depends on the choice of the charge $s$. The representation of $\frakg$ on $\bfF(s)$ given in
\S \ref{subsec:xgrading} is recovered by specializing the parameter $u$ to $1$.
Uglov has constructed a remarkable $A$-basis $\bfB^+_u(s)=\{\bfb_u^+(\mu,s)\,|\,\mu\in\scrP^l\}$ of $\bfF_{\!A}(s)$ 
which is a lower global basis for the representation of $U_u(\frakg)$ on $\bfF_{\!u}(s)$.

\smallskip

Next, we consider the pairing $(\bullet,\bullet)$ on $\bfF_{\!u}(s)$ defined in \cite[sec.~4.3]{Y05}.
Then, let $\bfB_u^\vee(s)=\{\bfb_u^\vee(\mu,s)\,|\,\mu\in\scrP^l\}$ be the $\bbC(u)$-basis of
$\bfF_{\!u}(s)$ dual to $\bfB_u^+(s)$ with respect to the bilinear form $(\bullet,\bullet)$.
It is an upper global basis of $\bfF_{\!u}(s).$ 
Let $\bfB^\vee(s)=\{\bfb^\vee(\mu,s)\,|\,\mu\in\scrP^l\}$ be the specialization at $u=1$ of
$\bfB_u^\vee(s)$, with the obvious labeling of its elements. 
It is a perfect basis of the usual Fock space $\bfF(s) =  \bfF_u(s)|_{u=1}$ by Proposition \ref{prop:upper-perfect}.

\smallskip

We equip the set of $l$-partitions
with the abstract crystal structure $B(s)=\big(\scrP^l,\widetilde{e_i},\widetilde{f_i}\big)$ defined in \cite{JMMO}.
Let $B(s)=\{b(\mu,s)\,|\,\mu\in\scrP^l\}$ be the obvious labeling.
The map $\bfb^\vee(\mu,s) \in \bfB^\vee(s)\longmapsto b(\mu,s) \in B(s)$ is a crystal isomorphism.

\medskip

\section{Unipotent representations}
\label{part:modular}
In this section we record standard results on unipotent representations of finite
reductive groups in non-defining characteristic. A good reference is \cite{CE}. 

\medskip 

\subsection{Basics}\label{sec:basics}\hfill\\

By an \emph{$\ell$-modular system} we will mean a triple $(K,\scrO,\k)$ where 
$K$ is a field of characteristic zero, $\scrO$ is a complete discrete valuation
ring with fraction field $K$, and $\k$ is the residue field of $\scrO$ with
$\text{char}(\k)=\ell$. When working with representations of a finite group $\Gamma$,
we will always assume that $(K,\scrO,\k)$ is a \emph{splitting $\ell$-modular system
for $\Gamma$}, which means that $K$ and $\k$ are splitting fields for all subgroups
of $\Gamma$. When $\Gamma$ comes from an algebraic group in characteristic $p$,
we will in addition assume that $\ell \neq p$. This case is usually referred to as
the \emph{non-defining characteristic} case.

\smallskip

Let $R$ be any commutative domain (with 1) and $\Gamma$ be a finite group.
We will assume that $p$ is invertible in $R$ and $\scrO$.
Let $R\Gamma$ denote the group ring of $\Gamma$ over $R.$ 
For any subset $S\subseteq\Gamma$ such that $|S|$ is invertible in $R$, let $e_S$ 
be the element $e_S=|S|^{-1}\sum_{g\in S}g$ in $R\Gamma$. When $S$ is a subgroup, $e_S$ is an idempotent. If $R$ is not a field,
an $R\Gamma$-module which is free as an $R$-module will be called an \emph{$R\Gamma$-lattice}.

\smallskip

We shall identify isomorphism classes of irreducible representations over $K$ with their character. When $R$ is a field containing $K$, we denote by $R\,\Irr(K\Gamma)$
or $R\,\Irr(\Gamma)$ the $R$-module of class functions $\Gamma\to R$. It is endowed with the canonical scalar product $\langle\bullet,\bullet\rangle_\Gamma$, for which the set of irreducible characters
$\Irr(K\Gamma)$ of $K\Gamma$ is an orthonormal basis.

\medskip

\subsection{Unipotent $KG$-modules\label{ssec:unip-modules-gen}}\hfill\\

By a \emph{rational group} we will mean a pair $(\bfG, F)$ consisting of an algebraic group $\bfG$ over $\overline\bbF_q$ 
and a Frobenius endomorphism $F$ of $\bfG$.
Assume that $\bfG$ is a connected reductive group.
Fix a parabolic subgroup $\bfP$ of $\bfG$ and an $F$-stable Levi complement
$\bfL$ of $\bfP$. We \emph{do not} assume $\bfP$ to be $F$-stable.
Write $L=\bfL^F$ and $G=\bfG^F$.

\smallskip

Let $R^{\bfG}_{\bfL\subset\bfP}$ and ${}^*\!R^{\bfG}_{\bfL\subset\bfP}$ 
denote respectively the Lusztig induction and restriction maps
from $\bbZ\,\Irr(KL)$ to $\bbZ\,\Irr(KG)$. 
We will assume that the Mackey formula holds for $R^{\bfG}_{\bfL\subset\bfP}$ and
${}^*\!R^{\bfG}_{\bfL\subset\bfP}$, which we know for the groups we will focus on 
later (see \cite{BM11} for more details). 
Under this condition, the Lusztig induction and restriction do not depend on the 
choice of the parabolic subgroup $\bfP$. We abbreviate
$R^{\bfG}_{\bfL\subset\bfP}=R^{\bfG}_{\bfL}$ and ${}^*\!R^{\bfG}_{\bfL\subset\bfP}={}^*\!R^{\bfG}_{\bfL}$.

\smallskip

Let $\bfB$ be an $F$-stable Borel subgroup of $\bfG$ and
$\bfT$ be an $F$-stable maximal torus of $\bfG$ with $\bfT \subset\bfB$. 
Let $\bfN$ be the normalizer of $\bfT$ in $\bfG$.
Write $B=\bfB^F$, $T=\bfT^F$ and $N=\bfN^F$.
The groups $\bfB$, $\bfN$ form a reductive $BN$-pair of $\bfG$ with Weyl group 
$\bfW=\bfW_\bfG$ given by $\bfW=\bfN/\bfT$.
Since $\bfB$, $\bfN$ are stable by $F$ and $\bfG$ is connected, 
the finite groups $B$, $N$ form a split BN-pair of $G$ whose Weyl group
$W=W(T)$ is given by $W=\bfW^F=N/T$.

\smallskip

The $G$-conjugacy classes of
$F$-stable maximal tori of $\bfG$ are parametrized by the $F$-conjugacy classes in $\bfW$.
For each $w\in\bfW$ let  $\bfT_{w}$ be
an $F$-stable maximal torus in the $G$-conjugacy class parametrized by the $F$-conjugacy class of $w$. 
Under conjugation by some element of $\bfG$, the pair $(\bfT_{w}, F)$ is identified with the pair $(\bfT, wF)$. 
In particular, we have $T_w\simeq\bfT^{wF}$ and $W(T_{w})\simeq\bfW^{wF}$.
The virtual characters $R^\bfG_{\bfT_{w}}(1)$ obtained by induction of the 
trivial representation of the tori $T_w$ are called the \emph{virtual Deligne-Lusztig
characters}. 

\begin{definition}
An irreducible $KG$-module is \emph{unipotent} if its character, say $\chi$, occurs 
as a constituent of a virtual Deligne-Lusztig character $R^\bfG_{\bfT_{w}}(1)$ for some element
$w\in W$, \emph{i.e.}, if we have $\langle\chi,R^\bfG_{\bfT_{w}}(1)\rangle_G\neq 0$.
\end{definition}

We denote by $KG\umod$ the full subcategory of $KG\mod$ consisting of the modules
which are sums of irreducible unipotent modules. The objects of this category are
the \emph{unipotent $KG$-modules}.

\begin{remark}\label{rem:twist}
Fix an element $\sigma$ in the normalizer of $\bfL$ in $\bfG$ such that
the conjugation by $\sigma$ yields an element in the normalizer of $\bfW_\bfL$ in $\bfW$.
Let us denote it again by the symbol $\sigma$.
Then, we can equip the subgroup $\bfL\subseteq\bfG$ with the rational structure with the Frobenius endomorphism $\sigma F$
and we have Lusztig induction and restriction maps $R^{\bfG\,,F}_{\bfL\,,\,\sigma F}$ and ${}^*\!R^{\bfG\,,F}_{\bfL\,,\,\sigma F}$ between $\bbZ\,\Irr(K\bfL^{\sigma F})$ and $\bbZ\,\Irr(KG)$.
In the particular case $\sigma=1$ considered above we have $R^{\bfG\,,F}_{\bfL\,,\,F}=R^{\bfG}_{\bfL}$ and 
${}^*\!R^{\bfG\,,F}_{\bfL\,,\,F}={}^*\!R^{\bfG}_{\bfL}.$
If $\bfT_w$ is as above then we have also $R^{\bfG\,,\,F}_{\bfT\,,\,wF}=R^\bfG_{\bfT_{w}}$
and ${}^*\!R^{\bfG\,,\,F}_{\bfT\,,\,wF}={}^*\!R^\bfG_{\bfT_{w}}$.
\end{remark}

\medskip

\subsection{Unipotent $\k G$-modules and $\ell$-blocks}\hfil\\\label{sec:unipblocks}

As a result of the lifting of idempotents, the blocks of $\scrO G$ and
$\k G$ correspond by reduction. Both are usually called the \emph{$\ell$-blocks} of
$G$. For $R=\scrO$ or $\k$, any block $B$ of $R G$ is of the form $B=R G\cdot b,$
where $b$ is a central primitive idempotent of $R G$.
The unit $b$ of $B$ is called the \emph{block idempotent of  $B$}.
We will also call \emph{block of $R G\mod$ associated with $B$}
the Serre subcategory generated by the simple modules on which $b$ acts non-trivially.
The $\ell$-blocks of $G$ induce a partition of $\Irr(KG)$ such that the piece
associated with $B$ is the set of all irreducible characters $\chi$ of $KG$ with
$\chi(b)=\chi(1)$. If $\chi\in\Irr(KG)$, we will write $B(\chi)\subseteq \Irr(KG)$
for the piece containing $\chi$. When there is no risk of confusion, we will also
call $B(\chi)$ an $\ell$-block of $G$.

\begin{definition}\label{def:lblock}
An $\ell$-block of $\scrO G$ is \emph{unipotent} if it contains 
at least one unipotent $KG$-module.
A simple $\k G$-module is \emph{unipotent} if it lies in a unipotent block of $\k G$.
\end{definition}

We denote by $\k G\umod$ the Serre subcategory of $\k G\mod$ generated by
the simple unipotent $\k G$-modules. It corresponds to the sum of unipotent
blocks of $\k G\mod$. The \emph{unipotent $\k G$-modules} are by definition the objects
of this category.

\smallskip

Recall that $(K,\scrO,\k)$ is a splitting $\ell$-modular system for $G$. To this system
one can associate a decomposition map $d_{\scrO G}:[KG\mod]\to[\k G\mod]$.
By \cite[thm. 3.1]{H90}, a simple $\k G$-module is unipotent if and only if it is a
constituent of the $\ell$-reduction of a unipotent $K G$-module, see
also \cite{BM}. In other words, the classes of unipotent modules are 
exactly the image of unipotent characters through the decomposition map.
We will denote by $d_{\scrU}:[KG\umod]\to[\k G\umod]$ the restriction
of this map to unipotent characters.

\begin{theorem}[Geck-Hiss \cite{GH91}, Geck \cite{G93}] \label{prop:isodecmap}
Assume $\ell$ is good for $\bfG$ and that $\ell \nmid |Z(\bfG)/Z(\bfG)^\circ|$.
Then the map $d_{\scrU}$ is a linear isomorphism $[KG\umod]\simto[\k G\umod]$.
\qed
\end{theorem}

Note that if $\bfG$ is semisimple with no component of type $A$, then the condition
$\ell \nmid |Z(\bfG)/Z(\bfG)^\circ|$ is superfluous since every prime divisor of 
$|Z(\bfG)/Z(\bfG)^\circ|$ is a bad prime.

\smallskip

Given a positive integer $f$, let $\Phi_f$ be the $f$th cyclotomic polynomial.
If $\bfT \subset \bfG$ is an $F$-stable torus, then the order of $T:= \bfT^F$ is
a polynomial in $q$ with coefficients in $\mathbb{Z}$. We say that $\bfT$ (and by extension $T$) is a \emph{$\Phi_f$-torus} if this polynomial is a power of $\Phi_f$.
An $F$-stable Levi subgroup $\bfL\subseteq\bfG$ is \emph{$f$-split} 
if it is the centralizer in $\bfG$ of a $\Phi_f$-torus.
A \emph{unipotent $f$-pair} is a pair $(\bfL,\chi)$ where $\bfL$ is an $f$-split
Levi subgroup and $\chi$ is an irreducible unipotent $KL$-module.
The pair $(\bfL,\chi)$ is \emph{$f$-cuspidal} if for every proper
$f$-split Levi subgroup $\bfM\subseteq\bfL$ we have ${}^*\!R^\bfL_\bfM(\chi)=0$.

\smallskip

Now assume that $f$ is the smallest positive integer such that $\ell$ divides $q^f-1$. 
Under the assumption that $\ell$ is good, unipotent $\ell$-blocks correspond
to $f$-cuspidal $f$-pairs.

\begin{theorem}[Brou\'e-Malle-Michel \cite{BMM}, Cabanes-Enguehard \cite{CE}]\label{prop:fcuspidalpairs}
Assume $\ell$ is good for $\bfG$, and $\ell \neq 3$ if $\bfG$ has a component
of type ${}^3 D_4$. 
Then there is a bijection between the $G$-conjugacy classes of unipotent
$f$-cuspidal $f$-pairs and the set of unipotent
$\ell$-blocks of $G$ which takes the class of $(\bfL,\chi)$ to the $\ell$-block
$B_{\bfL,\chi}$ such that the irreducible unipotent characters in $B_{\bfL,\chi}$
are exactly the irreducible constituents of $R^\bfG_\bfL(\chi).$
\qed
\end{theorem}

\medskip

\subsection{Harish-Chandra series}\label{sec:HCseries}\hfill\\

Assume now that the parabolic subgroup $\bfP\subseteq\bfG$ considered in
\S\ref{ssec:unip-modules-gen} is $F$-stable.
In that case the group $L$ is $G$-conjugate to a standard Levi subgroup of $G$.
Let $R_L^{G}$ and ${}^*\!R^{G}_L$ be the corresponding Harish-Chandra induction and restriction functors between $RL\mod$ and $RG\mod$.  
Let $P=\bfP^F$ and $U=\bfU^F$, where $\bfU\subset \bfP$ is the unipotent radical of $\bfP$. 
Notice that the Harish-Chandra induction is the special case of Lusztig induction for 1-split Levi subgroups.

\smallskip

The order of $U$ is a power of $q$, hence it is invertible in $R$.
Thus, for all $M\in RL\mod$, $N\in RG\mod$ we have 
  $$R_L^{G}(M)=RG\cdot e_U\otimes_{RL}M \quad \text{and} \quad
  {}^*\!R_L^{G}(N)=e_U\cdot RG\otimes_{RG}N.$$
We will say that the functors $R_L^{G}$ and ${}^*\!R^{G}_L$ are \emph{represented} by the
$(RG,RL)$-bimodule $RG\cdot e_U$ and the $(RL,RG)$-bimodule $e_U\cdot RG$ respectively.

\smallskip

Here are some well-known basic properties of the functors $R_L^{G}$, ${}^*\!R^{G}_L$. 
\begin{itemize}[leftmargin=8mm]
  \item[(a)] $R_L^{G}$, ${}^*\!R^{G}_L$ do not depend on $P$,
  \item[(b)] $R_L^{G}$, ${}^*\!R^{G}_L$ are exact and left and right adjoint
  to one another,
  \item[(c)] if $L\subseteq M\subseteq G$ are Levi subgroups there are isomorphisms of
  functors $R_L^{G}=R_M^{G}R_L^{M}$ and ${}^*\!R_L^{G}={}^*\!R_L^{M}{}^*\!R_M^{G}$.
\end{itemize}

Let $R=K$ or $\k$. An irreducible $RG$-module $E$ is \emph{cuspidal} if 
${}^*\!R^{G}_L(E)=0$ for all proper standard Levi subgroups $L\subsetneq G$.
A \emph{cuspidal pair of $RG$} is a pair $(L,E)$ where $L$ is as above and
$E\in\Irr(RL)$ is cuspidal. 
Since the group $L$ is uniquely recovered from $E$, from now on we may omit it from the notation.
Then, the set $\Irr(R G,E)\subseteq\Irr(R G)$ consisting of the constituents of the head of $R^G_L(E)$ 
is equal to the set of constituents of the socle of $R^G_L(E)$ and
is called the \emph{Harish-Chandra series} of $(L,E)$.
The $R$-algebra $\scrH(RG,E)=\End_{RG}(R^{G}_L(E))^\op$ is
its \emph{ramified Hecke algebra}.
We have the following facts: 
\begin{itemize}[leftmargin=8mm]
\item[(d)] the Harish-Chandra series form a partition of $\Irr(R G)$,
\item[(e)] the functor $\frakF_{R^{G}_L(E)}$ as defined in \S\ref{ssec:rings-cat}
yields a bijection 
$$\Irr(R G,E)
\mathop{\longleftrightarrow}\limits^{1:1} \Irr(\scrH(RG,E)).$$
\end{itemize}

\smallskip

\begin{proposition}[\cite{BoRo03}]\label{prop:stable} If $R=K$ or $\k$, then 
the Harish-Chandra induction and restriction functors preserve the category
of unipotent $RG$-modules.  \qed
\end{proposition}

\medskip

\section{Finite unitary groups}\label{part:unitarygroups}

This section is devoted to the construction of categorical actions on the
category of unipotent representations of finite unitary groups $\mathrm{GU}_n(q)$.

\smallskip

Let $R$ be a commutative domain with unit. 
Under mild assumptions on $R$, we construct in \S\ref{sec:rep-datum} a representation
datum on the abelian category 
 $$ RG\mod := \bigoplus_{n \geqslant 0} \mathrm{GU}_n(q)\mod$$
given by Harish-Chandra induction and restriction. 
When $R$ is an extension of $\bbQ_\ell$ or $\bbF_\ell$, the categorical datum
restricts to the category $\scrU_R$ of unipotent representations of $RG\mod$. 
On this smaller category, the eigenvalues of $X$ are powers of $-q$, which hints
that the Lie algebra $\frakg$ that should act on $[\scrU_R]$ corresponds
to the quiver with vertices $(-q)^\bbZ$ and arrows given by multiplication by
$q^2$. 

\smallskip

We prove that the representation datum lifts indeed to a categorical action of $\frakg$
on $\scrU_R$.  The situation depends on the parity of $e$, the order of
$-q$ modulo $\ell$. When $e$ is even, i.e., $\ell$ is a linear prime, $\frakg$ is a subalgebra of
$(\widehat{\fraks\frakl_{e/2}})^{\oplus 2}$ and each ordinary Harish-Chandra series
categorifies a level 2 Fock space for $\frakg$. When $e$ is odd (unitary prime
case), $\frakg \simeq  \widehat{\fraks\frakl_\e}$ and weight spaces correspond
to unipotent $\ell$-blocks, which are now transverse to the ordinary Harish-Chandra series.

\smallskip

For studying the weight space decomposition of $[\scrU_R]$ as well as the action of 
the affine Weyl group we use the action of a bigger Lie algebra $\frakg_\circ$, which 
comes from Harish-Chandra induction and restriction for general linear groups. 
Going from linear to unitary groups by Ennola duality introduces signs for this action
(see Lemma \ref{lem:isofocks}). It is to be expected that the action of $\frakg_\circ$
on $[\scrU_R]$ lifts to an action by triangulated functors on $D^b(\scrU_R)$ coming
from Lusztig induction, although we will not use it. See \cite{DVV2} for more details.

\medskip

\subsection{Definition}\label{sec:def-unitary}\hfill\\

Fix a positive integer $n$.
We equip the reductive algebraic group $\bfG\bfL_n=\GL_n(\overline\bbF_q)$ with
the standard Frobenius map $F_q:\bfG\bfL_n\to\bfG\bfL_n$, $(a_{ij})\mapsto (a_{ij}^q)$
given by raising every coefficient to the $q$th power. The finite general linear group
$\mathrm{GL_n}(q)$ is given by the fixed points of $\bfG\bfL_n$ under $F_q$. 
In this section we will work with a twisted version of this group obtained by
twisting the Frobenius map.  Let $J_n$ be the
$n\times n$ matrix with entry $1$ in $(i,n-i+1)$ and zero elsewhere. We will often
write $J = J_n$ when there is no risk of confusion on the size of the matrices.
We define a new Frobenius map $F$ on $\bfG\bfL_n$, called the \emph{twisted} Frobenius
map, by setting $F=F_q\circ\alpha$ where $\alpha(g)=J\cdot{}^tg^{-1}\cdot J$ for each
$g\in \bfG\bfL_n$. The \emph{finite unitary group $G_n=\GU_n(q)$} is then given by
$$G_n=(\bfG\bfL_n)^F=\{g\in\bfG\bfL_n\,;\,F(g)=g\}.$$
Since $F^2=F_{q^2}$ we have $G_n\subset GL_n$, where we abbreviate 
$$GL_n=\bfG\bfL_n(q^2):=(\bfG\bfL_n)^{F^2}.$$ By convention we also define $G_0 = \{1\}$
to be the trivial group.

\smallskip
We equip $\bfG\bfL_n$ with the standard split BN-pair such that $\bfB$ is the
subgroup of upper triangular matrices and $\bfN$ is the subgroup of all monomial
matrices. Since $\bfB$, $\bfN$ are stable by $F$ and $\bfG\bfL_n$ is connected, the
groups $B=\bfB^F$, $N=\bfN^F$ form a split BN-pair of the finite group $G_n$.
Let $\bfT$ be the diagonal torus in $\bfG\bfL_n$ and $T=\bfT^F$.
Let $\bfW=\bfW_n$ be the Weyl group of $\bfG\bfL_n$ and $W=W_n$ be the Weyl group of
$G_n$. We have $\bfW\simeq\frakS_n$, and $F$ induces on $\bfW$ the automorphism
given by conjugation with the longest element $w_0$. We will embed $\bfW$ in $G_n$
using permutation matrices, so that $w_0$ corresponds to $J$. We have $W=\bfW^F=
C_\bfW(w_0)$. It is a Weyl group of type $B_m$ if $n=2m$ or $2m+1$.

\smallskip

Let $\varepsilon_1,\dots,\varepsilon_n$ be the characters of $\bfT$ such that
$t=\diag(\varepsilon_1(t), \varepsilon_2(t), \dots, \varepsilon_n(t))$ for all $t \in \bfT$.
The roots (resp. simple roots) of $\bfG\bfL_n$ are given by $\{\varepsilon_i-\varepsilon_j \, | \, i \neq j\}$ (resp. $\{\varepsilon_r-\varepsilon_{r+1}\}$).
Let $s_r=(r,r+1)$ be the simple reflection in $\bfW$ associated with the simple root
$\alpha_r = \varepsilon_r-\varepsilon_{r+1}$. The action of $F$ on the roots
induces an automorphism $\sigma$ of the Dynkin diagram of $\bfG\bfL_n$ such that
 $$ F \circ \alpha_r^\vee= q \sigma (\alpha_r)^\vee\quad \text{with} \quad \sigma(\alpha_r)=\alpha_{n-r}.$$
For each root $\alpha$, let $\bfU_\alpha$ and $\alpha^\vee\in\Hom(\bbG_m,\bfT)$
be the corresponding root subgroup and cocharacter. We also choose an
isomorphism $u_\alpha \, : \, \mathbb{G}_a \simto \bfU_\alpha$ such that
$F(u_{\sigma(\alpha)}(t))=u_{\alpha}(t^q)$. Note that a one-parameter subgroup
of $T$ has either $(q^2-1)$ or $(q+1)$ elements, and a root subgroup of $G_n$ has either $q$ or $q^2$ elements.

\smallskip

The standard Levi subgroups $L_{r,\lambda}$ of $G_n$ are parametrized by pairs
$(r,\lambda)$ where $r$ is a non-negative integer and $\lambda=(\lambda_1,\lambda_2,\dots,\lambda_s)$ is a tuple of
positive integers such that $n=r+2\sum_{u=1}^s\lambda_u$. The group $L_{r,\lambda}$ consists
of all matrices of $\bfG\bfL_n$ which belong to $G_n$ and are of block-diagonal type 
$$\prod_{u=1}^s\bfG\bfL_{\lambda_u}\times\bfG\bfL_r\times \prod_{u=s}^1\bfG\bfL_{\lambda_u}.$$ 
Consequently we have a group isomorphism $L_{r,\lambda}\simeq G_r\times\prod_uGL_{\lambda_u}.$
If $m$ is a positive integer we abbreviate $L_{r,1^m}=L_{r,(1^m)}$ and
$L_{r,m}=L_{r,(m)}$.

\medskip

\subsection{The representation datum on $RG$-mod}
\label{sec:rep-datum}\hfill\\

Let $R$ be a commutative domain with unit. We assume that $q(q^2-1)$ in invertible
in $R$. Using parabolic induction and restriction, we show in this section how
to construct a representation datum on 
 $$ RG\mod:=\bigoplus_{n\in\bbN} RG_n\mod.$$

Fix a positive integer $n$. Parabolic (or Harish-Chandra) induction provides functors
between $L\mod$ and $G_n\mod$ for any standard Levi subgroup $L = L_{r,m}
\subset G_n$. Since we want functors between $G_{r}\mod$ and $G_n\mod$
we will consider a slight variation of the usual parabolic induction.

\smallskip 

Let $0 \leqslant r < n$. We denote by $V_r$ the unipotent radical of the standard
parabolic  subgroup $P_{r,1} \subset G_{r+2}$ with Levi complement $L_{r,1}$.
Let $U_{r}\subset G_{r+2}$ be the subgroup given by
\begin{align}\label{U}
  U_{r}=V_{r}\rtimes\bbF^\times_{q^2}= 
  \begin{pmatrix}*&*&\cdots&\cdots&*\\
	&1&&(0)&\vdots\\
	&&\ddots&&\vdots\\
	&(0)&&1&*\\
	&&&&*
  \end{pmatrix}
\end{align}
so that $P_{r,1} =V_{r} \rtimes  L_{r,1}  \simeq U_r \rtimes G_r$.
If $n-r$ is even, we set $U_{n,r}=U_{n-2}\rtimes\dots\rtimes U_{r}$ and we define
$e_{n,r}=e_{U_{n,r}}$ to be the corresponding idempotent. In particular,
we have $e_{r+2,r}=e_{U_{r}}$. We embed $G_r$ into the Levi subgroup
$L_{r,1^m}=G_r\times GL_1^{m}$ in the obvious way. This yields an embedding
$G_r\subset G_n$ and functors
$$\begin{aligned}
  F_{n,r}=& \, RG_n\cdot e_{{n,r}}\otimes_{RG_r}\bullet\,:\,RG_r\mod\longrightarrow RG_n\mod,\\
  E_{r,n}=& \, e_{{n,r}}\cdot RG_n\otimes_{RG_n}\bullet\,:\,RG_n\mod\longrightarrow RG_r\mod.
\end{aligned}$$
Note that $F_{n,r}$ can be seen as the composition of the inflation $G_r\mod
\longrightarrow L_{r,1^m}\mod$ with the parabolic induction from $L_{r,1^m}$ to $G_n$. 

\smallskip

An endomorphism of the functor $F_{n,r}$ can be represented by an $(RG_n,RG_r)$-bimodule
endomorphism of $RG_n\cdot e_{{n,r}}$, or equivalently by an element of 
$e_{{n,r}}\cdot RG_n\cdot e_{{n,r}}$ centralizing $RG_r$. Thus, the elements
$$ X_{r+2,r}= (-q)^r e_{{r+2,r}} (1,r+2)\, e_{{r+2,r}},\qquad
  T_{r+4,r}= q^2 e_{{r+4,r}}(1,2)(r+3,r+4)\,e_{{r+4,r}}$$
define respectively natural transformations of the functors $F_{r+2,r}$ and $F_{r+4,r}$.
We set 
$$F=\bigoplus_{r\geqslant 0}F_{r+2,r},\qquad 
  X=\bigoplus_{r\geqslant 0}X_{r+2,r},\qquad
  T=\bigoplus_{r\geqslant 0}T_{r+4,r}.$$

\begin{proposition}\label{prop:pre-cat}
  The tuple $(E,F,X,T)$ defines a representation datum with parameter $q^2$ on $RG\mod =
  \bigoplus_{n\in\bbN} RG_n\mod$. \qed
\end{proposition}

\begin{remark}
The reader might argue that we did not check that $X$ was invertible. 
Nevertheless, we shall only be working with unipotent modules over
a field, in which case the eigenvalues of $X$ are powers of $-q$, hence
non-zero (see Theorem \ref{thm:cat}). This will ensure that the
restriction of $X$ to this category is indeed invertible. 
\end{remark}

\medskip

\subsection{The categories of unipotent modules $\scrU_K$ and $\scrU_\k$}
\label{sec:uK-and-uk}\hfill\\

From now on, we fix a prime number $\ell$ such that $\ell \nmid q\,(q^2-1)$,
and an $\ell$-modular system $(K,\scrO,\k)$ with $\bbQ_\ell\subset K$ and 
$\overline{\bbF}_\ell \subset \k$. Although this modular system is not large enough
for all the finite groups encountered, it will be enough for our purpose since
we will be working with unipotent characters only, which are defined over $\bbQ_\ell$ (see below).

\smallskip

Throughout the following sections, we will denote by $d$, $e$ and $f$ the order
of $q^2$, $-q$ and $q$ in $\k^\times$. In particular $e \neq 1,2$. If $e$ is odd,
then $d=e$ and $f=2e$; if $e$ is even, then $d = e/2$ and either $f= e/2$ if
$e \equiv 2$ mod $4$ or $f=e$ if $e \equiv 0$ mod $4$. 

\smallskip

\subsubsection{The category $\scrU_K$}\label{subsec:uK}
Fix a positive integer $n$. By \cite{LS77}, the irreducible unipotent $KG_n$-modules
are labelled by partitions of $n$. Their character can be directly constructed from
the virtual Deligne-Lusztig characters. Namely, for each $w\in\frakS_n$, fix an $F$-stable
maximal torus $\bfT_{w}\subseteq\bfG\bfL_n$ in the $G_n$-conjugacy class parametrized by 
the $\frakS_n$-conjugacy class of $wF$ in the coset $\frakS_n F$, with the convention that $\bfT_1 = \bfT$. Then the class function
\begin{align}\label{unip.char}  \chi_\lambda=|\frakS_n|^{-1}\sum_{w\in\frakS_n}\phi_\lambda(ww_0)
  R_{\bfT_{w}}^{\bfG\bfL_n}  (1)\in\bbZ\Irr(KG_n)
\end{align}
is, up to a sign, an irreducible unipotent character. 
Since the cohomology groups are defined over $\bbQ_\ell \subset K$ and since any character of $\frakS_n$ has rational values, either $\chi_\lambda$ or $- \chi_\lambda$  is the character of some irreducible $KG_n$-module which we denote by $E_\lambda$. By abuse of notation we will still denote by $E_\lambda$ its isomorphism class. 
Let $\varepsilon_\lambda$ be the sign determined by the relation $\chi_\lambda=\varepsilon_\lambda\,E_\lambda$.

\smallskip

Recall that $G_0=\{1\}$. We call the category of unipotent $KG$-modules the category  
  $$\scrU_K:=\bigoplus_{n\in\bbN}K G_n\umod.$$
This category is abelian semisimple. 
From the previous discussion we have 
$$\Irr(\scrU_K)=\{E_\lambda\,|\,\lambda\in\scrP\},$$
where by convention $\Irr(KG_0)=\{E_\emptyset\}$.

\smallskip

\subsubsection{The category $\scrU_\k$}\label{subsec:uk}
Using the $\ell$-modular system we have decomposition maps $d_{\scrO G_n}$
which by Proposition \ref{prop:isodecmap} restrict to linear isomorphisms 
 $$ d_{\scrO G_n}:[KG_n\umod]\simto[\k G_n\umod]. $$
Note that we do not need to assume that $KG_n$ is a split algebra since we are
only interested in unipotent $KG_n$-modules, which are all absolutely irreducible.
\smallskip

For unitary groups, this map is actually unitriangular with respect 
to the basis of irreducible modules and the dominance order.

\begin{theorem}[\cite{G91}]\label{prop:unitriangularunip} 
  There is a unique labelling $\{D_\lambda\}_{\lambda\in\scrP_n}$ of the unipotent
  simple $\k G_n$-modules such that
    $d_{\scrO G_n}([E_\lambda]) \in \ [D_\lambda] + \sum_{\nu > \lambda} \bbZ[D_\nu]$
  where $>$ is the dominance order on partitions of $n$. 
  \qed
\end{theorem}

We call the category of unipotent $\k G$-modules the category  
  $$\scrU_\k:=\bigoplus_{n\in\bbN}\k G_n\umod.$$
This is an abelian category which is not semisimple. The isomorphism classes of
simple objects are $\Irr(\scrU_\k)=\{D_\lambda\}_{\lambda\in\scrP}$ and the
decomposition map yields a $\bbZ$-linear isomorphism $d_\scrU : [\scrU_K]\simto[\scrU_\k]$.
The following result is a consequence of the unitriangularity of this map. 

\begin{proposition}\label{prop:liftingstandards}
Given a partition $\lambda \vdash n$ there is a unique $\scrO G_n$-lattice
  $\widetilde E_\lambda$ such that 
  \begin{itemize}[leftmargin=8mm]
   \item[$\mathrm{(a)}$] $K\widetilde E_\lambda=E_\lambda$ as a $K G_n$-module,
   \item[$\mathrm{(b)}$] $V_\lambda:=\k\widetilde E_\lambda$ is an indecomposable $\k G_n$-module with
   head isomorphic to $D_\lambda$,
   \item[$\mathrm{(c)}$] $d_\scrU([E_\lambda])=[V_\lambda].$
   \qed
  \end{itemize}
\end{proposition}

\smallskip

\subsubsection{Blocks of $\scrU_\k$}\label{subsec:blocksofuk}
A block of $\scrU_\k$ is an indecomposable summand of $\scrU_\k$. Therefore 
blocks of $\scrU_\k$ correspond to the unipotent blocks of $\k G_n$
where $n$ runs over $\bbN$. 
Recall that $e$ is the order of $-q$ modulo $\ell$. 

\begin{proposition}[\cite{FS82}]\label{prop:lblocks}
The map $E_\lambda \longmapsto (\lambda_{[e]},w_e(\lambda))$
yields a bijection between unipotent 
$\ell$-blocks and pairs $(s,w)$ where $s \in \bbZ^e(0)$ and $w \in \bbN$.
\qed
\end{proposition}

Recall that $\nu \mapsto \nu_{[e]}$ induces a bijection between $e$-cores
and $\bbZ^e(0)$. Given $\nu$ an $e$-core we will denote by $B_{\nu,w}$ or $B_{\nu_{[e]},w}$
the unipotent $\ell$-block containing all the unipotent characters $E_\lambda$ such that 
$\nu$ is the $e$-core of $\lambda$ and $w_e(\lambda) =w$. It is a block of $\k G_n$
with $n = |\nu|+ew$.

\smallskip

It also follows from the classification of blocks that when $w  < e$ and $e \neq 1$,
the defect group of $B_{\nu,w}$ is an elementary abelian $\ell$-group of rank $w$ (see
\cite{BMM}). In particular, when $w=0$ the module $D_\lambda$ is simple and projective, 
isomorphic to $V_\lambda$, and when $w =1$ the defect group of $B_{\nu,w}$ is a cyclic
group. The structure of such blocks was determined in \cite{FS}.

\smallskip

\subsubsection{The weak Harish-Chandra series}\label{subsec:defweakhc}
For this section we assume that $R$ is one of the fields $K$ or $\k$. 

\smallskip

Let $r,m \geqslant 0$ and $n = r+2m$. The inflation from $G_r$ to
$L_{r,1^m}$ yields an equivalence between $RG_r\umod$ and $RL_{r,1^m}\umod$, since
the Deligne-Lusztig varieties depend only on the semisimple type of the reductive
group. This equivalence intertwines the functors $E_{n,r}$, $F_{n,r}$ with the
parabolic restriction and induction ${}^*R_{L_{r,1^m}}^{G_n}$, $R_{L_{r,1^m}}^{G_n}$.
Therefore working with $\scrU_R$ and the functors $E$ and $F$ is the same as
working in the usual framework of unipotent representations and Harish-Chandra
induction/restriction from Levi subgroups. Note however that we do not
consider all the standard Levi subgroups, but only the ones that are conjugate
to $L_{r,1^m}$. Therefore we need to consider a slight variation of the
usual Harish-Chandra theory. The following definitions and properties are taken from
\cite{GHJ}.

\begin{definition} Fix a non-negative integer $n$. 
 \begin{itemize}[leftmargin=8mm]
  \item[(a)] An $R G_n$-module $D$ is \emph{weakly cuspidal} if ${}^*\!R^{G_n}_L(D)=0$
   for any Levi subgroup $L\subsetneq G_n$ which is $G_n$-conjugate to a subgroup of
   the form $L_{r,1^m}$.
  \item[(b)] A \emph{weakly cuspidal pair} of $R G_n$ is a pair which is $G_n$-conjugate
  to $(L_{r,1^m}\,,\, D)$ for some $n=r+2m$ and some weakly cuspidal irreducible $R G_r$-module
  $D$ which is viewed as a $\k L_{r,1^m}$-module by inflation.
  \item[(c)] The \emph{weak Harish-Chandra series} $\wIrr(R G_n, D)$ of $R G_n$
  determined by the weakly cuspidal pair $(L,D)$ is the set of the constituents of the head
  of $R_L^{G_n}(D)$. By \cite{GHJ}, it coincides with the set of the constituents of its socle.
 \end{itemize}
\end{definition}

If $D$ is a unipotent $RG_r$-module, then $D$ is weakly cuspidal if and only if $E(D) =0$.
Moreover, if $D$ is irreducible, the weak Harish-Chandra series coincides with the set
of irreducible constituents in the head of $F^m D$. Therefore it makes sense to define
the weak Harish-Chandra series of $D$ in $\scrU_R$ by
  $$\wIrr(R G,D)=\bigsqcup_{n\geqslant r}\wIrr(R G_n,D).$$
It is the set of irreducible consituents in the head of $F^kD$ for some $k \geqslant 0$
(or equivalently in the socle). 
As in the case of the usual theory, the weak Harish-Chandra series of $\scrU_R$ form
a partition of $\Irr(\scrU_R)$, i.e., we have 
$\Irr(\scrU_R) =\bigsqcup \wIrr(R G,D)$
where the sum runs over the set of isomorphism classes of weakly cuspidal unipotent
irreducible modules $D$. 

\smallskip

When $R=K$, weakly cuspidal unipotent modules coincide with cuspidal modules and were
determined in \cite{L77}. One of the main results
in this paper is the classification of the weakly cuspidal modules in the case where $R =\k$.

\medskip

\subsection{The $\frakg_\infty$-representation on $\scrU_K$}
\label{sec:g-infty}\hfill\\

In this section we show how the categorical datum defined in \S \ref{sec:rep-datum} yields
a categorical representation on $\scrU_R$ in the case where $R=K$. 
This is achieved by translating into our framework the theory of Lusztig and Howlett-Lehrer
on endomorphism algebras of representations induced from a cuspidal module in \cite{L77} and \cite{HL80}. 

\smallskip

\subsubsection{Action of $E$ and $F$}
Since every unipotent character is a linear combination of (virtual) Deligne-Lusztig characters
$-$ we call such a class function a \emph{uniform function} $-$ it is well-known how to
compute the action of $E$ and $F$ on the category $\scrU_K$. 

\begin{lemma}\label{lem:indresunip}
  Let $\lambda$ be a partition. Then 
  $ [E] ([E_\lambda]) = \sum_{\mu} [E_\mu]$ and $ [F] ([E_\lambda]) = \sum_{\mu'} [E_{\mu'}]$
  where $\mu$ (resp. $\mu'$) runs over the partitions that are obtained from $\lambda$
  by removing (resp. adding) a $2$-hook. \qed
\end{lemma}

We deduce that
 a unipotent $KG_n$ module $E_\lambda$ is weakly cuspidal if and only if there exists
 $t \geqslant 0$ such that $n=t(t+1)/2$ and $\lambda = \Delta_t = (t,t-1,\ldots,1)$.
We will set $E_t := E_{\Delta_t}$. Consequently, the weakly cuspidal pairs are,
up to conjugation, the pairs $(L_{r,1^m},E_t)$ with $r =t(t+1)/2$. The partition
into series is thus
  $\Irr(\scrU_K) =\bigsqcup_{t \in \bbN} \wIrr(K G,E_t).$

\smallskip

\subsubsection{The ramified Hecke algebra}\label{subsec:hl-iso}
Let $r,m$ be non-negative integers and $n = r+2m$. As mentioned in \S \ref{subsec:defweakhc},
the inflation
from $G_r$ to $L_{r,1^m}$ yields an equivalence between $KG_r\umod$ and
$KL_{r,1^m}\umod$ which intertwines the functor $F_{n,r}$ with the parabolic
induction $R_{L_{r,1^m}}^{G_n}$. In particular, we have a canonical isomorphism
  $$\scrH(KG_n,E_t):= \End_{KG_n}(F^m(E_t))^\op \simto 
  \End_{K G_n}(R^{G_n}_{L_{r,1^m}}(E_t))^\op.$$
Now recall from \S \ref{subsec:hecke} that to the categorical datum $(E,F,X,T)$ is attached
a map $\phi_{F^m}:\bfH_{K,m}^{q^2}\to\End(F^m)$. The evaluation of this map at the
module $E_t$ yields a $K$-algebra homomorphism
  $$\phi_{K,m}:\bfH^{q^2}_{K,m}\to\scrH(KG_n,E_t),\quad 
  X_k\mapsto X_k(E_t),\quad T_l\mapsto T_l(E_t).$$
By \cite{HL80}, $\scrH(KG_n,E_t)$ is isomorphic to a Hecke algebra $\bfH^{Q_t;\,q^2}_{K,m}$
of type $B_m$ with 
\begin{align}\label{Q}
  Q_t=\begin{cases}\big((-q)^{-1-t}\,,\,(-q)^t\big)&\ \text{if}\ t\ \text{is\ even},\\
  \big((-q)^t\,,\,(-q)^{-1-t}\big)&\ \text{if}\ t\ \text{is\ odd}.
  \end{cases}
\end{align}
The previous map provides indeed such an isomorphism.

\begin{theorem}\label{thm:cat}
Let $t,m \geqslant 0$ and $n = t(t+1)/2+2m$. Then the map $\phi_{K,m}$ factors through
a $K$-algebra isomorphism 
  $\bfH^{Q_t;\,q^2}_{K,m}\mathop{\longrightarrow}\limits^\sim\scrH(KG_n,E_t).$
  \qed
\end{theorem}

\smallskip

\subsubsection{Parametrization of the weak Harish-Chandra series of $\scrU_K$}
\label{subsec:hcseries-char0}
Let $W(B_m)$ be the Weyl group of type $B_m$, and $t_0,t_1,\ldots,t_{m-1}$ be the
generators corresponding to the following Dynkin diagram
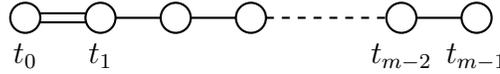
\begin{align*}
\begin{split}
\begin{tikzpicture}[scale=.4]
\draw[thick] (0 cm,0) circle (.5cm);
\node [below] at (0 cm,-.5cm) {$t_0$};
\draw[thick] (0.5 cm,-0.15cm) -- +(1.5 cm,0);
\draw[thick] (0.5 cm,0.15cm) -- +(1.5 cm,0);
\draw[thick] (2.5 cm,0) circle (.5cm);
\node [below] at (2.5 cm,-.5cm) {$t_1$};
\draw[thick] (3 cm,0) -- +(1.5 cm,0);
\draw[thick] (5 cm,0) circle (.5cm);
\draw[thick] (5.5 cm,0) -- +(1.5 cm,0);
\draw[thick] (7.5 cm,0) circle (.5cm);
\draw[dashed,thick] (8 cm,0) -- +(4 cm,0);
\draw[thick] (12.5 cm,0) circle (.5cm);
\node [below] at (12.5 cm,-.5cm) {$t_{m-2}$};
\draw[thick] (13 cm,0) -- +(1.5 cm,0);
\draw[thick] (15 cm,0) circle (.5cm);
\node [below] at (15 cm,-.5cm) {$t_{m-1}$};
\end{tikzpicture}
\end{split}
\end{align*}
Let us first recall the construction of the irreducible characters of $W(B_m)$. We denote by 
$\sigma_m$ the linear character of $W(B_m)$ such that $\sigma_m(t_0) = -1$
and $\sigma_m(t_i) = 1$ for all $i>0$. Given $\lambda \vdash m$ a partition 
of $m$, we write $\widetilde{\phi}_\lambda$ for the inflation to $W(B_m)$
of the irreducible character of $\frakS_m$ corresponding to $\lambda$. 
Given $a,b$ such that $a + b = m$, one can consider the subgroups
$W(B_a)\times \frakS_b \subset W(B_a) \times W(B_b)$ of $W(B_m)$ where
$W(B_a)\times \frakS_b$ is a parabolic subgroup generated by $\{t_0,t_1,\ldots,t_{a-1}\}\cup \{t_{a+1},\ldots,t_{m-1}\}$, and $W(B_a) \times W(B_b)$ is obtained by adding the
reflection $t_{a}\cdots t_0 \cdots t_a$. The irreducible character of $W(B_m)$
associated to a bipartition $(\lambda,\mu)$ of $m$ is
 $$ \scrX_{\lambda,\mu} = \mathrm{Ind}_{W(B_{|\lambda|}) \times W(B_{|\mu|})}^{W(B_m)}
  \, (\widetilde{\phi}_\lambda \boxtimes \sigma_{|\mu|} \widetilde{\phi}_\mu).$$

\smallskip

By Tits deformation theorem, the specialization $q \mapsto 1$ yields a bijection 
$$ \Irr( \bfH^{Q_t;\,q^2}_{K,m}) \, \simto \, \Irr(W(B_m)) $$
from which we obtain a canonical labelling of the irreducible representations of
$\bfH^{Q_t;\,q^2}_{K,m}$. We write $\Irr( \bfH^{Q_t;\,q^2}_{K,m}) = 
\{S(\lambda,\mu)^{Q_t;\,q^2}_K\}_{(\lambda,\mu) \in \scrP^2_m}$.
Setting
$T_0 = (-1)^t q^{-1-t} X$, we have now the quadratic relation
 $$(T_0 + 1)(T_0-q^{2t+1}) = 0.$$
Using this generator instead of $X$, we obtain the usual presentation for a Hecke algebra
of type $B_m$ with parameters $(q^{2t+1},q^2)$. The endomorphism of $KW(B_m)$ which is obtained
from the renormalization $\bfH^{Q_t;\,q^2}_{K,m} \simto  \bfH^{(q^{2t+1},1);\,q^2}_{K,m}$
at $q=1$ is the identity on $\frakS_m$ but sends $t_0$ to $(-1)^tt_0$. Therefore 
this renormalization sends $S(\lambda,\mu)$ to $S(\lambda,\mu)$ if $t$ is even, and
to $S(\mu,\lambda)$ if $t$ is odd. 

\begin{corollary}\label{cor:bijection} 
Let $t,m \geqslant 0$ and $n = t(t+1)/2+2m$. Then the map $\phi_{K,m}$ and the functor $\frakF_{F^m(E_t)}$
induce a bijection
  $ \wIrr(K G_n,E_t) \mathop{\longleftrightarrow}\limits^{1:1} \Irr( \bfH^{Q_t;\,q^2}_{K,m})$
sending $E_\lambda$ to $S(\lambda^{[2]})^{Q_t,q^2}_K$ for all partitions $\lambda \vdash n$ with
$2$-core $\Delta_t = (t,t-1,\ldots,1) $. 
\qed
\end{corollary}

\smallskip

\subsubsection{The $\frakg_\infty$-representation on $\scrU_K$}
\label{subsec:g-infty}
The functors $E,F$ preserve the subcategory $\scrU_K$ by Proposition \ref{prop:stable},
hence $(E,F,X,T)$ yields a representation datum on $\scrU_K$. In order to extend it to a
categorical representation on $\scrU_K$, one should consider the quiver
$\scrI(q^2)$ with vertices given by the various eigenvalues of $X$ and arrows
$i \longrightarrow q^2i$.
In this section we will view the integer $q$ as an element
of $K^\times$ in the obvious way. 

\begin{definition}\label{def:g-infty} 
Let $\scrI_\infty$ denote the subset $(-q)^\bbZ$ of $K^\times.$
We define $\frakg_{\infty}$  to be the (derived)
Kac-Moody algebra associated to the quiver  $\scrI_\infty(q^2)$.
\end{definition}

To avoid cumbersome notation, we will write for short $\scrI_{\infty} = \scrI_\infty(q^2)$,
and $(-)_\infty=(-)_{\scrI_\infty}$.
We denote by $\{\Lambda_i\}$, $\{\alpha_i\}$ and $\{\alpha_i^\vee\}$ the fundamental
weights, simple roots and simple coroots of $\frakg_\infty$. Here $\X_\infty$ coincides with
$\P_{\infty} = \bigoplus \bbZ \Lambda_i$. Consequently, there is a Lie algebra isomorphism $(\fraks\frakl_\bbZ)^{\oplus 2} \simto \frakg_{\infty}$ such that $(\alpha^\vee_d,0)\longmapsto\alpha^\vee_{-q^{2d-1}}$ and $(0,\alpha^\vee_d)\longmapsto\alpha^\vee_{q^{2d}}$.

\smallskip

For any $t,m,n\in\bbN$, let $(KG_n,E_t)\mod$ be the Serre subcategory of $\scrU_K$ generated by
the modules $F^m(E_t)$ with $n=r+2m$ and $r=t(t+1)/2$. We define
  $$\scrU_{K,t}=\bigoplus_{n\geqslant 0}\, (KG_n,E_t)\mod.$$
Then  $\scrU_K=\bigoplus_{t\geqslant 0}\scrU_{K,t}$. We obtain our first categorification result.

\begin{theorem}\label{thm:char0}
 Let $t \geqslant 0$ and $Q_t$ be as in \eqref{Q}. 
 \begin{itemize}[leftmargin=8mm]
  \item[$\mathrm{(a)}$] 
  The Harish-Chandra induction and restriction functors yield a representation of
  $\frakg_{\infty}$ on $\scrU_{K,t}$ which is isomorphic to $\scrL(\Lambda_{Q_t})_\infty$. 
  \item[$\mathrm{(b)}$] The map $|\mu,Q_{t}\rangle_\infty\longmapsto [E_{\varpi_t(\mu)}]$ induces an isomorphism
  of $\frakg_\infty$-modules $\bfF({Q_t})_{\infty} \simto [\scrU_{K,t}].$
  \qed
 \end{itemize}
\end{theorem}

\medskip

\subsection{The $\frakg_e$-representation on $\scrU_\k$.}
\label{sec:g-e}\hfill\\

By Proposition \ref{prop:stable}, 
the representation datum $(E,F,X,T)$ on $\k G\mod$ induces a representation datum
on $\scrU_\k$. Since the abelian category $\scrU_\k$ is not semisimple,
to extend the representation datum to a categorical $\frakg$-representation
one needs to prove
that weight spaces of $\scrU_\k$ are sums of blocks. This will be done combinatorially by studying a
representation of a bigger Lie algebra $\frakg_{e,\circ}$.

\smallskip

\subsubsection{The Lie algebras $\frakg_e$ and $\frakg_{e,\circ}$}
\label{subsec:g-e}
The eigenvalues of $X$ on $E$ and
$F$ are all powers of $-q$. If we denote again by ${q}$ the image of $q$ under the canonical map 
$\scrO \twoheadrightarrow \k$, then
the eigenvalues of $X$ on $\k E$ and $\k F$ 
belong to the finite set $(-{q})^\mathbb{Z} \subset \k^\times$. This set has exactly $e$ elements, where $e$ is the order of $-{q}$ in $\k^\times$. 

\begin{definition}
We define $\scrI_e$ to be the subset $(-{q})^\bbZ$ of $\k^\times$. We denote by 
$\scrI_e$ and $\scrI_{e,\circ}$ the finite quivers $\scrI_e({q}^2)$ and
$\scrI_e(-{q})$. 
\end{definition}

The quivers $\scrI_e$ and $\scrI_{e,\circ}$ have the same set of vertices, but the arrows in 
$\scrI_e$ are the composition of two consecutive arrows in $\scrI_{e,\circ}$. The quiver 
$\scrI_{e,\circ}$ is cyclic, whereas the quiver $\scrI_{e}$ is cyclic if $e$ is odd,
and is a union of two cyclic quivers if $e$ is even. Therefore the corresponding
Kac-Moody algebras are isomorphic to $\widehat{\fraks\frakl}_e$ or
$(\widehat{\fraks\frakl}_{e/2})^{\oplus2}$.

\smallskip

To avoid cumbersome notation, we will write  $(\bullet)_{e,\circ} = (\bullet)_{\scrI_{e,\circ}}$
and $(\bullet)_{e} = (\bullet)_{\scrI_{e}}$. 
We must introduce Lie algebras $\frakg_{e,\circ}$ and $\frakg_{e}$ such that 
$\frakg_{e,\circ}' = [\frakg_{e,\circ},\frakg_{e,\circ}]$ and $\frakg_{e}' = 
[\frakg_{e},\frakg_{e}]$. 
The Chevalley generators of $\frakg_{e,\circ}'$ and $\frakg_{e}' $ are $e_{i,\circ}, f_{i,\circ}$ and 
$e_i,f_i$ respectively, for $i \in (-q)^\bbZ$. There is a Lie algebra homomorphism $\kappa : \frakg_e' \longrightarrow
\frakg_{e,\circ}'$ defined by 
 $$\kappa(e_i) = [e_{-qi,\circ},e_{i,\circ}] \quad \text{and} \quad
    \kappa(f_i)= [f_{-qi,\circ},f_{i,\circ}]. $$
It restricts to a map between the coroot lattices sending $\alpha_i^\vee$ to $\alpha_{i,\circ}^\vee
+\alpha_{-qi,\circ}^\vee$

\smallskip

We denote by $\frakg_{e,\circ}$ the Kac-Moody algebra associated with the lattices
$\X_{e,\circ} = \P_{e,\circ} \oplus \bbZ\delta_\circ$ and $\X_{e,\circ}^\vee = 
\Q_{e,\circ}^\vee \oplus \bbZ \partial_\circ$, where $\delta_\circ = \sum \alpha_{i,\circ}$, $\partial_\circ = \Lambda_{1,\circ}^\vee$
and the pairing $\X_{e,\circ}^\vee
\times \X_{e,\circ} \longrightarrow \bbZ$ is given by
 $$ \langle \alpha_{j,\circ}^\vee, \Lambda_{i,\circ} \rangle_{e,\circ} = \delta_{ij}, \quad
 \langle \partial_\circ,\Lambda_{i,\circ} \rangle_{e,\circ} = \langle \alpha_{j,\circ}^\vee,
 \delta_\circ \rangle_{e,\circ} = 0, \quad
  \langle \partial_\circ,\delta_\circ \rangle_{e,\circ} = 1.$$
Then $\frakg_{e,\circ}$, $\frakg_{e,\circ}'$ are isomorphic to
$\widehat{\fraks\frakl}_e$, $\widetilde{\fraks\frakl}_e$ respectively. 

\smallskip

Let  $\widetilde\frakg_e$ be the usual Kac-Moody algebra associated with $\scrI_e$.
Its derived Lie subalgebra is equal to $\frakg'_e$.
Let $\widetilde \X_e$ and $\widetilde \X_e^\vee$ be the lattices corresponding to $\widetilde\frakg_e$.
The map
$\kappa : \frakg_e' \longrightarrow
\frakg_{e,\circ}'$ may not extend to a morphism of Lie algebras $\widetilde\frakg_e \longrightarrow \frakg_{e,\circ}$.
For this reason we will define $\frakg_e$ to be  the Lie subalgebra $\frakg'_e\oplus\bbC\partial$ of $\widetilde\frakg_e,$
where $\partial$ is the element given by $\partial=\Lambda_1^\vee+\Lambda_{-q^{-1}}^\vee$.
We can view $\frakg_e$ as the Kac-Moody algebra associated with the lattice 
$\X_{e}^\vee =  \Q_{e}^\vee \oplus \bbZ \partial$ above and a lattice $\X_{e} = \P_{e} \oplus \bbZ\delta$ that we define
below case-by-case.

\smallskip

If $e$ is odd, then $\scrI_e$ is a cyclic
quiver and $\frakg_e=\widetilde\frakg_e$ is isomorphic to $\widehat{\fraks\frakl}_e$. 
Let $\alpha_1$ be the affine root, then we have $\widetilde \X_e = \P_e \oplus \bbZ \tilde\delta$ and 
$\widetilde \X_e^\vee = \Q_e^\vee \oplus \bbZ \widetilde \partial$
with $\tilde\delta = \sum \alpha_i$ and $\widetilde\partial = \Lambda_1^\vee$. 
Then, we set $\X_e = \P_{e} \oplus \bbZ\delta$ and $\delta=\tilde\delta/2$.

\smallskip

If $e$ is even, then $\scrI_e$ is the disjoint union
of two cyclic quivers and $\tilde\frakg_e$ is isomorphic
to $(\widehat{\fraks\frakl}_{e/2})^{\oplus2}$. 
Let $\alpha_1$ and $\alpha_{-q^{-1}}$
be the affine roots, then we have $\widetilde \X_e = \P_e \oplus
 \bbZ \delta_1 \oplus \bbZ\delta_2$ and $\widetilde \X_e^\vee = \Q_e^\vee \oplus \bbZ \partial_1
\oplus \bbZ \partial_2$ with $\delta_1 = \sum_{j \text{ odd}} \alpha_{-q^{j}}$, 
$\delta_2 = \sum_{j \text{ even}} \alpha_{q^{j}}$, $\partial_1 = \Lambda_{-q^{-1}}^\vee$ and
$\partial_2 = \Lambda_{1}^\vee$. Then,  we set $\X_e = \widetilde \X_e / (\delta_1-\delta_2)$ and $\delta=(\delta_1+\delta_2)/2$.

\smallskip

In both cases the perfect pairing $\X_e^\vee\times\X_e\longrightarrow\bbZ$ 
is induced in the obvious way by the pairing $\widetilde\X_e^\vee\times\widetilde\X_e\longrightarrow\bbZ$. 

\begin{lemma}\label{lem:kappa}
There is a well-defined morphism of Lie algebras
$\frakg_e \longrightarrow \frakg_{e,\circ}$ which extends $\kappa$ and whose restriction to $\X_{e}^\vee$ 
is given by  
$ \kappa(\alpha_i^\vee) = \alpha_{i,\circ}^\vee + \alpha_{-qi,\circ}^\vee$ and
 $\kappa(\partial)= \partial_\circ.$
The restriction $\kappa : \X^\vee_e\longrightarrow\X^\vee_{e,\circ}$ has an adjoint $\kappa^* : \X_{e,\circ} \longrightarrow \X_{e}$ such that
$ \kappa^*(\Lambda_{i,\circ})\equiv\Lambda_{i} + \Lambda_{-q^{-1}i}\ \rm{mod}\ \delta$ and
 $\kappa^*(\delta_\circ)= \delta.$
  \qed
\end{lemma}

\smallskip

\subsubsection{The $\frakg'_e$-action on $[\scrU_\k]$}
\label{subsec:derivedaction-uk}
The quotient map $ \scrO \twoheadrightarrow \k$ induces a morphism of quivers
$\sp : \scrI_\infty \longrightarrow \scrI_e$ and a surjective morphism of abelian groups
$\sp : \P_\infty \twoheadrightarrow \P_e$  sending $\Lambda_i$ to $\Lambda_{\sp(i)}$.
In addition any integrable representation $V$ of $\frakg_\infty$ can be ``restricted''
to an integrable representation of the derived algebra $\frakg'_e$, where
$e_i \in \frakg'_e$ (resp. $f_i \in \frakg'_e$) act as $\sum_{\sp(j)=i} e_j$ (resp.
$\sum_{\sp(j)=i} f_j$). From the definition of the action of $\frakg_\infty$ and 
$\frakg_e$ on Fock spaces, see \eqref{eq:EF}, we deduce that the map
$|\mu,Q_{t}\rangle_\infty\mapsto|\mu,Q_t\rangle_e$ induces the following isomorphism
of $\frakg'_e$-modules
 \begin{equation}\label{eq:fockinftyande} 
 \sp \, : \, \mathrm{Res}_{\frakg'_e}^{\frakg_\infty} \, \bfF(Q_t)_\infty
 \, \mathop{\longrightarrow}\limits^\sim\, \bfF(Q_t)_e.\end{equation}
Under the decomposition map, this isomorphism endows $[\scrU_\k]$ with
a structure of $\frakg'_e$-module which is compatible with the one coming
from the representation datum.

\begin{proposition}\label{prop:lreduction}
For each $i \in \scrI_e$, let $\k E_i$ (resp.~$\k F_i$)
be the generalized $i$-eigenspace of $X$ on $\k E$ (resp.~$\k F$). 
  \begin{itemize}[leftmargin=8mm]
    \item[$\mathrm{(a)}$] $[\k E_i], [\k F_i]$ endow $[\scrU_\k]$ with a structure
    of $\frakg'_e$-module.
    \item[$\mathrm{(b)}$] $d_\scrU : [\scrU_K] \longrightarrow [\scrU_k]$
     induces a $\frakg'_e$-module isomorphism 
     $ \mathrm{Res}_{\frakg'_e}^{\frakg_\infty} \, [\scrU_K]\simto [\scrU_k].$
      \item[$\mathrm{(c)}$]
     $|\mu, Q_t\rangle_e \longmapsto [V_{\varpi_t(\mu)}]$ induces a
$\frakg'_e$-module isomorphism
$\bigoplus_{t\in\bbN}\bfF({Q_t})_{e} \, \mathop{\longrightarrow}\limits^\sim \, [\scrU_\k].$
     \qed
  \end{itemize}
\end{proposition}

\smallskip

\subsubsection{The $\frakg_e$-action on $[\scrU_\k]$}\label{subsec:action-uk}
We now define an action of $\frakg_e$ on $[\scrU_\k]$ by extending the action from $\frakg'_e$
to $\frakg_e$ on $\bigoplus_{t\in\bbN}\bfF({Q_t})_{e}$. This amounts to extending the grading from $\P_e$ to $\X_e = \P_e \oplus 
\bbZ \delta$.

\smallskip

Any integrable $\frakg_{e,\circ}$-representation 
(resp. $\frakg'_{e,\circ}$-representation) can be ``restricted'' to an 
integrable $\frakg_e$-representation (resp. $\frakg'_{e}$-representation)
through $\kappa$. We denote by $\mathrm{Res}_{\frakg_{e}}^{\frakg_{e,\circ}}$ and 
$\mathrm{Res}_{\frakg'_{e}}^{\frakg'_{e,\circ}}$
the corresponding operations. 

\smallskip

Let $a\,:\,\scrP\to\bbN$ be Lusztig's $a$-function, see \cite[4.4.2]{Lu84}. 
Recall that $a(\lambda)=\sum_i(i-1)\lambda_i$ if $\lambda=(\lambda_1,\lambda_2,\dots)$.

\begin{lemma}\label{lem:isofocks}
The assignment $|\mu,Q_t\rangle_e \longmapsto (-1)^{a(\varpi_t(\mu))}
|\varpi_t(\mu),1\rangle_{e,\circ}$
induces an isomorphism of $\frakg'_{e}$-modules
  $ \bigoplus_{t \in \bbN} \bfF(Q_t)_e \ \mathop{\longrightarrow}
  \limits^\sim\ \mathrm{Res}_{\frakg'_{e}}^{\frakg'_{e,\circ}}
  \, \bfF(1)_{e,\circ}.$
\qed
\end{lemma}

Now, consider the Fock space $\bfF(1)_{e,\circ}$ as a charged Fock space for the charge $s = 0$.
This endows $\bfF(1)_{e,\circ}$ with an integrable representation of $\frakg_{e,\circ}$ as in
\S \ref{subsec:xgrading}. 
Consequently, 
we can endow $\bigoplus_{t\in\bbN} \bfF(Q_t)_e$ and $[\scrU_\k]$ with
integrable representations of $\frakg_e$ such that there are $\frakg_e$-module isomorphisms
\begin{align}\label{eq:xgrading}
\begin{split}
&\bigoplus_{t\in\bbN}\bfF(Q_t)_e\ \mathop{\longrightarrow}
  \limits^\sim\  [\scrU_\k] \ \mathop{\longrightarrow}
  \limits^\sim\ \mathrm{Res}_{\frakg_{e}}^{\frakg_{e,\circ}}
  \, \bfF(1)_{e,\circ}\\
&|\mu,Q_t\rangle_e\mapsto [V_{\varpi_t(\mu)}] \mapsto (-1)^{a(\varpi_t(\mu))}|\varpi_t(\mu),1\rangle_{e,\circ}.
\end{split}
\end{align}

\smallskip

\subsubsection{The $\frakg_e$-action on $\scrU_\k$}
The following lemma holds.

\begin{lemma}\label{lem:sameblocksameweight}
Let $\lambda$ and $\mu$ be partitions of $n$. If $V_\mu$,
$V_\lambda$ belong to the same block of $\scrU_\k$ then $[V_\lambda]$, 
$[V_\mu]$ have the same weight for the action of $\frakg_e$.
\qed
\end{lemma}

We deduce that the classes of the simple modules $[D_\lambda]$ are
also weight vectors. Given $\omega \in \X_e$, we define
$\scrU_{\k,\omega}$ to be the Serre subcategory of $\scrU_\k$ generated by the
simple modules $D_\lambda$ such that the element $[D_\lambda]$ in $[\scrU_\k]$ has weight $\omega$. We have
$ \scrU_\k = \bigoplus_{\omega \in \X_e} \scrU_{\k,\omega}.$
We can now formulate the following.

\begin{theorem}\label{thm:charl}
The representation datum associated with Harish-Chandra induction and restriction and the decomposition $\scrU_\k = \bigoplus_{\omega \in \X_e} \scrU_{\k,\omega}$
yield a categorical representation of $\frakg_e$ on $\scrU_\k$. Furthermore, the map
$ |\mu,Q_t\rangle_e \longmapsto [V_{\varpi_t(\mu)}]$ induces a $\frakg'_e$-module
isomorphism $\bigoplus_{t \in \bbN} \bfF(Q_t)_e \simto [\scrU_\k]$.
\qed
\end{theorem}

\medskip

\subsection{Derived equivalences of blocks of $\scrU_\k$}\label{sec:broue}\hfill\\

Recall that $d$, $e$ and $f$ denote respectively the order of $q^2$, $-q$ and $q$ modulo 
$\ell$.

\smallskip

\subsubsection{Characterization of the blocks of $\scrU_\k$}
\label{subsec:block=weight}
 Here we investigate
which block can occur in a given weight space. 


Assume $e$ is even. Given $t\in\bbN$ and $\omega \in \X_e$, we define
$\scrU_{\k,t}$ (resp. $\scrU_{\k,t,\omega}$) to be the Serre
subcategory of $\scrU_{\k}$ generated by simple modules $D_\lambda$
where  $\lambda$ has $2$-core $\Delta_t $ (resp. with 
in addition $\omega_\lambda = \omega$). With $e$ being even, any pair of
partitions with the same $e$-core have the same $2$-core, therefore
$\scrU_{\k,t}$ is a direct summand of $\scrU_\k$.

\begin{proposition} \label{prop:main-blocks}
Let $\omega\in\X_e$.
\begin{itemize}[leftmargin=8mm]
  \item[$\mathrm{(a)}$] If $e$ is odd, then the category $\scrU_{\k,\omega}$ is an indecomposable
  summand of $\scrU_\k$.
  \item[$\mathrm{(b)}$] If $e$ is even, then the category $\scrU_{\k,t,\omega}$ is an 
  indecomposable summand of $\scrU_{\k,t}$ and $\scrU_\k$ for all $t\in\bbN$.
  \qed
\end{itemize}
\end{proposition}


\smallskip

\subsubsection{Derived equivalences of blocks of $\scrU_\k$} \label{subsec:broue}
It is not difficult to compute the orbits of the affine Weyl group on the weight spaces of
$[\scrU_\k]$. Hence, we can apply Proposition
\ref{prop:main-blocks} and the work of Chuang and Rouquier \cite{CR} to produce
derived equivalences between blocks of $\scrU_\k$ in the same orbit.
Using the results of Livesey \cite{Li12} on the structure of good blocks for linear primes,
we deduce that Brou\'e's abelian defect group conjecture holds for unipotent blocks
when $e$ is even.

\begin{theorem}\label{thm:broue}
Assume $e$ is even. Let $B$ be a unipotent block of $G_n$ over $\k$ or $\scrO$,
and $D$ be a defect group of $B$. If the group $D$ is abelian, then the block $B$ is derived equivalent to its Brauer
 correspondent in $N_{G_n}(D)$. \qed
\end{theorem}

\medskip

\subsection{The crystals of $\scrU_K$ and $\scrU_\k$}\label{sec:isocrystals}\hfill\\

Now, we want to compare the crystals of the categorical representations
$\scrU_K$ and $\scrU_\k$ (which are related to Harish-Chandra induction and restriction)
with the crystals of the Fock spaces related to $[\scrU_K]$ and $[\scrU_\k]$. 
This solves the main conjecture of Gerber-Hiss-Jacon \cite{GHJ} and gives a combinatorial way to
compute the (weak) Harish-Chandra branching graph and the Hecke algebras associated
to the weakly cuspidal unipotent modules.

\smallskip

\subsubsection{Crystals and Harish-Chandra series\label{ssec:crystals-and-series}}
Recall that to any categorical representation one can associate a perfect basis, and
hence an abstract crystal.
In the previous sections we constructed a categorical action on the categories
of unipotent representations over $K$ (denoted by $\scrU_K$) and over $\k$ (denoted
by $\scrU_k$). From these two categorical representations we get the following abstract crystals
\begin{itemize}[leftmargin=8mm]
  \item[(a)] $B(\scrU_K)=\big(\Irr(\scrU_K),\widetilde E_i,\widetilde F_i\big)$
  where $\Irr(\scrU_K)=\{[E_\lambda]\,|\,\lambda\in\scrP\}$,
  
  \item[(b)] $B(\scrU_\k)=\big(\Irr(\scrU_\k),\widetilde E_i,\widetilde F_i\big)$
  where $\Irr(\scrU_\k)=\{[D_\lambda]\,|\,\lambda\in\scrP\}.$
\end{itemize}
The (uncolored) crystal graph associated with $B(\scrU_R)$ coincides with the weak
  Harish-Chandra branching graph, and its connected components with the weak Harish-Chandra series.

\begin{proposition}\label{prop:wHC} Let $R = K$ or $\k$, and 
  $\scrI = (-q)^{\bbZ} \subset R^\times$. Let $D,M,N \in \Irr(\scrU_R)$.
  \begin{itemize}[leftmargin=8mm]
    \item[$\mathrm{(a)}$]
   $D$ is weakly cuspidal if and only if $\widetilde E_i (D) =0$
    for all $i \in \scrI$.
    \item[$\mathrm{(b)}$]  $N$ appears in $\head(F(M))$ if and only if
    there exists $i \in \scrI$ such that $N \simeq \widetilde F_i M$.
    \item[$\mathrm{(c)}$] If $D$ is weakly cuspidal, then 
    $$\wIrr(RG,D)= \{\widetilde F_{i_1}\cdots\widetilde F_{i_m}(D)\,\mid\,m\in\bbN,\,
    i_1,\dots,i_m\in \scrI\}.$$
  \end{itemize}
 
  \qed
  \end{proposition}

Note that by \cite{L77}, the ordinary Harish-Chandra series and weak Harish-Chandra series
on unipotent modules coincide when $R = K$.

\smallskip

\subsubsection{Comparison of the crystals}\label{sec:crystal-unitary}
In this section we will assume that $e$ is odd.
To any charged Fock space one can associate an abstract crystal, see \S \ref{subsec:crystalfock}. 
We now show how to choose the charge for each Fock space $\bfF(Q_t)_e$ so that the crystal
will coincide with the Harish-Chandra branching graph. 

\smallskip

Define 
$s_t=(s_1,s_2)=-{1\over 2}(e+1,0)+\sigma_t$ where $\sigma_t$ is as in \eqref{sigma}.
We have $Q_p=q^{2s_p}$ for each $p=1,2$.  We denote by $B(s_t)_e$ the corresponding
abstract crystal of $\bfF(Q_t)_e$, with the canonical labeling
$B(s_t)_e=\{b(\mu,s_t)\,|\,\mu\in\scrP^2,\,t\in\bbN\}.$
Finally, we set $B_e=\bigsqcup_{\,t\in\bbN}B(s_t)_e$.

\begin{theorem}\label{thm:wHC}
 The map $b(\mu,s_{t})\mapsto [D_{\varpi_t(\mu)}]$ is a
crystal isomorphism $B_e\simto B(\scrU_\k)$. \qed
\end{theorem}

\smallskip


We deduce the following corollaries.

\begin{corollary}\label{cor:1}
 The modules $D_\lambda$ and $D_\nu$ lie in the same weak
Harish-Chandra series if and only if the corresponding vertices of the abstract crystal $B_e$ belong to the same connected component.
In particular, if this holds then $\lambda$ and $\nu$ have the same $2$-core. \qed
\end{corollary}

\begin{corollary}\label{cor:2} 
Let $\lambda$ be a partition of $r>0$ such that $D_\lambda$ is weakly cuspidal.
\begin{itemize}[leftmargin=8mm]
\item[$\mathrm{(a)}$] The $e$-core of $\lambda$ is a 2-core $\Delta_t = (t,t-1,\ldots,1)$,
\item[$\mathrm{(b)}$] the weight of the class $[D_\lambda]$ with respect to the $\frakg'_e$-action on $[\scrU_\k]$ is $\Lambda_{Q_t}$,
\item[$\mathrm{(c)}$] for each $m \geqslant 0$ and $n = r+2m$,
the map $\phi_{\k,m}\, :\, \bfH^{q^2}_{\k,m}\longrightarrow\End_{\k G_n}(F^m)^\op$ factors through an algebra isomorphism
  $\bfH^{Q_t;\,q^2}_{\k,m}\mathop{\longrightarrow}\limits^\sim \scrH(\k G_n,D_\lambda).$
  \qed
\end{itemize}
\end{corollary}

\medskip

\section{The representation of the Heisenberg algebra on $\scrU_\k$}\label{sec:heisenberg}
Every cuspidal $\k G_n$-module is weakly cuspidal. Therefore, every Harish-Chandra
series in $\Irr(\scrU_\k)$ is partitioned into weak Harish-Chandra series. Proposition
\ref{prop:wHC} and Theorem \ref{thm:wHC} yield a complete (combinatorial) description of
the partition of $\Irr(\scrU_\k)$ into weak Harish-Chandra series using
the decomposition of $[\scrU_\k]$ for the action of $\frakg_e$. 
In this section we generalize this result to usual Harish-Chandra series by adding 
the action of some Heisenberg algebra on $[\scrU_\k]$.
Throughout this section we will assume that $e$, the order of $-q$ modulo $\ell$ 
is odd, and hence equal to $d$, the order of $q^2$ modulo~$\ell$.

\medskip

\subsection{The Heisenberg action on $[\scrU_\k]$}\label{sec:Heisenberg}
\subsubsection{The Heisenberg algebra}
Set $\pmb\Lambda=\bigoplus_{m\in\bbN}\bbC\Irr(K\frakS_m)$ where
$\frakS_0$ is the trivial group.
We equip the vector space $\pmb\Lambda$ with the bilinear form 
$\langle\bullet,\bullet\rangle_{\pmb\Lambda}=\bigoplus_{m\in\bbN}\langle\bullet,\bullet\rangle_{\frakS_m}$
and the vector space
$\pmb\Lambda\otimes\pmb\Lambda$ with the tensor square of $\langle\bullet,\bullet\rangle_{\pmb\Lambda}$.
The induction and restriction yield a pair of linear maps $\Ind_{n,m}^{n+m}$ and
$\Res_{n,m}^{n+m}$ between the $\bbC$-vector spaces $\pmb\Lambda$ and
$\pmb\Lambda\otimes\pmb\Lambda$ which are adjoint with respect to these bilinear forms.

\smallskip

Fix a positive integer $a$. A level $a$ representation of the Heisenberg algebra, or an $\frakH_a$-module,
is a $\bbC$-vector space $V$ with a family of 
endomorphisms $b_n,$ $ b^*_n$ labelled by integers $n\in a\bbZ_{>0}$ satisfying the relations
$[b_n,b_m]=[b^*_n, b^*_m]=0$ and $[b_n,b^*_m]=-n\delta_{n,m}$.
For each partition $\lambda$ we write $a\lambda=(a\lambda_1,a\lambda_2,\dots)$ and
$b_{a\lambda}=\prod_ib_{a\lambda_i}$, 
$b^*_{a\lambda}=\prod_ib^*_{a\lambda_i}$.

\smallskip

If $\lambda$ is a partition of $m$ and
$w\in\frakS_m$ is of \emph{cycle-type} $\rho(w)=\lambda$, we denote by $c_\lambda=c_w$ the conjugacy class of $w$. 
The characteristic function of $c_\lambda$ is a class function on $\frakS_m$, 
and as such it can be viewed as an element in $\pmb\Lambda$. 
We abbreviate $c_m=c_{(m)}$ for the conjugacy class of the $m$-cycles.
Let $z_{\lambda}$ be the number of elements of the centralizer of $w$ in $\frakS_m$.
We have $\langle c_\lambda,c_\mu \rangle = z_\lambda^{-1} \delta_{\lambda,\mu}$.

\smallskip

There is a unique representation of $\frakH_1$ on $\pmb\Lambda$ such that
for every $\phi\in\bbC\Irr(\frakS_r)$ and $\psi\in\bbC\Irr(\frakS_{n})$ with $n=r+m$ we have
\begin{align}\label{Heisenberg}
b_{\lambda}(\phi)=\Ind_{r,m}^{n}(\phi\otimes z_{\lambda}c_{\lambda}),\quad
b^*_{\lambda}(\psi)=\langle\Res_{r,m}^{n}(\psi)\,,\,z_{\lambda}c_{\lambda}\rangle_{\pmb\Lambda}.
\end{align}

\smallskip

\subsubsection{The Heisenberg action on $\bfF(Q)$}
Recall that $e$ is the order of $-q$ modulo $\ell$. We defined in \S \ref{subsec:g-e}
a Kac-Moody algebra $\frakg_e$ corresponding to a quiver $I_e$. 
From now one we will assume that $e$ is odd and we will abbreviate
$\frakg = \frakg_e$. Then $\frakg$ is isomorphic to $\widehat{\fraks\frakl}_e$.
We can now define the representation of $\frakH_{le}$ on $\bfF(Q)$, where
$Q=(Q_1,\ldots,Q_l) \in I_e^l$.
The Fock space $\bfF(Q)$ admits a representation of $\frakH_{l}$ 
which is identified with the $l$-th tensor power of the representation of $\frakH_1$ on $\pmb\Lambda$
in \eqref{Heisenberg} under the $\bbC$-linear isomorphism
$\bfF(Q)\simto\pmb\Lambda^{\otimes l}$ such that $|\mu,Q\rangle\mapsto\phi_{\mu^1}\otimes\cdots\otimes\phi_{\mu^l}.$
Then, by \eqref{eq:EF}, the representation of $\frakg'$ on  $\bfF(Q)$ commutes with the action of the subalgebra $\frakH_{le}$ 
of $\frakH_l$. We will apply this construction to $l=1$ or $2$.

\smallskip

\subsubsection{The Heisenberg action on $[\scrU_\k]$}\label{heisenberg-grot}

The previous section yields a representation of $\frakH_{2e}$ in $\bfF(Q_t)_e$ for each $t\in\bbN$ and a representation of
$\frakH_e$ in $\bfF(1)_{e,\circ}$.
We want to compare them.

\smallskip

To do that, consider the $\bbC$-linear isomorphism $\bigoplus_{t\in\bbN}\bfF(Q_t)_e\simto\bfF(1)_{e,\circ}$
studied in Lemma \ref{lem:isofocks}, which is given by $|\mu,Q_t\rangle_e\mapsto(-1)^{a(\lambda)}|\lambda,1\rangle_{e,\circ}$ with $\lambda=\varpi_t(\mu).$

\begin{lemma}\label{lem:isomFockHeisenberg} The $\frakg_e$-linear isomorphism 
in Lemma \ref{lem:isofocks} is an isomorphism of $\frakH_{2e}$-modules
$ \bigoplus_{t\in\bbN}\bfF(Q_t)_e\mathop{\longrightarrow}\limits^\sim \mathrm{Res}_{\frakH_{2e}}^{\frakH_e} 
\bfF(1)_{e,\circ}.$
\end{lemma}

\begin{proof}
Fix an integer $t\geqslant 0$. 
We must check that the linear map
$$\bfF(Q_t)_e\simeq\pmb\Lambda\otimes\pmb\Lambda\to\pmb\Lambda$$ such that
$$|\mu,Q_t\rangle_e\mapsto\phi_{\mu^1}\otimes\phi_{\mu^2}\mapsto(-1)^{a(\lambda)}\phi_{\lambda},
\quad\forall\mu\in\scrP^2$$
intertwines the operator $b_n\otimes 1+1\otimes b_n$ on $\pmb\Lambda\otimes\pmb\Lambda$ 
with the operator $b_{2n}$ on $\pmb\Lambda$ if $n\in e\bbZ_{>0}$.
Here, we have set $\lambda=\varpi_t(\mu)$.

\smallskip

To do this, note that we have
\begin{align}\label{form2bis}
b_{2n}(\phi_\lambda)=\sum_{\nu}(-1)^{N(\lambda,x,2n)}\,\phi_\nu,
\end{align}
where the sum runs over all partitions $\nu$ obtained by adding a $2n$-hook $(x,x+2n)$ to
the charged partition $(\lambda,0)$
and $N(\lambda,x,2n)$ is the number of elements in $\beta_0(\lambda)\cap(x,x+2n)$.
Recall that $\sigma_t=(\sigma_1,\sigma_2)$. 
Since
$$\beta_0(\lambda)=(-1+2\beta_{\sigma_1}(\mu^1))\sqcup(2\beta_{\sigma_2}(\mu^2)),$$
we have
$$N(\lambda,x,2n)=\begin{cases}
N(\mu^1,y-\sigma_1,n)+|(2\beta_{\sigma_2}(\mu^2))\cap(x,x+2n)|&\text{if}\ x=-1+2y,\\
|(-1+2\beta_{\sigma_1}(\mu^1))\cap(x,x+2n)|+N(\mu^2,y-\sigma_2,n)&\text{if}\ x=2y.
\end{cases}$$
Furthermore, adding a $2n$-hook $(x,x+2n)$ to $\lambda$ is equivalent to adding an
$n$-hook to $\mu^1$ if $x$ is odd or to $\mu^2$ if $x$ is even. 
So we are reduced to check that if $\lambda$, $\nu$, $x$ are as in the sum \eqref{form2bis}, then
$a(\nu)-a(\lambda)$ has the same parity as the integer $b(\lambda,x)$ given by
\begin{itemize}
\item $|(2\beta_{\sigma_2}(\mu^2))\cap(x,x+2n)|$ if $x$ is odd,
\item $|(-1+2\beta_{\sigma_1}(\mu^1))\cap(x,x+2n)|$ if $x$ is even.
\end{itemize}
Now a partition $\nu$ obtained by adding a $2n$-hook to $\lambda$ can also be 
obtained by adding $n$ successive $2$-hooks to $\lambda$. Therefore it is enough to check
that the parity of $a(\nu)-a(\lambda)$ and $b(\lambda,x)$ are the same when $n =1$. 
Since $a(\lambda)=\sum_i(i-1)\lambda_i$, the integer $a(\nu)-a(\lambda)$ 
is even if $\nu\setminus\lambda$ is a vertical $2$-hook
(increasing $\lambda_i$ by $2$ for some $i$) or odd if $\nu\setminus\lambda$ is a horizontal 2-hook 
(increasing $\lambda_i$ and $\lambda_{i+1} = \lambda_i$ by $1$ for some $i$). 
On the other hand, we have $b(\lambda,x) = 0$ (resp. $b(\lambda,x)=1$) 
for a vertical $2$-hook (resp. for a horizontal $2$-hook).
\end{proof}

\smallskip

We equip $[\scrU_\k]$ with the unique representation of $\frakH_{2e}$ such that the isomorphism
$\bigoplus_{t\in\bbN}\bfF(Q_t)_e\simto[\scrU_\k]$ in Theorem \ref{thm:charl} is $\frakH_{2e}$-equivariant.
From Lemmas
\ref{lem:isofocks} and \ref{lem:isomFockHeisenberg} we deduce that  there is a $\frakg_e\times\frakH_{2e}$-equivariant isomorphism 
given by
\begin{align}\label{isomB}
[\scrU_\k]\simto\bfF(1)_{e,\circ},\quad [V_\lambda]\mapsto(-1)^{a(\lambda)}|\lambda,1\rangle_{e,\circ}.
\end{align}

\smallskip

\begin{remark}\label{rem:epsilon}
Set $A(\lambda)=m(m-1)/2-a(\lambda^*)$ with $m=|\lambda|$.
By \cite[lem.~4.2]{DM} we have $(-1)^{A(\lambda)}=\varepsilon_\lambda$ where $\varepsilon_\lambda$ was defined in \S\ref{subsec:uK}.
From this we can check that
$$\lambda_{[2]}=\nu_{[2]}\Rightarrow(-1)^{a(\lambda)}\varepsilon_\lambda=(-1)^{a(\nu)}\varepsilon_\nu.$$ 
Hence, the endomorphism of $[\scrU_\k]$ such that
$[V_\lambda]\mapsto(-1)^{a(\lambda)}\varepsilon_\lambda\,[V_\lambda]$ commutes with the action of $\frakg_e$ and $\frakH_{2e}$. 
This shows that we could have used the sign $\varepsilon_\lambda$ instead of $(-1)^{a(\lambda)}$ in the isomorphism \eqref{isomB}.
 \end{remark}

\medskip


\subsection{The modular Harish-Chandra series of $GL_m$}\hfill\\

Recall that $GL_n := \mathrm{GL}_n(q^2)$ denotes the finite general linear group over a finite field with $q^2$ elements. 

\subsubsection{The unipotent modules of $GL_n$}

The set of unipotent characters of $K GL_n$ is
\begin{align}\label{L}
\Unip(KGL_n)=\{L_\lambda\,|\,\lambda\in\scrP_n\},
\end{align}
where
\begin{align}\label{L2}
L_\lambda= |\frakS_n|^{-1}\sum_{w \in \frakS_n}\phi_\lambda(w) R_{\bfT_w}^{\bfG\bfL_n}(1)
\end{align}
and $\bfT_w$ is an $F_{q^2}$-stable maximal torus of $\bfG\bfL_n$ of cycle-type $\rho(w)$,
with the convention that $\bfT_1=\bfT$ is the split torus of diagonal matrices. 
There is a unique labeling of the simple unipotent modules of $\k GL_n$ given by
$$\Unip(\k GL_n)=\{S_\lambda\,|\,\lambda\in\scrP_n\}$$
where $$d_{\scrO GL_n}([L_\lambda])=[S_\lambda]\ \text{modulo}\ \bigoplus_{\mu>\lambda}\bbZ\,[S_\mu].$$
As in the case of unitary groups (see  Proposition \ref{prop:liftingstandards}), there exists a unique indecomposable unipotent $\scrO GL_n$-lattice $\widetilde W_\lambda$ with character
$L_\lambda$ and such that $W_\lambda:= k \widetilde W_\lambda$ is indecomposable with simple head $S_\lambda$.

\smallskip

\subsubsection{The modular Steinberg character and Harish-Chandra series}
First, let us introduce the following notation.
Given non-negative integers $m_j$ with $j\geqslant -1$, we will say that a partition 
$\lambda$ of $m$ is \emph{cuspidal of type $(m_j)$} if $\lambda$ is of the following form
$$\lambda=(1^{(m_{-1})},e^{(m_0)},(e\ell)^{(m_1)},(e\ell^2)^{(m_2)},\dots)$$
in exponential notation.

\smallskip

For a partition $\lambda$ as above, we consider the Levi subgroup of $\bfG\bfL_m$ given by
$$\bfG\bfL_\lambda=(\bfG\bfL_1)^{m_{-1}}\times\prod_{j\geqslant 0}(\bfG\bfL_{e\ell^j})^{m_j}.$$
Recall that $GL_\lambda=\bfG\bfL_\lambda(q^2)$. 
The unipotent module $S_{(1^m)}$ is the head of $W_{(1^n)}$. Since  $W_{(1^n)}$ is a modular reduction of the Steinberg 
character $L_{(1^m)}$ of $KGL_m$, we call
$S_{(1^m)}$ the \emph{modular Steinberg representation} of $\k GL_m$ and we write $St_m=S_{(1^m)}$. 
Let $St_\lambda$ be the $\k GL_\lambda$-module given by
$$St_\lambda=(St_1)^{\otimes m_{-1}}\otimes\bigotimes_{j\geqslant 0}(St_{e\ell^j})^{\otimes m_j}.$$

\smallskip

The following result is due to Dipper, James and Dipper-Du \cite{DiDu}, see also \cite[sec.~19]{CE04}.

\begin{proposition}\label{prop:cusp-GL}\hfill
\begin{itemize}[leftmargin=8mm]
   \item[$\mathrm{(a)}$]
The cuspidal pairs of $\k GL_m$ are the $GL_m$-conjugates of the pairs $(GL_\lambda,St_\lambda)$ where
$\lambda$ runs over the set of all cuspidal partitions of $m$. 
 \item[$\mathrm{(c)}$] If $\lambda$ is cuspidal of type $(m_j)$, then we have a $\k$-algebra isomorphism
$$\scrH(\k GL_m,St_\lambda)\simeq\bfH^{1;\,q^2}_{\k,m_{-1}}\otimes\bigotimes_{j\geqslant 0}\bfH^{1;1}_{\k,m_j}.$$
 \qed
 \end{itemize}
\end{proposition}

Recall that here $\bfH^{1;\,u}_{\k,m}$ denotes the Hecke algebra of type $\frakS_m$ over $\k$
 with parameter $u$. In particular, we have $\bfH^{1;\,1}_{\k,m} = \k \frakS_m$.

\medskip

\subsection{The Heisenberg functors}\hfill\\

For $R=K$ or $\k$ we call the category of unipotent $RGL$-modules the category
$$\scrV_R=\bigoplus_{m\in\bbN}R GL_m\umod.$$
Recall that the vector space $\pmb\Lambda$ is equipped with the bilinear form 
$\langle\bullet,\bullet\rangle_{\pmb\Lambda}$. It induces a  bilinear form 
 $\langle\bullet,\bullet\rangle_{\scrV}$ on $[\scrV_\k]$ via the isomorphism of vector spaces
 \begin{align}\label{isomC0}
[\scrV_\k]\simto\pmb\Lambda,\quad [W_\lambda]\mapsto\phi_\lambda.
\end{align}

\smallskip

Now, let us come back to the category $\scrU_\k$ of unipotent modules of the unitary group.
Let $\langle\bullet,\bullet\rangle_\scrU$ be the unique bilinear form on $[\scrU_\k]$ such that the map \eqref{isomB} is an isometry.

\smallskip

Given $n,r,m>0$ such that $n=r+2m$, 
we write $L_{r,m} \simeq G_r \times GL_m$.
Given a unipotent $R GL_{m}$-module $X$ we consider the following functors
\begin{align*}
&B_X\,:\,RG_r\umod\to RG_{n}\umod,\ M\mapsto R_{L_{r,m}}^{G_{n}}(M\otimes X),\\
&B^*_X\,:\,RG_{n}\umod\to RG_{r}\umod,\ M\mapsto 
\Hom_{R  GL_m}\big(X,{}^*\!R_{L_{r,m}}^{G_{n}}(M)\big).
\end{align*}
The functor $B^*_X$ is right adjoint to $B_X$.
In the particular case where $m=1$ and $X =R$ is the trivial module, we recover the functors
$F$, $E$ defined in \S\ref{sec:rep-datum}.
The functor $B_X$ is exact, hence it induces a $\bbC$-linear endomorphism $[B_X]$ of $[\scrU_\k]$. However, the functor $B_X^*$ 
is not exact and the module $X$ may not have a finite projective dimension in $\k GL_m\umod$.
Hence, the right derived functor  of $B_X^*$
may not yield any  endomorphism of $[\scrU_\k]$.
Nevertheless we can avoid this problem by defining the $\bbC$-linear map
$$[B^*_X]=\langle\,[X]\,,\,[{}^*\!R_{L_{r,m}}^{G_{n}}](\bullet)\,\rangle_{\scrV} \, : 
[\k G_n\umod] \rightarrow [\k G_r\umod].$$
The Harish-Chandra induction relative to the subgroups $GL_{n_1}\times GL_{n_2}$ of $GL_n$ for $n_1+n_2=n$
yields a bifunctor $\bullet\odot\,\bullet\,:\,\scrV_R\otimes\scrV_R\to\scrV_R.$

\begin{proposition} \label{prop:B} Let $X, X'$ be objects in $\scrV_\k$. 
\begin{itemize}[leftmargin=8mm]
  \item[$\mathrm{(a)}$] 
There is an isomorphism of functors $B_{X}B_{X'}\simeq B_{X\odot X'}\simeq B_{X'}B_{X}$.
 \item[$\mathrm{(b)}$] 
There is an equality of linear maps $[B^*_{X}]\circ[B^*_{X'}] = [B^*_{X\odot X'}] = [B^*_{X'}]\circ[B^*_{X}].$
 \item[$\mathrm{(c)}$]  
 The linear maps $[B_X]$, $[B^*_X]$ are adjoint relatively to the bilinear form on $[\scrU_\k]$. 
  \item[$\mathrm{(d)}$]  
  If $M\in\scrU_\k$ is cuspidal, then $[B^*_{X}]([M])=0$.
\qed
\end{itemize}
\end{proposition}

\medskip

\subsection{The categorification of the Heisenberg action on $[\scrU_\k]$}\hfill\\

Recall that $e$ is odd. The representation of $\frakH_{2e}$ on $[\scrU_\k]$ in \S\ref{heisenberg-grot} 
is determined by linear operators
$b_n,$ $b^*_n$ with $n$ a positive multiple of $2e$ satisfying the relations
$[b_n,b_m]=[b^*_n, b^*_m]=0$ and $[b_n,b^*_m]=-n\,\delta_{n,m}\,\Id$. 
In this section we will re-interpret \emph{some} of these operators in representation theoretic terms using the category $\scrU_\k$.

\smallskip

Now, fix a positive integer $m$ and a cuspidal partition $\lambda$ of $m$. Consider the module
$$X_\lambda=R^{GL_m}_{GL_\lambda}(St_\lambda) = R^{GL_m}_{GL_\lambda}\Big( (St_1)^{\otimes m_{-1}}\otimes\bigotimes_{j\geqslant 0}(St_{e\ell^j})^{\otimes m_j}\Big).$$
If we denote by $\bfT_\lambda$ an $F_{q^2}$-stable maximal torus of $\bfG\bfL_m$ of 
cycle-type~$\lambda$, the class of $X_\lambda$ coincides with the Deligne-Lusztig
character induced from $\bT_\lambda$. 

\begin{lemma}\label{lem:X}
If $\lambda$ is a cuspidal partition of $m$ then
$[X_{\lambda}]=R_{\bfT_{\lambda}}^{\bfG\bfL_m}(1)$ in $[\scrU_\k]$.
\end{lemma}

\begin{proof}The finite group $T_\lambda=\bfT_{\lambda}(q^2)$ is isomorphic to
$\bfG\bfL_1(q^2)^{m_{-1}}\times\prod_{j\geqslant 0}\bfG\bfL_1(q^{2e\ell^j})^{m_j}$.
By \cite[thm.~19.18]{CE04}, for each integer $j\geqslant 0$ there is 
a character $\theta$ of $K \bfG\bfL_1(q^{2e\ell^j})$ of order a power of $\ell$ such that we have
$$R_{\bfT_{(e\ell^j)}}^{\bfG\bfL_{e\ell^j}}(1)=d_{\scrO GL}R_{\bfT_{(e\ell^j)}}^{\bfG\bfL_{e\ell^j}}([\theta])
=(-1)^{e\ell^j-1}[St_{e\ell^j}]=[St_{e\ell^j}].$$
We deduce that
$R_{\bfT_{\lambda}}^{\bfG\bfL_{\lambda}}(1)=[St_{\lambda}].$
By transitivity of the Lusztig induction we get 
\begin{align*}
[X_{\lambda}]=R^{GL_m}_{GL_{\lambda}}([St_{\lambda}])=R_{\bfT_{\lambda}}^{\bfG\bfL_m}(1).
\end{align*}
\end{proof}

Given a cuspidal partition $\lambda$, we set 
$B_\lambda=B_{X_\lambda}.$ Note that the endofunctor $F$ of $\scrU_\k$ is equal to $B_{1}$.
The functor $B^*_{\lambda}$ is not well-defined, but we may still consider the linear endomorphism 
$[B^*_{\lambda}]$ of $[\scrU_\k]$. They are related to the Heisenberg operators considered earlier by the following proposition
(which will not be used in the rest of the paper).

\begin{proposition} \label{prop:functors} Assume that $\ell\neq 2$. 
Let $\lambda$ be a cuspidal partition of $m$ of type $(m_j)$ with $m_{-1}=0$.
Then, we have the equalities $b_{2\lambda}=[B_{\lambda}]$
and $b_{2\lambda}^*=[B^*_{\lambda}]$ in $\End([\scrU_\k])$.
\end{proposition}

To prove the proposition we need more material.
Fix a positive integer $a.$ We abbreviate $GL_{m}^{(a)}=\bfG\bfL_{m}(q^{2a})$ and
we consider the category
$$\scrV_R^{(a)}:=\bigoplus_{m\in\bbN}R GL_m^{(a)}\umod.$$

\smallskip

First, we view $GL_m^{(a)}$ as the group of $\bbF_{q^{2a}}$-rational points of the rational group $(\bfG\bfL_m\,,\,F_{q^{2a}})$.
Unipotent irreducible modules of $K GL_m$ and $K GL_m^{(a)}$ are both parametrized
by partitions of $m$. Therefore there is an isomorphism $[\scrV_\k] \simto [\scrV_\k^{(a)}]$ sending $[W_\mu]$ to $[W_\mu^{(a)}]$ for any partition $\mu$ of $m$.
More generally we will denote by $[V^{(a)}]$ the image of $[V]$ under this map. 
By composition we obtain $\bbC$-linear isomorphisms given by
\begin{align}\label{isomC}
[\scrV_K^{(a)}]\simto[\scrV_\k^{(a)}]\simto\pmb\Lambda,\quad 
[L_\mu^{(a)}]\mapsto[W_\mu^{(a)}]\mapsto\phi_\mu.
\end{align} 

\smallskip

Next, we view $GL_m^{(a)}$ as the group of rational points of $\bfG\bfL_m^{(a)}=((\bfG\bfL_m)^{a}\,,\,\sigma_a F_{q^2}),$
where $\sigma_a$ is the automorphism of $(\bfG\bfL_m)^{a}$ such that
$(g_1,g_2,\dots,g_a)\mapsto (g_2,\dots,g_a,g_1)$.
Note that $\bfG\bfL_m^{(a)}$ is a Levi subgroup of $\bfG\bfL_{am}=(\bfG\bfL_{am}\,,\,F_{q^2})$, but the rational structures are different.
The Lusztig induction and restriction yield $\bbC$-linear maps 
\begin{align*}
R^{\bfG\bfL_{am}}_{\bfG\bfL_m^{(a)}}&\,:\,[\k GL_{m}^{(a)}\umod]\to [\k GL_{am}\umod],\\
{}^*R^{\bfG\bfL_{am}}_{\bfG\bfL_m^{(a)}}&\,:\,[\k GL_{am}\umod]\to [\k GL_{m}^{(a)}\umod].
\end{align*}

\smallskip

Let  $\scrF^{a}$ be the $\bbC$-linear map $\pmb\Lambda\to\pmb\Lambda$ such that 
$\scrF^a(z_\lambda\, c_\lambda)=z_{a\lambda}\,c_{a\lambda}$ for each partition $\lambda$.
Let
$\scrF_{\!a}$ be the transposed map relative to the bilinear form on $\pmb\Lambda$.
We have
$\scrF_{\!a}(c_\lambda)=c_{\mu}$
if $\lambda=a\mu$ and $\scrF_{\!a}(c_\lambda)=0$
if $\lambda\neq a\mu$ for all $\mu$.


\begin{lemma}\label{lem:ind/res} \hfill
\begin{itemize}[leftmargin=8mm]
  \item[$\mathrm{(a)}$]
The isomorphism \eqref{isomC} identifies the maps $R^{\bfG\bfL_{am}}_{\bfG\bfL_m^{(a)}}$ and ${}^*R^{\bfG\bfL_{am}}_{\bfG\bfL_{m}^{(a)}}$
with the endomorphisms $\scrF^{a}$ and $\scrF_{\!a}$ of $\pmb\Lambda$.
\item[$\mathrm{(b)}$]
The isomorphism \eqref{isomC0} identifies the maps $R_{GL_{m}\times GL_n}^{GL_{m+n}}$ and ${}^*R^{GL_{m+n}}_{GL_{m}\times GL_n}$
with the maps $\Ind^{m+n}_{m,n}$ and $\Res^{m+n}_{m,n}$ between $\pmb\Lambda\otimes\pmb\Lambda$ and $\pmb\Lambda$.
 \item[$\mathrm{(c)}$]
 If $n=r+2m$, the isomorphism $[\scrU_\k]\simto\pmb\Lambda$ sending $\varepsilon_\lambda[V_\lambda]$ to $\phi_\lambda$ and the isomorphism
$[\scrV_\k]\simto\pmb\Lambda$ in \eqref{isomC} identify the maps
$R^{G_n}_{L_{r,m}}$ and ${}^*R^{G_n}_{L_{r,m}}$ with the maps
$\Ind^n_{r,2m}(\bullet\otimes\scrF^2(\bullet))$ and
$(1\otimes\scrF_{2})\Res^n_{r,2m}$ between $\pmb\Lambda\otimes\pmb\Lambda$ and $\pmb\Lambda$.
\end{itemize}
\end{lemma}

\begin{proof}
The proof of the lemma is standard, see, e.g.,  \cite[prop.~6.1]{S14} or \cite[sec.~3.A]{BMM}.
Let us give some details.

\smallskip

First, we prove part (a).
We abbreviate $n=am$.
Consider the rational groups $(\bfG\bfL_{n},F_{q^2})$ and
$((\bfG\bfL_m)^{a}\,,\,\sigma_a F_{q^2})$.
We will denote them by $\bfG$ and $\bfG^{(a)}$ respectively.
The Weyl groups of $\bfG$ and $G$ are identified with $\frakS_n$.
The Weyl group of $\bfG^{(a)}$ is $(\frakS_m)^a.$ 
Let $\sigma_a$ denote also the permutation in $\frakS_n$ of cycle-type $(a^m)$ which normalizes the subgroup
$(\frakS_m)^a$ and induces the automorphism such that
$(w_1,w_2,\dots,w_a)\mapsto (w_2,\dots,w_a,w_1)$.
The action of the Frobenius homomorphism on the Weyl group of $\bfG^{(a)}$ is
given by the conjugation by $\sigma_a$.
The Weyl group of $G^{(a)}$ is $\frakS_m$.

\smallskip

By \eqref{L2} we have
\begin{align*}
[L_\lambda]= \sum_{\nu\,\vdash \,n}\langle \phi_\lambda\,,\,c_\nu\rangle_{\pmb\Lambda}\, R_{\bfT_\nu}^{\bfG}(1).
\end{align*}
Hence, the isometry $[\scrV_\k]\simto\pmb\Lambda$
such that $[L_\lambda]\mapsto\phi_\lambda$ maps the virtual module $R_{\bfT_\lambda}^{\bfG}(1)$
to $z_\lambda\,c_\lambda$ for each partition $\lambda$ of $n$.
Similarly, the isometry 
$[\scrV_\k^{(a)}]\to\pmb\Lambda$
in \eqref{isomC} maps
$[L_\mu^{(a)}]$ to $\phi_\mu$ for each partition $\mu$ of $m$, hence it maps the virtual module
$\big(R_{\bfT_\mu}^{\bfG\bfL_m}(1)\big)^{(a)}$ to $z_\mu\,c_\mu$. 

\smallskip

With the notation from Remark \ref{rem:twist} for Lusztig induction, we have
\begin{align*}
\big(R_{\bfT_\mu}^{\bfG\bfL_m}(1)\big)^{(a)}
&=\big(R_{\bfT_\mu\,,\,F_{q^2}}^{\bfG\bfL_m\,,\,F_{q^2}}(1)\big)^{(a)}\\
&=R_{\bfT_\mu\,,\,F_{q^{2a}}}^{\bfG\bfL_m\,,\,F_{q^{2a}}}(1)\\
&=R_{\bfT_w\,,\,\sigma_a F_{q^2}}^{(\bfG\bfL_m)^{a}\,,\,\sigma_a F_{q^2}}(1),
\end{align*}
where $w=(w_1,w_2,\dots,w_a)$ is any element in $(\frakS_m)^a$ such that
$w_1w_2\cdots w_a$ is of cycle-type $\mu$.
By transitivity of Lusztig induction we deduce that
$$R^{\bfG}_{\bfG^{(a)}}\big(R_{\bfT_\mu}^{\bfG\bfL_m}(1)\big)^{(a)}=R_{\bfT_{w\sigma_a}}^{\bfG}(1)=R_{\bfT_{a\mu}}^{\bfG}(1),$$
because $w\sigma_a$ is of cycle-type $a\mu$.
Therefore, up to the isometries above, the map
$$R^{\bfG}_{\bfG^{(a)}}\,:\,[\scrV_K^{(a)}]\to[\scrV_K]$$
is identified with the linear endomorphism $\scrF^a$ of $\pmb\Lambda$.
This proves the first identity in part (a).
The second one follows  by adjunction.

\smallskip

Now, we prove part (b).
Consider the $\bbF_{q^2}$-rational group $(\bfG\bfL_{m+n},F_{q^2})$ with its $F_{q^2}$-stable Levi subgroup $(\bfG\bfL_{m}\times\bfG\bfL_n,F_{q^2})$.
The isometry 
$[K GL_m\umod]\to\bbC\Irr(K\frakS_m)$
in \eqref{isomC0} maps the virtual module
$R_{\bfT_\mu}^{\bfG\bfL_m}(1)$ to $z_\mu\,c_\mu$ for each partition $\mu$ of $m$.
By the transitivity of Lusztig induction we  get
$$R^{GL_{m+n}}_{GL_m\times GL_n}R_{\bfT_\mu\times\bfT_\lambda}^{\bfG\bfL_m\times\bfG\bfL_n}(1)=R_{\bfT_\mu\times\bfT_\lambda}^{\bfG\bfL_{m+n}}(1).$$
This proves the first identity in part (b).
The second one follows  by adjunction.

\smallskip

It remains to prove part (c).
Consider the $\bbF_q$-rational group $(\bfG\bfL_n,F)$ with its $F$-stable Levi subgroup $\bfL_{r,m}$.
We will abbreviate $\bfG_n=(\bfG\bfL_n,F)$ and
we will identify the rational groups $(\bfL_{r,m},F)$ and $(\bfG\bfL_r\times(\bfG\bfL_m)^2, (F,\sigma_2F_q))$,
where $\sigma_2$ is an element of $\frakS_{2m}$ of cycle-type $(2^m)$.
Let $w_n$ denote the longest element in $\frakS_n$.
The Weyl groups of $\bfG\bfL_n$ and $G_n$ are identified with $\frakS_n$ and $W_n=C_{\frakS_n}(w_n)$.
The Weyl group of $\bfL_{r,m}$ is identified with 
$\frakS_r\times(\frakS_m)^2$ with
the action of the Frobenius endomorphism $F$ given by conjugation by $(w_r,\sigma_2)$.
Here, we embed $\frakS_r\times(\frakS_m)^2$ into $\frakS_n$ via the map $(u,v)\mapsto (w_mv_1w_m,u,v_2)$ if $v=(v_1,v_2)$
with $v_1,v_2\in\frakS_m$.
The Weyl group of $L_{r,m}$ is identified with $W_r\times\frakS_m$.

\smallskip

From \eqref{unip.char} we deduce that the isometry $[\scrU_K]\simto\pmb\Lambda$
such that $\chi_\lambda\mapsto\phi_\lambda$ maps the virtual module $R_{\bfT_\lambda}^{\bfG_r}(1)$
to $z_\lambda\,c_\lambda$ for each partition $\lambda$ of $r$.
In the proof of (a) we already observed that
the isometry 
$[\scrV_K]\simto\pmb\Lambda$
in \eqref{isomC} maps the virtual module
$R_{\bfT_\mu\,,\,F_{q^2}}^{\bfG\bfL_m\,,\,F_{q^2}}(1)=
R^{(\bfG\bfL_m)^2\,,\,\sigma_2F_q}_{\bfT_\mu\times\bfT_1\,,\,\sigma_2F_q}(1)$ to $z_\mu\,c_\mu$
for each partition $\mu$ of $m$.
By transitivity of Lusztig induction we  get
\begin{align*}
R^{\bfG\bfL_n\,,\,F}_{\bfL_{r,m}\,,\,F}\big(
R^{\bfG\bfL_r\times(\bfG\bfL_m)^2\,,\,(F,\sigma_2F_q)}_{\bfT_\lambda\times \bfT_\mu\times\bfT_1\,,\,(F,\sigma_2F_q)}(1)\big)=
R^{\bfG\bfL_n\,,\,F}_{\bfT_\mu\times\bfT_\lambda\times \bfT_1\,,\,F}(1).
\end{align*}
We deduce that the map $R^{G_n}_{L_{r,m}}$
is identified with the map $\bbC\Irr(\frakS_r)\otimes\bbC\Irr(\frakS_m)\to\bbC\Irr(\frakS_n)$ given by
$\Ind^n_{r,2m}(\bullet\otimes\scrF^2(\bullet))$.
This proves the first identity in part (c).
The second one follows  by adjunction.
\end{proof}

\smallskip

\begin{proof}[Proof of Proposition \ref{prop:functors}] 
Let $n=r+2m$. Under the isomorphism \eqref{isomC0} and the isomorphism
$[\scrU_\k]\simto\pmb\Lambda$ given by $\varepsilon_\lambda[V_\lambda]\mapsto\phi_\lambda$,
Lemma \ref{lem:ind/res} identifies the maps
$R^{G_n}_{L_{r,m}}$ and ${}^*R^{G_n}_{L_{r,m}}$ 
between $[\scrU_\k\otimes\scrV_\k]$ and $[\scrV_\k]$ with the $\bbC$-linear maps
$\pmb\Lambda\otimes\pmb\Lambda\to\pmb\Lambda$ and $\pmb\Lambda\to\pmb\Lambda\otimes\pmb\Lambda$ given by
\begin{equation}\label{eq:ind-in-lambda}
\phi_\nu\otimes\phi_\lambda\mapsto\Ind^n_{r,2m}(\phi_\nu\otimes\scrF^2(\phi_\lambda)),\qquad
\phi_\nu\mapsto(1\otimes\scrF_2)\Res^n_{r,2m}(\phi_\nu).
\end{equation}
Now, assume that $\lambda$ is cuspidal of type $(m_j)$.
By Lemmas \ref{lem:X}, \ref{lem:ind/res} the class of $X_\lambda$ in $[\scrV_\k]$ is identified with the element 
$z_\lambda c_\lambda\in\pmb\Lambda$. 
From \eqref{Heisenberg} we deduce that the operator
$[B_\lambda]=\bigoplus_rR^{G_{r+2m}}_{L_{r,m}}(\bullet\otimes[X_{\lambda}])$ on $[\scrU_\k]$ is identified
with the operators $b_{2\lambda}$ on $\pmb\Lambda$.
Hence, the proposition follows from Lemma \ref{lem:isomFockHeisenberg} and Remark \ref{rem:epsilon}.
\end{proof}

\smallskip

\medskip

\subsection{Cuspidal modules and highest weight vectors}\hfill\\

By Proposition \ref{prop:wHC},
the set of isomorphism classes of weakly cuspidal modules in $\scrU_\k$ is a basis of the space
$[\scrU_\k]^{0}$ of all elements in $[\scrU_\k]$ which are killed by the map $[E_i]$ for each $i\in I_e$. We define
\begin{align*}
[\scrU_\k]^\hw
&=\{x\in[\scrU_\k]^{0}\,|\,b_n^*(x)=0,\,\forall n\in 2e\bbZ_{>0}\}.
\end{align*}
Then, we have the following inclusion.  

\begin{lemma} $\{[D]\,$\,|\,$D\in\Irr(\scrU_\k)$ is cuspidal $\}\subseteq [\scrU_\k]^\hw.$
\end{lemma}

\begin{proof}
For any partition $\lambda$ we set 
$$a_{2\lambda}=\sum_{\nu}\chi^\lambda_\nu\,z_{\nu}^{-1}\,b_{2\nu},\quad
a^*_{2\lambda}=\sum_{\nu}\chi^\lambda_\nu\,z_{\nu}^{-1}\,b^*_{2\nu},$$
where 
$\chi^\lambda_\nu=\langle\phi_\lambda\,,\,z_\nu\, c_\nu\rangle_{\pmb\Lambda}$ 
is the value of the irreducible character $\phi_\lambda$
on a permutation of cycle-type $\nu$.
Since the $\scrP\times\scrP$-matrix with entries $\chi^\lambda_\nu$ is invertible, we have
\begin{align*}
[\scrU_\k]^\hw
&\supseteq\{x\in[\scrU_\k]^{0}\,|\,a_{2\lambda}^*(x)=0,\,\forall\lambda\in\scrP\}.
\end{align*}

Now, for any partition $\lambda$ we set $A_{\lambda}=B_{W_{\lambda}}$ and $[A^*_{\lambda}]=[B_{W_{\lambda}}^*]$.
Under the isomorphism $[\scrU_\k]\simeq\pmb\Lambda$, 
the operator $[A_\lambda]$ 
on $[\scrU_\k]$ is identified by \eqref{eq:ind-in-lambda} with the operator
$\bigoplus_r\Ind^{r+2m}_{r,2m}(\bullet\otimes\scrF^2(\phi_\lambda))$ on $\pmb\Lambda$.
Since $\phi_\lambda=\sum_{\nu}\chi^\lambda_\nu\,c_{\nu},$ the formula \eqref{Heisenberg} yields
$$\bigoplus_r\Ind^{r+2m}_{r,2m}(\bullet\otimes\scrF^2(\phi_\lambda))
=\sum_\nu \chi^\lambda_\nu\,z_{\nu}^{-1}\,\,b_{2\nu}.$$
We deduce that we have 
$a_{2\lambda}=[A_{\lambda}]$ and $a^*_{2\lambda}=[A^*_{\lambda}]$
in $\End([\scrU_\k])$.
Now, to prove the lemma it is enough to prove that
$a_{2\lambda}^*([D])=0$ for all partition $\lambda$ and all cuspidal module $D$ in $\Irr(\scrU_\k)$.
This follows from Proposition \ref{prop:B}.
\end{proof}

\subsubsection{The parameters of the ramified Hecke algebras}
It is not obvious that the vector space $[\scrU_\k]^\hw$ is spanned by classes of irreducible unipotent modules.
Our goal is to determine precisely the cuspidal modules. Let us first recall the following basic facts.

\begin{lemma}\label{lem:cuspidalpairs}\hfill
\begin{itemize}[leftmargin=8mm]
  \item[$\mathrm{(a)}$]
Any unipotent cuspidal pair of $\k G_n$  is $G_n$-conjugate to a pair $(G_r\times GL_\lambda\,,D_\nu\otimes St_\lambda)$ where $n=r+2m$, 
$\lambda$ is a cuspidal partition of $m$ and $D_\nu$ is a unipotent cuspidal $\k G_r$-module.
\item[$\mathrm{(b)}$]
Given a cuspidal pair $(G_r\times GL_\lambda\,,D_\nu\otimes St_\lambda)$ with $\lambda$ of type $(m_j)$, the ramified Hecke algebra 
$\scrH(\k G_n, D_\nu\otimes St_\lambda)$ is isomorphic to the $\k$-algebra
$\bfH_{\k,m_{-1}}^{Q_t\,;\,q^2}\otimes \bigotimes_{j\geqslant 0}\bfH_{\k,m_{j}}^{1,1\,;\,1}$, where the integer $t\in\bbN$ is such that 
the $e$-core of $\nu$ is $\Delta_t$.
\end{itemize}
\end{lemma}

\begin{proof}
Part (a) follows from Proposition \ref{prop:cusp-GL}.
To prove (b), let first observe that by \cite[prop.~4.4]{GHM96}
we have an algebra isomorphism 
$$\scrH(\k G_n, D_\nu\otimes St_\lambda)\simeq
\bfH_{\k,m_{-1}}^{P_{-1}\,;\,q^2}\otimes \bigotimes_{j\geqslant 0}\bfH_{\k,m_{j}}^{P_{j}\,;\,1},$$
where $P_j=(a_j,b_j)$ is a parameter in $\k^\times\times\k^\times$ for each $j\geqslant -1$.
Using the transitivity and faithfulness of induction, see \cite[prop.~1.23]{CE04} for details, we deduce from Corollary \ref{cor:2} that 
$P_{-1}=Q_t$ for some $t\geqslant 0$.
By \cite[prop.~2.3.5]{Gr96}\footnote{The reference \cite[prop.~2.3.5]{Gr96} has been indicated to us by G. Hiss.
} the parameter $P_j$ is of the form $(a_j,a_j)$ for some element $a_j\in\k^\times$
for each $j\geqslant 0$. 
\end{proof}

\subsubsection{The theorem} We can now formulate the theorem comparing cuspidal modules
and highest weight vectors.

\begin{theorem} \label{thm:cuspidal} 
The set $\{[D]$\,|\,$D\in\Irr(\scrU_\k)$ is cuspidal$\,\}$ is a basis of $[\scrU_\k]^\hw$.
\end{theorem}

\begin{proof}
First, let us fix some notation that we will use in the whole proof.
For every sequence $(m_{-1},m_0,m_1,\dots)$ of integers $\geqslant 0$ we set $m_+=\sum_{j\geqslant 0}m_j\ell^j$,
$m=m_{-1}+em_+$ and $n=r+2m$.
Let $\lambda$ be the cuspidal partition of $m$ of type $(m_{-1}, m_0,\dots)$
and $\lambda_+$ be the cuspidal partition of type $(m_0, m_1,\dots)$.
Let $\scrC\subseteq\scrP$ be the set of all partitions $\nu$ such that the module $D_\nu$ is cuspidal.
For each $\nu\in\scrC$ the $e$-core of $\nu$ is a $2$-core $\Delta_{t(\nu)}$.
We will abbreviate $\Delta_\nu=\Delta_{t(\nu)}$ and $Q_\nu=Q_{t(\nu)}$.

\smallskip

Now, we consider the Harish-Chandra series  of $\k G_n$. 
The crystal basis $B(\Lambda)$ of the $\frakg_e$-module $\bfL(\Lambda)$
is a disjoint union of bases $B(\Lambda)_\alpha$
of the weight subspaces $\bfL(\Lambda)_{\Lambda-\alpha}$
where $\alpha$ runs over $\Q^+_e$.
For each integer $d\geqslant 0$ we set
$$\bfL(\Lambda)_{d}=\bigoplus_{\height(\alpha)=d}\bfL(\Lambda)_{\Lambda-\alpha},
\quad
B(\Lambda)_{d}=\bigsqcup_{\height(\alpha)=d} B(\Lambda)_{\Lambda-\alpha}.$$
We first prove the following.

\begin{claim}\label{claim:A}
Given $\nu$ a partition, and $m_{-1}, m_+ \geq 0$, we have a bijection
\begin{align*}
\bigsqcup_{\lambda}\Irr(\k G_n, D_\nu\otimes St_\lambda)\, \mathop{\longleftrightarrow}\limits^{1:1}\,B(\Lambda_{Q_{\nu}})_{m_{-1}}\times\scrP_{m_+}
\end{align*}
where $\lambda$ runs over all cuspidal partitions with fixed $m_{-1}$ and $m_+$.
\end{claim}

From Lemma \ref{lem:cuspidalpairs} and the generalities in Section \ref{sec:HCseries}, we have a canonical bijection
\begin{align}\label{bijection1}
\Irr(\k G_n, D_\nu\otimes St_\lambda)\, \mathop{\longleftrightarrow}\limits^{1:1}\,\Irr(\bfH_{\k,m_{-1}}^{Q_{\nu};\,q^2})\times\prod_{j\geqslant 0}\Irr(\bfH_{\k,m_{j}}^{1,1\,;\,1}).
\end{align}
By \S \ref{subsec:minimalcat}, there is a canonical bijection
\begin{align}\label{bijection2}
B(\Lambda_{Q_{\nu}})_{m_{-1}}\, \mathop{\longleftrightarrow}\limits^{1:1}\,\Irr(\bfH_{\k,m_{-1}}^{Q_{\nu}\,;\,q^2}).
\end{align}
Further, the elements of $\Irr(\bfH_{\k,m_{j}}^{1,1\,;\,1})$ are labelled by
the $\ell$-restricted partitions of $m_j$, see, e.g., \cite[prop.~4.6.6]{GJ}.
Hence, the unicity of the $\ell$-adic expansion of partitions yields the following bijection 
\begin{align}\label{bijection3}
\scrP_{m_+}\, \mathop{\longleftrightarrow}\limits^{1:1}\,\bigsqcup_{\lambda_+}\prod_{j\geqslant 0}\Irr(\bfH_{\k,m_{j}}^{1,1\,;\,1}),
\end{align}
where $\lambda_+ = (e^{m_0},(e\ell)^{m_1},\ldots)$ runs over all cuspidal partitions 
of $m_+$ with parts divisible by $e$ (see, e.g., \cite[lem.~19.26]{CE}).
The Claim \ref{claim:A} follows from \eqref{bijection1}, \eqref{bijection2} and \eqref{bijection3}.

\smallskip

Now, we consider the complexified Grothendieck group $[\scrU_\k]$.
The modular Harish-Chandra theory yields bijections
\begin{align}\label{eqC}\Irr(\scrU_\k)\, \mathop{\longleftrightarrow}\limits^{1:1}\,\bigsqcup_{n,\nu,\lambda} \Irr(\k G_n, D_\nu\otimes St_\lambda).
\end{align}
Next, we identify $[\scrU_\k]$ with the $\frakg_e$-module $\bigoplus_t\bfF(Q_t)$.
Fix a basis $B$ of $[\scrU_\k]^\hw$ containing the family $\{[D_\nu]\,|\,\nu\in\scrC\}$.
By Corollary \ref{cor:2} we can identify the $\frakg'_e$-submodule of $[\scrU_\k]$ generated by $[D_\nu]$ with $\bfL(\Lambda_{Q_{\nu}})$.
We get a canonical isomorphism of $\frakH_{2e}\oplus\frakg'_e$-modules
\begin{align}\label{eqB}[\scrU_\k]=\bigoplus_{\nu\in\scrC}U(\frakH_{2e})(\bfL(\Lambda_{Q_{\nu}}))\oplus
\bigoplus_{b\in B\setminus\scrC}U(\frakH_{2e})U(\frakg'_e)(b).\end{align}
For each $n$ let $[\scrU_\k]_n$ be the subspace of $[\scrU_\k]$ given by
$$[\scrU_\k]_n=\bigoplus_{\nu,m_{-1},m_+,\gamma}b_{2e\gamma}(\bfL(\Lambda_{Q_{\nu}})_{m_{-1}}),$$
where $\gamma\in\scrP_{m_+}$ and $r, m_{-1}, m_+\in\bbN$ are such that $n=r+2m_{-1}+2em_+$.
Comparing \eqref{eqC} and \eqref{eqB}, we deduce that 
$$\sum_{\nu,\,\lambda}|\Irr(\k G_{n}, D_\nu\otimes St_\lambda)|\geqslant\dim([\scrU_\k]_n),$$
and that we have equality for all $n$ if and only if $B=\{[D_\nu]\,|\,\nu\in\scrC\}$.

\smallskip

Now, any non-zero highest-weight $\frakH_{2e}$-module is free as a module over the subalgebra generated by 
$\{b_{2en}\,|\,n\geqslant 1\}$.
Hence, for each $n$, the Claim \ref{claim:A} gives
\begin{align*}
\sum_{\nu,\,\lambda}|\Irr(\k G_n, D_\nu\otimes St_\lambda)|=\dim([\scrU_\k]_n),
\end{align*}
proving the theorem.
\end{proof}

\begin{remark}
The theorem implies that 
$$\{\,[D]\,|\,D\in\Irr(\scrU_\k) \text{\ is\ cuspidal}\}=
\Irr(\scrU_\k)\cap [\scrU_\k]^\hw.$$
\end{remark}

\smallskip

\subsubsection{Cuspidal modules and FLOTW $e$-partitions}\label{sec:FLOTW}

Recall that in Theorem \ref{thm:wHC} we have given an explicit isomorphism between the crystal basis $B(\scrU_\k)$,
whose underlying set is $\Irr(\scrU_\k)$, and the abstract crystal $B_e$ of the $\frakg_e\times\frakH_{2e}$-module
$\bigoplus_t\bfF(Q_t)_e$.
Using Theorem \ref{thm:cuspidal},  we deduce that the cuspidal 
modules in $\Irr(\scrU_\k)$ corresponds precisely to the elements of $B_e$ whose associated elements
in the upper global basis $\bfB^\vee_e:=\bigsqcup_t\bfB^\vee(s_t)$ of $\bigoplus_t\bfF(Q_t)_e$,
see \S \ref{subsec:crystalfock}, are highest weight vectors for the
action of both $\frakg_e$ and $\frakH_{2e}$.

\smallskip

These elements of $B_e$ have been computed recently in \cite{G16}.
To do so, one defines a class of $e$-partitions called FLOTW, see  \cite[def.~6.22]{G16}.
Then, the cuspidal modules are related to FLOTW $e$-partitions as in \cite[thm.~7.7]{G16}.

\medskip

\section{Types B and C}\label{sec:BC}

This section  is devoted to the construction of categorical actions on the category of unipotent representations of the finite groups
$\SO_{2n+1}(q)$ and  $\Sp_{2n}(q)$.
The arguments are similar to those in Section \ref{part:unitarygroups} and will yield:
 \begin{itemize}[leftmargin=8mm]
  \item[(1)] the Hecke algebras associated to weakly cuspidal representations, 
  \item[(2)] the branching graph for the parabolic induction and restriction,
  \item[(3)] derived equivalences between blocks.
 \end{itemize}
The main difference is the lack of a description in terms of level $1$ Fock spaces (given by Ennola duality, see \S\ref{part:unitarygroups}). In particular, the grading by the imaginary roots, which in Section
\ref{part:unitarygroups} was deduced from the grading in the $\GL_n$ case, is constructed in types B, C via some explicit combinatorics
in terms of Lusztig's symbols.

\subsection{Definitions} Throughout this section and the following, $q$ is any prime power.
 Note however that we will assume that $q$ is odd starting from Section \ref{sec:cat-uK-and-uk}.

\subsubsection{Odd-dimensional orthogonal groups} 
Let $J_n$ be as in \S \ref{sec:def-unitary}. It is the matrix of a non-degenerate quadratic form.
The odd-dimensional orthogonal group $\bfSO_{2n+1} = \SO_{2n+1}(\overline\bbF_q)$ is  
  $$  \mathrm{SO}_{2n+1}(\overline\bbF_q) = \{g \in \bfS\bfL_{2n+1} \, \mid 
  \, {}^tg  J_{2n+1} g =   J_{2n+1}\}.$$
It is a connected reductive group of type $B_n$. The standard Frobenius map $F=F_q$ on $\bfG\bfL_{2n+1}$
induces a Frobenius endomorphism on $\bfSO_{2n+1}$. The \emph{finite orthogonal group}  $\SO_{2n+1}(q)$ is
given by
  $$\SO_{2n+1}(q)= (\bfSO_{2n+1})^F.$$
The subgroup of diagonal matrices $\bfT$ (resp. upper-triangular matrices $\bfB$) in $\bfSO_{2n+1}$ is
a split maximal torus (resp. an $F$-stable Borel subgroup). Note that $\bfT = \{ \mathrm{diag}(t_n,
\ldots,t_1,1,t_1^{-1},\ldots,t_n^{-1})\}$. The Weyl group $W$ of $(\bfSO_{2n+1},\bfT)$ is a Weyl group of
type $B_n$, and $F$ acts trivially on it. For numbering the simple reflections of $W$ we will take
the following convention:
\begin{align}\label{B}
\begin{split}
\begin{tikzpicture}[scale=.4]
\draw[thick] (0 cm,0) circle (.5cm);
\node [below] at (0 cm,-.5cm) {$s_1$};
\draw[thick] (0.5 cm,-0.15cm) -- +(1.5 cm,0);
\draw[thick] (0.5 cm,0.15cm) -- +(1.5 cm,0);
\draw[thick] (2.5 cm,0) circle (.5cm);
\node [below] at (2.5 cm,-.5cm) {$s_2$};
\draw[thick] (3 cm,0) -- +(1.5 cm,0);
\draw[thick] (5 cm,0) circle (.5cm);
\draw[thick] (5.5 cm,0) -- +(1.5 cm,0);
\draw[thick] (7.5 cm,0) circle (.5cm);
\draw[dashed,thick] (8 cm,0) -- +(4 cm,0);
\draw[thick] (12.5 cm,0) circle (.5cm);
\node [below] at (12.5 cm,-.5cm) {$s_{n-1}$};
\draw[thick] (13 cm,0) -- +(1.5 cm,0);
\draw[thick] (15 cm,0) circle (.5cm);
\node [below] at (15 cm,-.5cm) {$s_n$};
\end{tikzpicture}
\end{split}
\end{align}

For $i \neq 1$, the action of $s_i$ on $\mathrm{diag}(t_n,\ldots,t_1,1,t_1^{-1},\ldots,t_n^{-1})$
swaps $t_{i-1}$ and $t_i$, whereas $s_1$ swaps $t_1$ and $t_1^{-1}$. For $i \neq 1$, the simple
reflection $s_i$ can be lifted to $N_{\bfSO_{2n+1}}(\bfT)$ as the permutation matrix
$(n-i+1,n-i+2)(n+i,n+i+1)$ and $s_1$ as the signed permutation matrix 
\begin{equation*}
s_1= - \ \begin{blockarray}{(ccc|ccc|ccc)}
         \BAmulticolumn{3}{c|}{\multirow{2}{*}{$\Id_{n-1}$}}&&&&&&\\
         &&&&&&&&&\\
         \cline{1-9}
        &&&0&0&1&&&&\\
        &&&0&1&0&&&&\\
        &&&1&0&0&&&&\\
        \cline{1-9}
        &&&&&&\BAmulticolumn{3}{c}{\multirow{2}{*}{$\Id_{n-1}$}}\\
        &&&&&&&&&\\
    \end{blockarray}\,\, .
\end{equation*}

\smallskip

\subsubsection{Symplectic groups} 
Let $\widetilde J_{2n}$ be the $2n\times2n$ square matrix 
$$\widetilde J_{2n}=\left( \begin{array}{cc} 0&J_n\\-J_n&0\end{array}\right).$$
It is the matrix of a symplectic form. The symplectic group $\bfSp_{2n}=\Sp_{2n}(\overline\bbF_q)$ is
  $$ \bfSp_{2n} =  \{g \in \bfG\bfL_{2n} \, \mid \, {}^tg \widetilde J_{2n} g =  \widetilde J_{2n}\}.$$
It is a connected reductive groups of type $C_n$. The standard Frobenius map $F=F_q$ on $\bfG\bfL_{2n}$
induces a Frobenius endomorphism on $\bfSp_{2n}$.
The \emph{finite symplectic group} $\Sp_{2n}(q)$ is given by
$$\Sp_{2n}(q) = (\bfSp_{2n})^F.$$ 
The subgroup of diagonal matrices $\bfT$ (resp. upper-triangular matrices $\bfB$) in $\bfSp_{2n}$ is a split
maximal torus (resp. an $F$-stable Borel subgroup). Note that $\bfT = \{\mathrm{diag}(t_n,\ldots,t_1,t_1^{-1},
\ldots,t_n^{-1})\}$. The Weyl group $W$ of $(\bfSp_{2n},\bfT)$ is a Weyl group of type $B_n$, and $F$
acts trivially on it. For numbering the simple reflections of $W$ we will take the same convention as in \eqref{B}.

\smallskip

As in the case of odd-dimensional orthogonal groups, the action of $s_i$ on the element $\mathrm{diag}(t_n,
\ldots,t_1,t_1^{-1},\ldots,t_n^{-1})$ swaps $t_{i-1}$ and $t_i$ when $i \neq 1$, whereas $s_1$ swaps
$t_1$ and $t_1^{-1}$. For $i \neq 1$, we will lift $s_i$ to $N_{\bfSp_{2n}}(\bfT)$ as the permutation matrix
$(n-i+1,n-i+2)(n+i-1,n+i)$ and $s_1$ as the signed permutation matrix 
\begin{equation*}
s_1= \begin{blockarray}{(cccc|cc|cccc)}
         \BAmulticolumn{4}{c|}{\multirow{2}{20pt}{$\Id_{n-1}$}}&&&&&\\
         &&&&&&&&&\\
         \cline{1-10}
        &&&&0&-1&&&&\\
        &&&&1&0&&&&\\
        \cline{1-10}
        &&&&&&\BAmulticolumn{3}{c}{\multirow{2}{15pt}{$\Id_{n-1}$}}\\
        &&&&&&&&&
    \end{blockarray}\ .
\end{equation*}

\medskip

\subsection{The representation datum on $RG$-mod}\label{sec:repdatumB}\hfill\\

Throughout this section we assume that $R$ is any commutative domain in which $q(q-1)$ is invertible.

\smallskip

Given a positive integer $n$, we will denote by $\bfG_n$ and $G_n$ either the odd-dimensional
orthogonal groups $\bfSO_{2n+1}$ and $\SO_{2n+1}(q)$ or the symplectic groups $\bfSp_{2n}$
and $\Sp_{2n}(q)$. In addition, we set $G_0 = \{1\}$ by convention. 

\smallskip

Given $r,m\in\bbN$ such that $n=r+m$, let $\bfL_{r,m}$ and $\bfL_{r,1^m}$ be the standard Levi
subgroups of $\bfG_n$ corresponding to the sets of simple reflexions $\{s_k\,|\,k\neq r+1\}$,
$\{s_k\,|\,k\leq r\}$ respectively. The corresponding finite groups are $L_{r,m}
\simeq G_r \times \GL_m(q)$ and $L_{r,1^m}\simeq G_r \times \GL_1(q)^m.$ We will abbreviate
$\bfL_r=\bfL_{r,1}$. It is isomorphic to $\bfG_r \times \bfG\bfL_1$ via the map
$$ (g,\lambda) \longmapsto \left( \begin{array}{ccc} \lambda  & &  \\  & g & \\ & & 
  \lambda^{-1}   \end{array} \right) .
$$
Let $\bfP_{r}\subset\bfG_{r+1}$ be the corresponding parabolic subgroup and $P_r=\bfP_r^F$.
Let $\bfV_r$ be the unipotent radical of $\bfP_{r}$ and $V_r = \bfV_r^F$. Let us consider the subgroup $U_{r}\subset G_{r+1}$ given by
  $$U_{r}=V_{r}\rtimes\bbF^\times_{q}.$$
It is represented by the same figure as in \eqref{U}. For each $r<n$ we set $U_{n,r}=
U_{n-1}\rtimes\dots\rtimes U_{r}$ and $e_{n,r}=e_{U_{n,r}}$. The embedding of $G_r$ into $L_{r}$
yields an embedding of $G_r$ into $G_{r+1}$, and by induction, of $G_r$ into $G_n$. We obtain
functors
\begin{align*}
F_{n,r}=RG_n\cdot e_{{n,r}}\otimes_{RG_r}\bullet\,:\,RG_r\mod\to RG_n\mod,\\
E_{r,n}=e_{{n,r}}\cdot RG_n\otimes_{RG_n}\bullet\,:\,RG_n\mod\to RG_r\mod.
\end{align*}
An endomorphism of the functor $F_{n,r}$ can be represented by an $(RG_n,RG_r)$-bimodule
endomorphism of $RG_n\cdot e_{{n,r}}$, or, equivalently, by an element of $e_{{n,r}}\cdot
RG_n\cdot e_{{n,r}}$ centralizing $RG_r$. Thus, the elements
\begin{align}
& X_{r+1,r}= q^r e_{{r+1,r}} (s_{r+1}s_r \cdots s_1 \cdots s_r s_{r+1}) \, e_{{r+1,r}},\qquad
T_{r+2,r}= q e_{{r+2,r}} s_{r+2}\,e_{{r+2,r}}
\end{align}
define respectively natural transformations of the functors $F_{r+1,r}$ and $F_{r+2,r}$.
Indeed, with our convention $s_{r+1}s_r \cdots s_1 \cdots s_r s_{r+1}$ is one of the matrices
  $$ -\left(\begin{array}{ccc} & & 1 \\ & \Id_{G_r} & \\ 1 \\ \end{array}\right)
  \quad \text{or} \quad \left(\begin{array}{ccc} & & -1 \\ & \Id_{G_r} & \\ 1 \\ \end{array}\right)$$
which centralize $G_r$. We set 
$$F=\bigoplus_{r\geqslant 0}F_{r+1,r},\qquad 
X=\bigoplus_{r\geqslant 0}X_{r+1,r},\qquad
T=\bigoplus_{r\geqslant 0}T_{r+2,r}.$$

\begin{proposition}\label{prop:RD} 
The endomorphisms $X\in\End(F)$ and $T\in\End(F^2)$ satisfy the following relations:
\begin{itemize}[leftmargin=8mm]
  \item[$\mathrm{(a)}$] $1_FT\circ T1_F\circ 1_FT=T1_F\circ 1_FT\circ T1_F$,
  \item[$\mathrm{(b)}$] $(T+1_{F^2})\circ(T-q1_{F^2})=0$,
  \item[$\mathrm{(c)}$] $T\circ(1_FX)\circ T=qX1_F$.
\end{itemize}
\end{proposition}

\begin{proof} Similar to the proof of  \cite[prop.~4.1]{DVV}.
\end{proof} 

\medskip

\subsection{The categories of unipotent modules $\scrU_K$ and $\scrU_\k$\label{sec:cat-uK-and-uk}}\hfill\\

From now on, we fix an $\ell$-modular system $(K,\scrO,\k)$ which we assume to be 
big enough so that every indecomposable unipotent representation of 
$G_n$ (over $K$ or $\k$) is absolutely indecomposable. Such a modular system exists
since unipotent representations of $G_n$ are defined over $\bbQ_\ell$,
see \cite[cor.~1.12]{L02}. 
In addition, since we will be dealing with representations in non-defining
characteristic, we will always assume that $\ell \neq p$. We denote by $f$ (resp. $d$)
the order of $q$ (resp. $q^2$) in $\k^\times$. If $f$ is odd, we say that the prime
number $\ell$ is \emph{linear}, otherwise we say that $\ell$ is \emph{unitary}.
In the first case we have $f=d$, whereas in the second one $f=2d$. 

\smallskip

From now on, we will always assume that both $p$ and $\ell$ are odd, and that
$f > 1$. In particular $q(q-1) \in \scrO^\times$ and we can apply the previous constructions
with $R$ being any ring among $(K,\scrO,\k)$.

\smallskip

\subsubsection{Parametrization by symbols}\label{sec:symbols}
By \cite{L77}, the unipotent characters of $G_n$ are parametrized by symbols. For our
purpose it will be more convenient to work with a slightly different notion which 
we refer to as charged symbols. 

\smallskip

A \emph{charged symbol} $\Theta$ with charge $s =(s_1,s_2)$
is a pair of charged $\beta$-sets $\Theta = (\beta_{s_1}(\mu^1),\beta_{s_2}(\mu^2))$ for
some bipartition $\mu = (\mu^1,\mu^2)$. We abbreviate $\Theta = \beta_s(\mu)$. 
If $\beta_{s_1}(\mu^1) =: X =  \{x_1 > x_2 > \cdots\}$ and $\beta_{s_2}(\mu^2) =: Y
= \{y_1 > y_2 > \cdots\}$ we write
 $$ \Theta  = (X,Y) \, = \, \left( \begin{array}{ccccc}
 x_1 &  x_2  & \cdots  \\ y_1 & y_2 & \cdots  \\ \end{array} \right).$$
The components $X$ and $Y$ are called the \emph{first and second row} of the symbol.
The \emph{defect} of $\Theta$ is $D = s_1-s_2$ and its \emph{rank} is $|\mu| + \lfloor D^2/4 \rfloor$.

\smallskip

A \emph{$d$-hook} of $\Theta$ is a pair of integers $(x,x+d)$ which is either
a $d$-hook of $X$ or a $d$-hook of $Y$. The charged symbol obtained by deleting $x+d$ from $X$ (resp. $Y$)
and replacing it by $x$ is said to be gotten from $\Theta$ by \emph{removing the $d$-hook $(x,x+d)$}. 
A \emph{$d$-cohook} is a pair of integers $(x,x+d)$ such that $x+d\in X$ and $x\not\in Y$, or
$x+d\in Y$ and $x\not\in X$. 
The charged symbol obtained by deleting $x+d$ from $X$ (resp. $Y$) and adding $x$ to $Y$ (resp. $X$) is said to be gotten from $\Theta$ by \emph{removing the $d$-cohook $(x,x+d)$}.
The \emph{$d$-core} of $\Theta$ is obtained by removing recursively all $d$-hooks from $\Theta$.
A similar definition using $d$-cohooks gives the \emph{$d$-cocore} of $\Theta$.

\smallskip

We will denote by $\Theta^\dag = (Y,X)$ the charged symbol of charge $(s_2, s_1)$ obtained
by swapping the two $\beta$-sets. The defect of $\Theta^\dag$ is $-D$ but the rank is the same.
If one shifts simultaneously the charged $\beta$-sets $X$ and $Y$ by an integer $m$, one
obtains a symbol $\Theta[m]$ of charge $(s_1+m,s_2+m)$. This operation does not change
the defect nor the rank. 

\smallskip

\emph{Symbols}  are orbits of charged symbols
under the shift operator and the transformation $\Theta \mapsto \Theta^\dag$. 
We write 
$$ \{X,Y\} \, = \, \left\{ \begin{array}{ccccc}
 x_1 &  x_2  & \cdots  \\ y_1 & y_2 & \cdots  \\ \end{array} \right\}$$ for the symbol associated with $(X,Y)$. The rank of
the symbol is the rank of any charged symbol in its class whereas its defect is $|D|$
where $D$ is the defect of any representative. Removing and adding $d$-hooks or $d$-cohooks
are well-defined operations on symbols.

\smallskip

We denote by $\scrS$ the set of symbols and by $\scrS_{\mathrm{odd}}$ the set of symbols of odd defects.

\smallskip

\subsubsection{The unipotent modules over $K$}\label{sec:symbol}
Fix a positive integer $n$. The unipotent $KG_n$-modules were classified in \cite[thm.~8.2]{L77}.
They are parametrized in terms of symbols of odd defect and rank $n$. This parametrization follows
from the determination of cuspidal unipotent $KG_n$-modules and their ramified Hecke
algebras.

\begin{proposition}[Lusztig \cite{L77}]\label{prop:cuspidalBC}
 Up to isomorphism, there is at most one cuspidal unipotent module in $KG_r\mod$, and
 there exists one if and only if $r=t(t+1)$ for some $t\geqslant 0$. It is denoted
 by $E_t$.
 \qed
\end{proposition}

Since there do not exist any cuspidal unipotent character of $\mathrm{GL}_n(q)$ unless $n =1$,
we deduce that any cuspidal pair of $G_n$ is conjugate to a pair of the form $(L_{r,1^m}\,,\, E_t)$
with $n=r+m$ and $r=t(t+1)$ with $t\geqslant 0$.  The general theory recalled in \S \ref{sec:HCseries} implies that the irreducible characters lying in the
Harish-Chandra series above $E_t$ are in bijection with the irreducible representations of
the ramified Hecke algebra $\scrH(KG_n,E_t)$. To describe explicitly the latter, we
introduce, as in \eqref{Q}, the parameters
\begin{equation}\label{eq:chargeBC}
	Q_t= \big((-q)^t,(-q)^{-1-t}\big).
\end{equation}
Then $\scrH(KG_n,E_t)$ is isomorphic to $\bfH^{Q_t\,;\,q}_{K,m}$ by \cite[\S 5,8]{L77}.
There is a canonical choice for this isomorphism given by Theorem \ref{thm:HL-BC} below.
Consequently, the Harish-Chandra theory yields a canonical bijection 
\begin{align}\label{bijectionB}\Irr(KG_n,E_t)\, \mathop{\longleftrightarrow}\limits^{1:1}\, \Irr(\bfH^{Q_t\,;\,q}_{K,m})\end{align}
and hence a parametrization of $\Irr(KG_n,E_t)$. More precisely, given
$\mu=(\mu^1,\mu^2)$ a bipartition of $m$ and $t \geqslant 0$ we can associate the symbol
$$\Theta_t(\mu)  = \big\{\beta_t(\mu^1),\beta_{-t-1}(\mu^2)\big\}.$$
So, the defect and the rank of $\Theta_t(\mu)$ are
$$D(\Theta_t(\mu))=2t+1,\quad\rk(\Theta_t(\mu))=m+t(t+1).$$
Then we define $E_{\Theta_t(\mu)}$ to be the unipotent $KG_r$-module corresponding to 
$S(\mu)_K^{Q_t\,;\,q}$ via the bijection \eqref{bijectionB}. This yields a parametrization of the
irreducible unipotent $KG_n$-modules as 
$$\{E_{\Theta_t(\mu)}\,|\,t\in\bbN\,,\,\mu\in\scrP^2_m\,,\,m+t(t+1)=n\}.$$ 
We will abbreviate $\Theta_t=\Theta_t(\emptyset)$. Note that $E_t=E_{\Theta_t}$.

\begin{proposition}[Lusztig \cite{L77}]\label{prop:HC-B}
The irreducible unipotent characters are parametrized by symbols of odd defect. 
If $n=r+m$ and $r=t(t+1)$ with $t \geqslant 0$, then the unipotent module
$E_\Theta$ lies in the Harish-Chandra series $\Irr(KG_n,E_t)$ if and only if $\Theta$
has defect $2t+1$ and rank $n$. \qed
\end{proposition}

It is important to observe that this parametrization is exactly the one described in
\cite[thm. 8.2]{L77}, which was subsequently used for the determination of blocks in \cite{FS89}.
Indeed, the labelling of unipotent characters by bipartitions or symbols in
\cite{L77} (see also \cite[\S 13.8]{Car}) corresponds to the presentation of the Hecke algebra $\scrH(KG_n,E_t)$ in terms of the set of generators $(T_0,T_1,\ldots,T_{m-1})$
instead of $(X_1,T_1,\ldots,T_{m-1})$ with $T_0 = (-1)^t q^{t+1} X_1$ so that
$(T_0 -q^{2t+1}) (T_0+1) =0$. Now the isomorphism $\bfH^{Q_t\,;\,q}_{K,m} \simeq 
\bfH^{(q^{2t+1},-1)\,;\,q}_{K,m}$ induces a bijection on irreducible objects which
is just $S(\mu)^{Q_t} \longmapsto S(\mu)^{(q^{2t+1},-1)}$. 

\begin{remark}
Note that $Q_t^\dag=Q_{-1-t}$, $\Theta_t(\mu)=\Theta_{-1-t}(\mu^\dag)$ and $t(t+1)$ is invariant under the map $t \longmapsto -t-1$. 
We will work with symbols $\Theta_t(\mu)$ such that $t \geqslant 0$.
We might also have worked with $t<0$ using the symmetries above.
\end{remark}

Recall that $G_0=\{1\}$. We call the category of unipotent $KG$-modules the category  
  $$\scrU_K=\bigoplus_{n\in\bbN}K G_n\umod.$$
This category is abelian semisimple. From the previous discussion we have $$\Irr(\scrU_K)=
\{E_\Theta\,|\,\Theta\in\scrS_{\mathrm{odd}}\},$$ where by convention $\Irr(KG_0)=\{E_0\}$.

\smallskip

\subsubsection{The unipotent modules over $\k$}
Using the $\ell$-modular system we have decomposition maps $d_{\scrO G_n}$
which by Proposition \ref{prop:isodecmap} and since $\ell$ is odd, restrict to linear isomorphisms 
 $$ d_{\scrO G_n}:[KG_n\umod]\simto[\k G_n\umod]. $$

\smallskip

We call the category of unipotent $\k G$-modules the category
  $$\scrU_\k=\bigoplus_{n\in\bbN}\k G_n\umod.$$
This is an abelian category which is not semisimple. As above, the decomposition map yields a $\bbZ$-linear
isomorphism $d_\scrU : [\scrU_K]\simto[\scrU_\k]$. It is conjectured 
that this map is unitriangular on the basis of irreducible modules, yielding
a parametrization of unipotent simple $\k G$-modules (see Conjecture \ref{conj:unitriangular} below).

\smallskip

\subsubsection{The unipotent blocks}\label{sec:blockB}
Recall that $d$ is the order of $q^2$ modulo $\ell$. The partition of the unipotent
$KG_n$-modules into $\ell$-blocks was determined in \cite[(10B), (11E)]{FS89}. If $\ell$ is
a linear prime (resp. a unitary prime), i.e., if $f$ is odd (resp. even),
then two unipotent characters of $KG_n$ belong to the same $\ell$-block
if and only if their symbols have the same $d$-core (resp. $d$-cocore). 

\smallskip

In addition, Fong-Srinivasan described the structure of the unipotent blocks of $G_n$ with
cyclic defect groups. First, the unipotent $\ell$-block containing $E_\Theta$ has a cyclic
defect if and only if the symbol $\Theta$ has a unique $d$-hook if $\ell$ is linear
or a unique $d$-cohook if $\ell$ is unitary. In this case, let $\{X,Y\}$ be the $d$-core of
$\Theta$ if $\ell$ is linear, and the $d$-cocore of $\Theta$ if $\ell$ is unitary.
By \cite[(5A), (6A)]{FS}, the Brauer tree of the $\ell$-block containing $E_\Theta$ is 
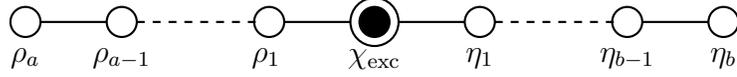
\begin{center}
\begin{tikzpicture}[scale=.4]
\draw[thick] (-1.7 cm,0) circle (.5cm);
\node [below] at (-1.7 cm,-.5cm) {$\rho_a$};
\draw[thick] (-1.2 cm,0) -- +(2.2 cm,0);
\draw[thick] (1.5 cm,0) circle (.5cm);
\node [below] at (1.5 cm,-.5cm) {$\rho_{a-1}$};
\draw[dashed,thick] (2 cm,0) -- +(4 cm,0);
\draw[thick] (6.4 cm,0) circle (.5cm);
\node [below] at (6.3 cm,-.5cm) {$\rho_1$};
\draw[thick] (6.9 cm,0) -- +(2.2 cm,0);
\draw[thick,fill=black] (9.9 cm,0) circle (.5cm);
\draw[thick] (9.9 cm,0) circle (.8cm);
\node [below] at (9.9 cm,-.55cm) {$\chi_{\text{exc}}$};
\draw[thick] (10.7 cm,0) -- +(2.2 cm,0);
\draw[thick] (13.4 cm,0) circle (.5cm);
\node [below] at (13.4 cm,-.5cm) {$\eta_1$};
\draw[dashed,thick] (13.9 cm,0) -- +(4 cm,0);
\draw[thick] (18.3 cm,0) circle (.5cm);
\node [below] at (18.3 cm,-.5cm) {$\eta_{b-1}$};
\draw[thick] (18.8 cm,0) -- +(2.2 cm,0);
\draw[thick] (21.5 cm,0) circle (.5cm);
\node [below] at (21.5 cm,-.5cm) {$\eta_b$};
\end{tikzpicture}
\end{center}
where
\begin{enumerate}
  \item[(a)] if $\ell$ is linear, then $f=d=a=b$ and $\rho_k$ (resp. $\eta_k$) is obtained by
  adding a $d$-hook to $X$ (resp. $Y$) for each $k=1,\dots,d$.
  \item[(b)] if $\ell$ is unitary, then $f=2d$,  $\rho_1,\dots,\rho_a$ are the unipotent characters
  corresponding to the symbols obtained by adding a $d$-cohook which increases $|X|$ and $\eta_1,\dots,\eta_b$ are
  the unipotent characters corresponding to the symbols obtained by adding a $d$-cohook which increases $|Y|$.
  In addition, $a=d+D$ and $b=d-D$ where $D$ is the defect of the symbol.
\end{enumerate}

\medskip

\subsection{The $\frakg_{\infty}$-representation on $\scrU_K$}\hfill\\

We show in this section that the representation datum on $KG\mod$ yields
a categorical action of $\frakg_{\infty}$ on $\scrU_K$. 

\subsubsection{The ramified Hecke algebra}
Let $r = t(t+1)$ and $n = r+m$ with $t,m \geqslant 0$. Recall that the inflation
from $G_r$ to $L_{r,1^m}$ yields an equivalence between $KG_r\umod$ and
$KL_{r,1^m}\umod$ which intertwines the functor $F_{n,r}$ with the parabolic
induction $R_{L_{r,1^m}}^{G_n}$. In particular, we have a canonical isomorphism
  $$\scrH(KG_n,E_t):= \End_{KG_n}(F^m(E_t))^\op \simto 
  \End_{K G_n}(R^{G_n}_{L_{r,1^m}}(E_t))^\op.$$
This algebra is isomorphic to the Hecke algebra of type $B_m$
with parameters $q$ and $Q_t$, giving a parametrisation of the constituents of (the head of) $F^m(E_t)$
in terms of bipartitions. More precisely, recall from \S \ref{sec:categoricalactions}
that to the categorical datum $(E,F,X,T)$ is attached a map $\phi_{F^m}:\bfH_{K,m}^{q}\to\End(F^m)$.
The evaluation of this map at the module $E_t$ yields a $K$-algebra homomorphism
  $$\phi_{K,m}:\bfH^{q}_{K,m}\to\scrH(KG_n,E_t),\quad 
  X_k\mapsto X_k(E_t),\quad T_l\mapsto T_l(E_t).$$
We show that it induces the aforementioned isomorphism.

\begin{theorem}\label{thm:HL-BC}
Let $t,m \geqslant 0$ and $n = t(t+1)+m$. Then the map $\phi_{K,m}$ factors through
a $K$-algebra isomorphism $\bfH^{Q_t\,;\,q}_{K,m}\mathop{\longrightarrow}\limits^\sim\scrH(KG_n,E_t).$
\end{theorem}

\begin{proof} The proof is similar to the proof of Theorem \ref{thm:cat}, see \cite[thm.~4.12]{DVV}. 
Write $Q_t=(Q_1,Q_2)$ and $X=X(E_t)$. We must check that the operator $X$ on $F(E_t)$
satisfies the relation 
$$(X-(-q)^{-1-t})(X-(-q)^{t})=0.$$
Then the invertibility of the morphism $\bfH^{Q_t\,;\,q}_{K,m}\to\scrH(KG_n,E_t)$
follows from \cite{HL80}. In fact, it is shown there
that $X$ satifies the relation
\begin{equation} \label{eq:xsign} 
(X-\epsilon_t\,(-q)^{-1-t})(X-\epsilon_t\,(-q)^{t})=0
\end{equation}
for some $\epsilon_t = \pm 1$. Therefore we must show that $\epsilon_t=1$ for all $t \geq 0$,
which we will do by induction on $t$. First observe that 
the eigenvalues of $X_{1,0}$ on $R_{G_0}^{G_1} (K)$ are $1$, $(-q)^{-1}$ and
thus are powers of $-q$. Now fix $t \geq 1$ and assume that for all $t>s\geqslant 0$
the eigenvalues of $X(E_s)$ on $F (E_s)$ are powers of $-q$. We will
show that the eigenvalues of $X(E_t)$ are also powers of $-q$ using
the modular representation theory of $G_n$. 

\smallskip

Recall that $K$ is chosen with respect to an $\ell$-modular system $(K,\scrO,\k)$.
Since the parametrization of unipotent characters does not depend on $\ell$
and $(E,F,X,T)$ are defined over $\bbZ[1/q(q-1)]$, we can first choose a specific prime
number $\ell$ and prove that the eigenvalues of $X(E_t)$ are powers of $-q$ modulo $\ell$. 
We choose $\ell$ to be odd and such that the order of $q$ in $\k^\times$ is $f:=4t$. 
Thus the order of $q^2$ is $d = 2t$. Note that $\ell$ is a unitary prime.

\smallskip

Set $r=t(t+1)$ and $n=r+1$. The cuspidal representation $E_t$ is attached to
the symbol 
$$ \Theta_t \, = \, \left\{ \begin{array}{ccccccc} t & t-1 & t-2 & \cdots & -t & -t-1 & \cdots \\ 
				 & & & & & -t-1 & \cdots \\ \end{array} \right\}$$
Since only one $d$-cohook can be removed from $\Theta_t$, the $\ell$-block of $E_t$ has
cyclic defect groups. Moreover, the $d$-cocore of $\Theta_t$ equals $\Theta_{t-1}$. 
Consequently, the Brauer tree of the $\ell$-block of $E_t$ is
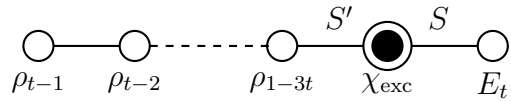
\begin{center}
\begin{tikzpicture}[scale=.4]
\draw[thick] (-1.7 cm,0) circle (.5cm);
\node [below] at (-1.7 cm,-.5cm) {$\rho_{t-1}$};
\draw[thick] (-1.2 cm,0) -- +(2.2 cm,0);
\draw[thick] (1.5 cm,0) circle (.5cm);
\node [below] at (1.5 cm,-.5cm) {$\rho_{t-2}$};
\draw[dashed,thick] (2 cm,0) -- +(4 cm,0);
\draw[thick] (6.4 cm,0) circle (.5cm);
\node [below] at (6.4 cm,-.5cm) {$\rho_{1-3t}$};
\draw[thick] (6.9 cm,0) -- +(2.2 cm,0);
\draw[thick,fill=black] (9.9 cm,0) circle (.5cm);
\node [above] at (8.3 cm,.2cm) {$S'$};
\draw[thick] (9.9 cm,0) circle (.8cm);
\node [below] at (9.9 cm,-.55cm) {$\chi_{\text{exc}}$};
\draw[thick] (10.7 cm,0) -- +(2.2 cm,0);
\node [above] at (11.6 cm,.2cm) {$S$};
\draw[thick] (13.4 cm,0) circle (.5cm);
\node [below] at (13.4 cm,-.5cm) {$E_t$};

\end{tikzpicture}
\end{center}
where a symbol $\Xi_k$ of $\rho_k$ is obtained by adding the $d$-cohook $(k,k+d)$ to $\Theta_{t-1}$
for $k \in \{1-3t,\ldots,t-1\}$. Explicitly, we have
\begin{align*}
\Xi_k\, = \, \left\{ \begin{array}{ccccccc} t-1 & t-2 & \cdots & \widehat{k} &  \cdots  & \cdots\\ 
				 & & k+2t &  -t & -t-1 & \cdots \\ \end{array} \right\}.
\end{align*}
Here the notation $\widehat k$ means that the integer $k$ has been removed. This symbol has defect
$|2t-3| = |2(t-2)+1|$. Therefore by Proposition \ref{prop:HC-B}, the unipotent characters $\rho_k$ all lie in the
Harish-Chandra series above $E_{t-2}$ (recall the convention $E_{-1} := E_0$). Furthermore, the
bipartition $\mu_k$ such that $\Theta_{t-2}(\mu_k) = \Xi_k$ is $((1^{t-1-k}),(k+3t-1))$,
except when $t =1$  in which case $\mu_k = ((k+2),(1^{-k}))$.

\smallskip

The $KG_r$-module $E_t$ is cuspidal, hence weakly cuspidal. It follows that all the composition
factors of any $\ell$-reduction of $E_t$ are weakly cuspidal as well. In particular $$E(S) =0.$$
For each $k$ we can compute the character $E(\rho_k)$ obtained by removing a 1-hook from $\Xi_k$. 
Two cases arise: if $k=1-3t$ or $k=t-1$, then $E(\rho_k)$ is irreducible,
and the corresponding symbols are equal respectively to
$$\begin{aligned}
& \left\{\begin{array}{ccccccc} t-1 & t-2 & \cdots & \widehat{-3t+2} &  \cdots  & \cdots\\ 
				 & & -t+1 &  -t & -t-1 & \cdots \\ \end{array} \right\}, \\
\hskip-2.5cm \text{and} \hskip3cm & \left\{ \begin{array}{ccccccc} t-2 & \cdots & \cdots  & \cdots & \cdots\\ 
				 & 3t-2 &  -t & -t-1 & \cdots \\ \end{array} \right\};
\end{aligned}$$
otherwise $E(\rho_k)$ has two constituents whose symbols are
$$\begin{aligned}
& \left\{ \begin{array}{ccccccc} t-1 & t-2 & \cdots & \widehat{k+1}   & \cdots\\ 
				 & k+2t &  -t & -t-1 & \cdots \\ \end{array} \right\},\\
\hskip-4cm \text{and} \hskip3.1cm & \left\{ \begin{array}{ccccccc} t-1 & t-2 & \cdots & \widehat{k} &   \cdots\\ 
				 &  k+2t-1 &  -t & -t-1 & \cdots \\ \end{array} \right\}.
\end{aligned}$$
We deduce that $\sum (-1)^k [E(\rho_k)] = 0$ in $[K G_{r-1}\mod]$.
Since $[S'] = \sum (-1)^k d_\scrU([\rho_k])$ in $[\k G_{r-1}\mod]$, this implies that $$E(S')=0.$$
Therefore, the two composition factors $S$, $S'$ of the 
$\ell$-reduction of the exceptional characters are weakly
cuspidal, which forces the exceptional characters to be weakly cuspidal as well.

\smallskip

Given a symbol $\Theta$, the module $F(E_\Theta)$ is the sum of the unipotent characters
associated to the symbols obtained from $\Theta$ by adding a 1-hook. 
Thus, we have $F(E_t)=E_{\Xi}\oplus E_{\Xi'}$, where
$$\begin{aligned}
\Xi = &\, \left\{\begin{array}{ccccccc} t & t-1 & \cdots & -t & -t-1 & -t-2 &  \cdots \\ 
				 & & & & -t& -t-2 & \cdots \\ \end{array} \right\},\\
\Xi'= &\, \left\{\begin{array}{ccccccc} t+1 & t-1 & t-2 & \cdots  & -t-1 & \cdots \\ 
				 & & & & -t-1 & \cdots \\ \end{array} \right\}.
\end{aligned}$$
Let $B$, $B'$ be the $\ell$-blocks of $E_{\Xi}$ and $E_{\Xi'}$ respectively and $b,$ $b'$ be the corresponding idempotents in $K G_n$. Since the $d$-cocores of $\Xi$ and $\Xi'$ are
different, the idempotents $b$ and $b'$ are orthogonal. Moreover $X(E_t)$ has the eigenvalues
$\epsilon_t(-q)^{-1-t}$ and $\epsilon_t(-q)^{t}$, on $bF(E_t)=E_{\Xi}$ and $b'F(E_t)=E_{\Xi'}$ respectively. 
Note that $\epsilon_t(-q)^{-1-t}$ and $\epsilon_t(-q)^{t}$ are not congruent modulo $\ell$
since $q^{2t+1} \equiv -q$. 

\smallskip

Let $\chi$ be an exceptional character. 
Recall that $S$ is isomorphic to any $\ell$-reduction of $E_t$. Since it is a composition
factor of the $\ell$-reduction of $\chi$, we deduce that both $bF(\chi)$ and $b'F(\chi)$ are
non-zero. Since $\chi$ is weakly cuspidal, by the Mackey formula $F(\chi)$ has at most two
irreducible constituents and $X(\chi)$ at most two eigenvalues whose product equals $-q^{-1}$.
We deduce that $F(\chi) = bF(\chi) \oplus b'F(\chi)$. In addition, the eigenvalues of
$X(\chi)$ must be congruent to the eigenvalues of $X(S)$ on $F(S)$, which are equal to
$\epsilon_t(-q)^{-1-t}$
and $\epsilon_t(-q)^{t}$. Now, $S'$ is a composition factor of the $\ell$-reduction of
$\chi$ so one of $bF(S')$ or $b'F(S')$ must be non-zero, and therefore $X(S')$ must have an eigenvalue
congruent to $\epsilon_t(-q)^{-1-t}$ or $\epsilon_t(-q)^{t}$. To conclude that $\epsilon_t=1$,
we must compute the possible eigenvalues of $X(\rho_{1-3t})$ on $bF(\rho_{1-3t})$ 
and $b'F(\rho_{1-3t})$ and use that $S'$ is a composition factor of the $\ell$-reduction of
$\rho_{1-3t}$.

\smallskip

We abbreviate $\rho = \rho_{1-3t}$. The Harish-Chandra induction of $\rho$ is $F(\rho)=E_\Theta\oplus E_{\Theta'}\oplus E_{\Theta''}$ where
\begin{align*}
\Theta\, &= \, \left\{ \begin{array}{ccccccc} t & t-2 & t-3 & \cdots & \widehat{1-3t} &  \cdots  \\ 
				 & & & -t+1 &  -t &  \cdots \\ \end{array} \right\},\\
\Theta'\, &= \, \left\{\begin{array}{ccccccc} t-1 & t-2 & \cdots & \widehat{-3t} &  \cdots  & \cdots\\ 
				 & & -t+1 &  -t & -t-1 & \cdots \\ \end{array} \right\},\\
\Theta''\, &= \, \left\{ \begin{array}{ccccccc} t-1 & t-2 & \cdots & \widehat{1-3t} &  \cdots  & \cdots\\ 
				 & & -t+2 &  -t & -t-1 & \cdots \\ \end{array} \right\}.
\end{align*}
We first observe that the modules $E_\Theta$, $E_{\Theta'}$ belong to the $\ell$-blocks $B$, $B'$
and that $E_{\Theta''}$ is projective. Indeed, the $d$-cocores of the symbols $\Theta$, $\Theta'$ and
$\Theta''$ are different and equal respectively to
\begin{align*}
\left\{ \begin{array}{ccccccc} t & t-2 & t-3 & \cdots & \cdots  \\ 
				 & & & -t &  \cdots \\ \end{array} \right\},\quad
\left\{\begin{array}{ccccccc} t-1 & t-2 & \cdots & \cdots &  \cdots  \\ 
				 & &  -t+1 &  -t-1 &  \cdots \\ \end{array} \right\},\quad
\Theta''.
\end{align*}
Since $\rho$ belong to the Harish-Chandra series above $E_{t-2}$, we can compute the eigenvalues of
$X(\rho)$ using the eigenvalues of $X(E_{t-2})$. Let $n'=r+m-4t+2$ for some integer $m\geqslant 0$. 
If $t \geq 2$ we can use the induction hypothesis and the map $\phi_{K,m}$ yields a $K$-algebra isomorphism $$\bfH^{Q_{t-2},\,q}_{K,m}\simto\scrH(K G_{n'},E_{t-2}).$$
By the discussion after \eqref{bijectionB}, the corresponding bijection
$$\Irr(KG_{n'},E_{t-2})\mathop{\longleftrightarrow}\limits^{1:1}\Irr(\bfH^{Q_{t-2},\,q}_{K,m})$$ takes the module $E_{\Theta_{t-2}(\mu)}$
to $S(\mu)_K^{Q_{t-2},\,q}$ for each 2-partition $\mu$ of $m$.
Under this parametrization, the character $\rho$ of $KG_r$ and the characters $E_\Theta$, $E_{\Theta'}$
and $E_{\Theta''}$ of $KG_n$ are mapped to the modules
$$S(\lambda)_K^{Q_{t-2},\,q},\quad
S(\mu)_K^{Q_{t-2},\,q},\quad
S(\mu')_K^{Q_{t-2},\,q},\quad
S(\mu'')_K^{Q_{t-2},\,q}$$ 
labelled by the following 2-partitions 
$$\lambda=((1^{4t-2}),\emptyset),\quad
\mu=((21^{4t-3}),\emptyset),\quad
\mu'=((1^{4t-1}),\emptyset),\quad
\mu''=((1^{4t-2}),(1)).$$
The $(Q_{t-2}, q)$-shifted residue of the boxes
$Y(\mu)\setminus Y(\lambda)$ and $Y(\mu')\setminus Y(\lambda)$
are $(-1)^t q^{t-1}$ and $(-1)^t q^{-3t}$ respectively. They are congruent to $(-q)^{-1-t}$ and $(-q)^t$ modulo $\ell$, 
because $q^{2t}$ is congruent to $-1$ modulo $\ell$. We deduce that the eigenvalues of
the operator $X(\rho_{1-3t})$ on $E_\Theta$, $E_{\Theta'}$ are congruent to $(-q)^{-1-t}$ and
$(-q)^t$ modulo $\ell$. By the above argument, at least one of these must be congruent to
$\epsilon_t (-q)^{-1-t}$ or $\epsilon_t (-q)^t$, which forces $\epsilon_t =1$. 

\smallskip
Finally, if $t <2$ (in which case $t=1$) we use the fact that $E_{t-2} =E_{1-t}$
and $Q_{t-2} = Q_{1-t}^\dag$. In particular, the bijection
$\Irr(KG_{n'},E_{t-2})\to\Irr(\bfH^{Q_{1-t},\,q}_{K,m})$ takes the module
$E_{\Theta_{t-2}(\mu)}$ to $S(\mu^\dag)_K^{(Q_{1-t},q)}$ for any bipartition $\mu$
of $m$. With $t=1$, we deduce that the unipotent characters $\rho$, 
$E_\Theta$, $E_{\Theta'}$ and $E_{\Theta''}$ correspond to Specht modules
of $\bfH^{Q_{0},\,q}_{K,m}$ labelled by the 2-partitions
$\lambda=(\emptyset,(1^2))$, $\mu=(\emptyset,(21))$, $\mu'=(\emptyset,(1^{3}))$ and
$\mu''=((1),(1^2)).$ With $Q_0 = (1,-q^{-1})$, 
the $(Q_{0}, q)$-shifted residues of the boxes
$Y(\mu)\setminus Y(\lambda)$ and $Y(\mu')\setminus Y(\lambda)$
are $-1$ and $-q^{-3}$ respectively, which are congruent to $q^{-2}$ and $-q$ since 
$q^{2}$ is congruent to $-1$ modulo $\ell$. Therefore $\epsilon_{1} = 1$. 
\end{proof}

\smallskip

\subsubsection{The $\frakg_\infty$-representation on $\scrU_K$} 
By Proposition \ref{prop:RD}, the Harish-Chandra induction and restriction yield a representation datum $(E,F,X,T)$ on the category
$$\scrU_K=\bigoplus_{n\in\bbN}K G_n\umod.$$
For any $t,m,n\in\bbN$, let $(KG_n,E_t)\mod$ be the Serre subcategory of $\scrU_K$ generated by the modules
$F^m(E_t)$ with $n=r+m$ and $r=t(t+1)$. Write $$\scrU_{K,t}:=\bigoplus_{n\geqslant 0}(KG_n,E_t)\mod.$$
We have $\Irr((KG_n,E_t)\mod)=\Irr(KG_n,E_t)$ and $\scrU_K=\bigoplus_{t\geqslant 0}\scrU_{K,t}.$ 

\smallskip

\begin{definition}\label{def:infty} 
Let $\scrI_\infty$ denote the subset $q^{\,\bbZ} \sqcup \big(-q^{\,\bbZ}\big)$ of $K^\times.$
We define $\frakg_\infty$ to be the (derived) Kac-Moody algebra associated to the quiver $\scrI_\infty(q)$. 
\end{definition}

\smallskip

To avoid cumbersome notation, we will write $I_\infty=I_\infty(q)$.
We denote by $\{\Lambda_i\}$, $\{\alpha_i\}$ and $\{\alpha_i^\vee\}$ the fundamental
weights, simple roots and simple coroots of $\frakg_\infty$. Then $\X_\infty=\P_{\infty} = \bigoplus \bbZ \Lambda_i$. 
There is a Lie algebra isomorphism $(\fraks\frakl_\bbZ)^{\oplus 2} \simto \frakg_{\infty}$ such that 
$(\alpha^\vee_k,0)\mapsto\alpha^\vee_{q^{k}}$ and $(0,\alpha^\vee_k)\mapsto\alpha^\vee_{-q^{k}}$.
Since the pair $Q_t$ belongs to $(\scrI_\infty)^2,$ the $\frakg_\infty$-module $\bfF(Q_t)_\infty$ is well-defined.

\smallskip

\begin{theorem}\label{thm:ginfinityB}
 Let $t \geqslant 0$ and $Q_t$ be as in \eqref{eq:chargeBC}. 
 \begin{itemize}[leftmargin=8mm]
  \item[$\mathrm{(a)}$] 
  The Harish-Chandra induction and restriction functors yield a representation of
  $\frakg_{\infty}$ on $\scrU_{K,t}$ which is isomorphic to $\scrL(\Lambda_{Q_t})_\infty$. 
  \item[$\mathrm{(b)}$] The map $|\mu,Q_{t}\rangle_\infty\mapsto [E_{\Theta_t(\mu)}]$ gives a
  $\frakg_\infty$-module isomorphism $\bfF({Q_t})_{\infty} \simto [\scrU_{K,t}].$
 \end{itemize}
\end{theorem}

\begin{proof}
The proof is similar to the proof of Theorem \ref{thm:char0} in \cite[thm.~4.15]{DVV}.
It is a consequence of Theorem \ref{thm:HL-BC} and the discussion after \eqref{bijectionB}.
\end{proof}

\medskip

\subsection{The $\frakg_{2d}$-representation on $\scrU_\k$\label{sec:uk-typeBC}}\hfill\\

We now turn to the case of unipotent representations in positive characteristic.
We first deduce from the previous construction a categorical action of the derived
Kac-Moody algebra $\frakg_{2d}'$ on $\scrU_\k$. We then show how to extend this action depending
on the parity of $f$. This reflects the difference between the case of linear primes
($f$ odd) and unitary primes ($f$ even). 

\smallskip

\subsubsection{The $\frakg'_{2d}$-representation on $\scrU_\k$\label{ssec:cataction2d}}
Recall that $(K,\scrO,\k)$ is an $\ell$-modular system with $K \supset \bbQ_\ell$ 
and $\k \supset \overline{\bbF}_\ell$. We still work under
the assumption that $\ell\nmid q$ and $\ell$ and $q$ are odd. Recall that
$d$ (resp. $f$) is the order of $q^2$ (resp. $q$) modulo $\ell$.

\begin{definition} Let $\scrI_{2d}$ be the quiver obtained from $\scrI_\infty$ by specialization $\scrO \to \k$.
We define $\frakg'_{2d}$ to be the derived Kac-Moody algebra associated to the quiver $\scrI_{2d}$. 
\end{definition} 

If $f$ is odd then $f=d$ and $-1$ cannot be expressed as a power of $q$ in $\k$. In that case
the quiver $\scrI_{2d}$ decomposes as $I_{2d}=I_{d,1}\,\sqcup\, I_{d,2}$ with
$I_{d,1}=q^\bbZ$ and $I_{d,2}=-q^{\bbZ}$ being both cyclic quivers of size $d$.
This yields a Lie algebra isomorphism 
$$
(\widetilde{\fraks\frakl}_d)^{\oplus 2}\simeq\frakg'_{d,1}\oplus\,\frakg'_{d,2}=\frakg'_{2d}$$
such that $(\alpha^\vee_k,0)\mapsto\alpha^\vee_{q^{k}}$ and $(0,\alpha^\vee_k)\mapsto\alpha^\vee_{-q^{k}}$.

\smallskip

If $f$ is even, then $f=2d$ and $q^d = -1$. Hence $I_{2d} = q^{\bbZ}$ is a cyclic quiver
of size $2d$ and we have an isomorphism $\widetilde{\fraks\frakl}_f\simeq\frakg'_{2d}$
such that $\alpha^\vee_k\mapsto\alpha^\vee_{q^{k}}$.

\smallskip

The specialization from $\scrO\subset K$ to $\k$ yields a morphism of quivers
$\sp\,:\,\scrI_\infty\to\scrI_{2d}$ and a morphism of abelian groups $\P_\infty\to\P_{2d}$
such that $\Lambda_i\mapsto\Lambda_{\sp(i)}$. The infinite sums $E_i=\bigoplus_{\sp(j)=i}E_j$
and $F_i=\bigoplus_{\sp(j)=i}F_j$ give well-defined operators on $\bfF({Q_t})_{\infty}$.
This yields a representation of $\frakg'_{2d}$ on $\bfF({Q_t})_{\infty}$ such that
the linear map $$\sp\,:\,\Res_{\frakg'_{2d}}^{\frakg_\infty}\bfF({Q_t})_{\infty}\to\bfF({Q_t})_{2d}$$ given by
$|\mu,Q_{t}\rangle_\infty\mapsto|\mu,Q_t\rangle_{2d}$ is a $\frakg'_{2d}$-equivariant isomorphism.

\smallskip

Under the map $d_\scrU:[\scrU_K]\to[\scrU_\k]$ and the isomorphism
$\bigoplus_{t\in\bbN}\bfF({Q_t})_{\infty} \simto [\scrU_{K}]$ in Theorem \ref{thm:ginfinityB}, the map $\sp$ endows $[\scrU_\k]$ with
a representation of $\frakg'_{2d}$ which is compatible with the representation associated with
the representation datum. More precisely, the following hold.

\begin{proposition}\label{prop:lreductionB}
For each $i \in \scrI_{2d}$, let $\k E_i$ and $\k F_i$
be the generalized $i$-eigenspace of $X$ on $\k E$ and $\k F$. Then
  \begin{itemize}[leftmargin=8mm]
    \item[$\mathrm{(a)}$] $[\k E_i],$ $ [\k F_i]$ endow $[\scrU_\k]$ with a structure
    of $\frakg'_{2d}$-module,
    \item[$\mathrm{(b)}$] $d_\scrU$
     yields a $\frakg'_{2d}$-module isomorphism 
     $ \mathrm{Res}_{\frakg'_{2d}}^{\frakg_\infty} \, [\scrU_K]\simto [\scrU_\k],$
      \item[$\mathrm{(c)}$] the map
     $|\mu, Q_t\rangle_{2d} \mapsto d_\scrU([E_{\Theta_t(\mu)}])$ yields a
$\frakg'_{2d}$-module isomorphism
$$\bigoplus_{t\in\bbN}\bfF({Q_t})_{2d} \simto [\scrU_\k].$$
     \qed
  \end{itemize}
\end{proposition}

\begin{proof}
The proposition is a direct application of Theorem \ref{thm:ginfinityB}, once we have proved that the decomposition map
$d_\scrU$ is a vector space isomorphism.
This is known by Theorem \ref{prop:isodecmap} since $\ell$ is odd.
\end{proof}

\smallskip

In order to prove that this, and the representation datum $(E,F,X,T)$ on 
$\scrU_\k$ introduced in \S \ref{sec:repdatumB},
endow $\scrU_\k$ with a $\frakg'_{2d}$-representation,
it remains to see that weight spaces are sums of blocks. This follows from the following
lemma.

\smallskip

\begin{lemma} Let $t,s \geqslant 0$ and $\mu,\lambda$ be bipartitions. If $E_{\Theta_t(\mu)}$
and $E_{\Theta_s(\lambda)}$ are in the same $\ell$-block then $|\mu,Q_t\rangle_{2d}$ and
$|\lambda,Q_s\rangle_{2d}$ have the same weight for the action of $\frakg'_{2d}$.
\end{lemma}

\begin{proof} 
Recall that if $f$ is even (resp. odd), two unipotent characters of $G_n$ lie in the
same $\ell$-block if and only if the corresponding symbols have the same $d$-cocore (resp.
$d$-core), see \S \ref{sec:blockB}. Therefore it is enough to show that adding/removing $d$-cohooks (resp. $d$-hooks)
does not affect the weight of $|\mu,Q_t\rangle_{2d}$ for the action of $\frakg'_{2d}$.

\smallskip

Write $Q_t = (Q_1,Q_2)$ and $\mu = (\mu^1,\mu^2)$. Using Proposition \ref{prop:tensorfock} and 
\eqref{eq:weightlevel1} one can compute the
weight of $|\mu,Q_t\rangle_{2d}$ from the $\beta$-set $\beta_0(\mu)$. 
For $i \in I_{2d}$ and $N$ large enough, we have
 $$\begin{aligned}\begin{split}
  \langle \mathrm{wt}\big(|\mu,Q_t\rangle_{2d} \big) ,\alpha_i^\vee \rangle =\,\,&
 |\{\beta \in \beta_0(\mu^1), \, \beta \geqslant -N \, | \,q^{\beta}Q_1\equiv i \text{ mod $\ell$}\}| \\ & - 
 |\{\beta \in \beta_0(\mu^1),\, \beta \geqslant -N+1 \, | \,q^{\beta}Q_1\equiv q i \text{ mod $\ell$}\}| \\
 & + |\{\beta \in \beta_0(\mu^2), \, \beta \geqslant -N \, | \,q^{\beta}Q_2\equiv i \text{ mod $\ell$}\}| \\ & - 
 |\{\beta \in \beta_0(\mu^2),\, \beta \geqslant -N+1 \, | \,q^{\beta}Q_2\equiv q i \text{ mod $\ell$}\}|.
 \end{split}
\end{aligned}$$
This expression becomes simpler when working with the charged $\beta$-sets
used in the definition of $\Theta_t(\mu)$ in \S \ref{sec:symbol}. Write 
$\Theta_t(\mu) = \{X,Y\}$ where  $X=\beta_t(\mu^1)$ and $Y=\beta_{-1-t}(\mu^2)$. 
Then by definition of $Q_t$ we have, for $N$ large enough 
$$\begin{aligned}\begin{split}
  \langle \mathrm{wt}\big(|\mu,Q_t\rangle_{2d} \big) ,\alpha_i^\vee \rangle =\,\,&
 |\{x \in X, \, x \geqslant -N \, | \,q^{x} \equiv (-1)^ti \text{ mod $\ell$}\}|  \\ & - 
 |\{x \in X,\, x \geqslant -N+1 \, | \,q^{x}\equiv (-1)^tq i \text{ mod $\ell$}\}| \\
 & + |\{y \in Y, \, y \geqslant -N \, | \,q^{y}\equiv (-1)^{t+1} i \text{ mod $\ell$}\}| \\& - 
 |\{y \in Y,\, y \geqslant -N+1 \, | \,q^{y}\equiv (-1)^{t+1}q i \text{ mod $\ell$}\}|.
 \end{split}
\end{aligned}$$

If $f$ is odd, then $q^d = q^f = 1$ in $\k$ and adding or removing a $d$-hook has the
effect of removing $x$ from $X$ (resp. $Y$) and adding $x \pm d$ to $X$
(resp. $Y$), which does not change the congruence of $q^x$. In particular, this operation
does not affect $\mathrm{wt}\big(|\mu,Q_t\rangle_{2d} \big)$.

\smallskip

If $f$ is even, then $q^d = -1$, $q^f = 1$ and adding or removing a $d$-cohook has
the effect of removing $x$ from $X$ (resp. $Y$) and adding $x\pm d$ to $Y$ (resp. $x$).
This changes $q^x$ to $q^{x\pm d} \equiv -q^x$. Therefore if $t' = t \pm 1$ and
$\Theta_{t'}(\mu')$ is the symbol obtained from $\Theta_t(\mu)$ by adding or
removing a $d$-cohook we have
$$\mathrm{wt}\big(|\mu,Q_t\rangle_{2d} \big) ,\alpha_i^\vee \rangle 
= \mathrm{wt}\big(|\mu',Q_{t'}\rangle_{2d} \big) ,\alpha_{-i}^\vee \rangle.$$
In particular, symbols corresponding to unipotent characters in the same $\ell$-block
have the same weight. 
\end{proof}

\smallskip

Notice that the lemma implies that the simple unipotent modules
are weight vectors of the $\frakg'_{2d}$-action.
From this and the results in the previous section we deduce the expected result. 

\smallskip

\begin{theorem}\label{thm:B} Recall that $\ell$ and $q$ are odd, and $\ell \nmid q(q^2-1)$.
The representation datum on $\scrU_\k$
yields a $\frakg'_{2d}$-representation such that
the map $d_\scrU\,:\,[\scrU_K]\to[\scrU_\k]$ intertwines the representations
of  $\frakg_{\infty}$ and $\frakg'_{2d}$.
There is a $\frakg'_{2d}$-module isomorphism
$\bigoplus_{t\in\bbN}\bfF({Q_t})_{2d}\simto[\scrU_\k]$ sending
$|\mu, Q_t\rangle_{2d}$ to $d_\scrU([E_{\Theta_t(\mu)}])$.
Further, the classes in $[\scrU_\k]$ of the simple unipotent modules
are weight vectors for the $\frakg'_{2d}$-action.
\qed
\end{theorem}

\smallskip

\subsubsection{The $\frakg_{2d}$-representation on $\scrU_\k$ in the linear prime case}
In this section we assume that $f$ is odd, so $f=d$. In that case
the Kac-Moody algebra $\widetilde \frakg_{2d}$ associated with the quiver $I_{2d}$
is isomorphic to $(\widehat{\fraks\frakl}_d)^{\oplus 2}$. As in the case of unitary
groups, see \S \ref{subsec:g-e}, the action of $\frakg_{2d}'$ on $\scrU_\k$ can be naturally
extended to an action of an algebra $\frakg_{2d}$ lying between $\frakg_{2d}'$ 
and $\widetilde\frakg_{2d}$. 

\smallskip

Let $\widetilde \X_{2d}$ and $\widetilde \X_{2d}^\vee$ be the lattices corresponding
to $\widetilde\frakg_{2d}$. Since $f$ is odd, $\scrI_{2d}$ is the disjoint union
of two cyclic quivers. We choose $\alpha_1$ and $\alpha_{-1}$ to be the affine roots
attached to each of these quivers. Then we have 
$$\widetilde \X_{2d} = \P_{2d} \oplus
 \bbZ \delta_+ \oplus \bbZ\delta_-,\quad \widetilde \X_{2d}^\vee = \Q_{2d}^\vee
\oplus \bbZ \partial_+ \oplus \bbZ \partial_-,$$ with $\delta_+ = \sum \alpha_{q^{j}}$, 
$\delta_- = \sum \alpha_{-q^{j}}$, $\partial_+ = \Lambda_{1}^\vee$ and
$\partial_- = \Lambda_{-1}^\vee$. We set $\partial=\partial_++\partial_-$ and $\delta=(\delta_++\delta_-)/2$. 
We define $\frakg_{2d}:=\frakg'_{2d}\oplus\bbC\partial$ and we view it as the Kac-Moody algebra associated with the lattices 
$$\X_{2d} := \P_{2d} \oplus \bbZ\delta \simeq \widetilde \X_{2d} / (\delta_+-\delta_-),\quad
\X_{2d}^\vee :=  \Q_{2d}^\vee \oplus \bbZ \partial.$$
The pairing $\widetilde\X_{2d}^\vee\times\widetilde\X_{2d}\longrightarrow\bbZ$ induces
in the obvious way a perfect pairing $\X_{2d}^\vee\times\X_{2d}\longrightarrow\bbZ$.
\smallskip

For $t \in \mathbb{N}$, the Fock space $\bfF(Q_t)_{2d}$ has a tensor product
decomposition into level~$1$ Fock spaces 
$$ \bfF(Q_t)_{2d} \simeq \bfF((-q)^t)_d \otimes \bfF((-q)^{-1-t})_d.$$
Out of the charged Fock spaces $(\bfF((-q)^{t})_d,t)$ and $(\bfF((-q)^{-1-t})_d,-1-t)$ and
the isomorphism $\widetilde \frakg_{2d} \simeq (\widehat{\fraks\frakl}_d)^{\oplus 2}$
(which depends on the parity of $t$) we can therefore equip $\bfF(Q_t)_{2d}$ with an
action of $\widetilde \frakg_{2d}$ which in turn restricts to an action of $\frakg_{2d}$.

\smallskip

Recall that two unipotent characters are in the same $\ell$-block if and
only if the corresponding symbols  have the same $d$-core. In particular,
the unipotent characters of a given unipotent block all lie in the same
Harish-Chandra series. In addition, two unipotent characters lying in
different Harish-Chandra series lie in different blocks. Consequently, for each $t\in\bbN$
we can form the category $\scrU_{\k,t}$ associated with
the Harish-Chandra series labelled by $t$, yielding
$$ \scrU_\k = \bigoplus_{t \in \bbN} \scrU_{\k,t} \quad \text{ with } \quad 
[\scrU_\k] \simeq \bfF(Q_t)_{2d}.$$
Using the action of $\frakg_{2d}$ on $ \bfF(Q_t)_{2d}$ defined above we
equip each $[\scrU_{\k,t}]$ with a structure of $\frakg_{2d}$-module
which extends the structure of $\frakg_{2d}'$-module defined in \S \ref{ssec:cataction2d}. 
The situation is completely similar to \S \ref{subsec:block=weight} and we get the following
theorem.

\begin{theorem}\label{thm:BB}
Recall that $\ell$ and $q$ are odd, and $\ell \nmid q(q^2-1)$.
Assume that $f$ is odd. For each $t \in \bbN$, the Harish-Chandra induction
and restriction functors yield a representation of $\frakg_{2d}$ on $\scrU_{\k,t}$ which categorifies
$\bfF(Q_t)_{2d}$. \qed
\end{theorem}

\smallskip

\subsubsection{Combinatorics of $d$-cohooks and $d$-cocores}
We now concentrate on the unitary prime case, i.e., the case where $f$ is even,
which is the most delicate. We wish to define an action of a (non derived)
Kac-Moody algebra $\frakg_{2d}$ on $[\scrU_\k]$ which extends the action of
$\frakg'_{2d}.$ Then, we will extend the grading of $[\scrU_\k]$ from $\P_{2d}$
to $\X_{2d} = \P_{2d}\, \oplus\, \bbZ \delta/2$. 
Before this, we must introduce some combinatorial tools related to symbols.

\smallskip

Let $\mu = (\mu^1,\mu^2)$ be a bipartition and $t \in \mathbb{Z}$. 
Consider the symbol $\Theta$ defined by 
$$\Theta = \{X\,,\,Y-d\} = \Theta_t(\mu) $$
where $X = \beta_t(\mu^1)$ and $Y=\beta_{d-1-t}(\mu^2)$.
With this notation, removing a $d$-cohook on $\Theta$ changes the pair $(X,Y)$ to
$$(X\smallsetminus\{x\},Y\sqcup \{x\})\quad\text{or}\quad(X\sqcup\{y-2d\},Y\smallsetminus \{y\}).$$

\smallskip 

We denote by $s = (s_1,\ldots,s_{2d})$ and $r = (r_1,\ldots,r_{2d})$ the $2d$-core of the charged partitions
$(\mu^1,t)$ and $(\mu^2,d-1-t)$ respectively, and by $(\mu^{1,1},\dots,\mu^{1,2d})$ and 
$(\mu^{2,1},\dots,\mu^{2,2d})$ their $2d$-quotients.
Hence, we have 
\begin{align*}
\tau_{2d}(\mu^1,t)&=(\mu^{1,1},\dots,\mu^{1,2d},s_1,\ldots,s_{2d}),\\
\tau_{2d}(\mu^2,d-1-t)&=(\mu^{2,1},\dots,\mu^{2,2d},r_1,\ldots,r_{2d}).
\end{align*}

\begin{lemma}\label{lem:A}
The symbol $\Theta$ is a $d$-cocore if and only if 
\begin{itemize}[leftmargin=8mm]
  \item[$\mathrm{(a)}$] 
$\mu^1$ and $\mu^2$ are $2d$-cores,
 \item[$\mathrm{(b)}$]  $r_p -s_p \in \{0,1\}$ for all $p \in \{1,\ldots,2d\}$.
\end{itemize}
\end{lemma}

\begin{proof} Let us write $X = \beta_t(\mu^1)$ and $Y=\beta_{d-1-t}(\mu^2)$.
Then $\Theta$ is a $d$-cocore if and only if $X \subset Y \subset X+2d$. In particular,
we must have $X \subset X+2d$ and $Y \subset Y+2d$, which is equivalent to $\mu^1$ and
$\mu^2$ being $2d$-cores. In that case we have 
$$ X   = \bigsqcup_{p=1}^{2d} (p-2d+2d\beta_{s_p}(\emptyset)) \quad \text{and} \quad
 Y = \bigsqcup_{p=1}^{2d} (p-2d+2d\beta_{r_p}(\emptyset)). 
$$
Using the definition of the $\beta$-set of an empty partition we deduce that
$X \subset Y$ if and only if $s_p \leqslant r_p$ for all $p$. 
For the same reason, we have $Y \subset X+2d$ if and only if $r_p \leqslant s_p + 1$ for all $p$.
The combination of the two conditions gives 
$$r_p - s_p \in \{0,1\},\, \text{ for all } p.$$
\end{proof}

\smallskip

\begin{lemma}\label{lem:B} Let $\Theta$, $\Theta'$ be two symbols. The following assertions are equivalent:
 \begin{itemize}[leftmargin=8mm]
  \item[$\mathrm{(i)}$] there exists a sequence of symbols $\Theta= \Theta_0, \Theta_1, \ldots, \Theta_m = \Theta'$
  where $\Theta_{i+1}$ is obtained from $\Theta_i$ by adding or removing a $d$-cohook,
  \item[$\mathrm{(ii)}$] $s+r=s'+r'$.
 \end{itemize}
\end{lemma}

\begin{proof} 
Let $\mu$ be the bipartition such that $\Theta = \Theta_t(\mu)$.
We write 
$$ \beta_t(\mu^1)   = \bigsqcup_{p=1}^{2d} (p-2d+2d\beta_{s_p}(\mu^{1,p})) \quad  \text{and} \quad
 \beta_{d-1-t}(\mu^2) = \bigsqcup_{p=1}^{2d} (p-2d+2d\beta_{r_p}(\mu^{2,p})). $$
Therefore adding or removing a $d$-cohook changes $(s_p,r_p)$ to $(s_p,r_p)\pm(1,-1)$
for some $p \in \{ 1,\ldots,2d\}$. This does not change $s+r$.

\smallskip
Let $\Theta_\circ$, $\Theta_\circ'$ be the $d$-cocore of $\Theta$, $\Theta'$,
and $(s_\circ,{r}_\circ)$, $(s_\circ',r_\circ')$
be the corresponding $2d$-cores. 
By the previous argument, we have 
$$s+r = s_\circ +{r}_\circ,\quad{s}'+{r}' = {s}_\circ' +{r}_\circ'.$$
Now, assertion (i) is equivalent to $\Theta_\circ = \Theta_\circ'$, which yields
in particular $${s+r} = {s}_\circ +{r}_\circ = {s}'+{r}'.$$
Conversely, if we assume (ii), then ${s}_\circ +{r}_\circ = {s}_\circ' +{r}_\circ'$ which we can write 
$${r}_\circ-{s}_\circ = {r}_\circ'-{s}_\circ'+2({r}_\circ-{r}_\circ').$$
But by Lemma \ref{lem:A}, both ${r}_\circ-{s}_\circ$ and ${r}_\circ'-{s}_\circ'$ are $2d$-tuples of integers equal to $0$ or $1$. 
This forces ${r}_\circ={r}_\circ'$ and therefore ${s}_\circ={s}_\circ'$. 
Since $t_\circ =\sum_p s_{\circ,p},$ this determines $\Theta_\circ$ uniquely and therefore we have
$\Theta_\circ= \Theta_\circ'$. \end{proof}

\smallskip

\subsubsection{The weight of a symbol}
We are still working under the assumption that $f =2d$ is even. In particular
$q^d=-1$ in $\k$. Fix an integer $t\in\bbZ$ and set
$$P_t = q^{c_t} = (q^t\,,\,-q^{-1-t}),\quad c_t = (t,d-1-t).$$
In particular $Q_t = (-1)^t P_t$. Let $\mu$ be a bipartition and let ${s}$ and ${r}$
be the $2d$-cores of the charged partitions $(\mu^1,t)$ and $(\mu^2,d-1-t)$. 
We consider the weight in 
$\text{X}=\text{P}\oplus\bbZ\,\delta/2$ given by
\begin{align}\label{wt}
\begin{split}
\wt(\mu,t) &= \Lambda_{P_t}- \sum_i n_i(\mu,P_t)\, \alpha_i - \nabla({s},{r})\,\delta,\\
\nabla({s},{r}) &= \Delta(t,2d)+\Delta(d-1-t,2d)+t/2,
\end{split}
\end{align}
where $\Delta(t,2d)$ is as in \eqref{Delta}. 
For any $2d$-tuple of integers ${x} = (x_1,\ldots,x_{2d})$, we define as in \S \ref{subsec:xgrading}
$$\pi_{x} = \sum_{p=1}^{2d} (x_p-x_{p+1})\, \Lambda_{q^p}.$$ 
Then, we have
\begin{equation}
\label{eq:classical}
\begin{split}
\wt(\mu,t)&\equiv \Lambda_{P_t}- \sum_i n_i(\mu,P_t)\, \alpha_i^\text{cl}\ \text{mod}\ \bbZ\,\delta, \\
&\equiv \pi_{{s} + {r}} + 2\Lambda_1\ \text{mod}\ \bbZ\,\delta.
\end{split}
\end{equation}
The effect on the weight of removing a $d$-cohook is studied in the following lemma.

\begin{lemma} \label{lem:C} Let $t,$ $t'$ be integers and $\mu,$ $\mu'$ be bipartitions. 
Set $\Theta = \Theta_t(\mu)$ and $\Theta' =\Theta_{t'}(\mu')$. 
If $\Theta'$ is obtained from $\Theta$ by removing a $d$-cohook, then
$\wt(\mu,t) = \wt(\mu',t') - \delta/2.$
\end{lemma}

\begin{proof}
By \eqref{eq:weightlevel1}, we have
$$ \begin{aligned}
\wt(\mu,t) & \, = \wt(|\mu^1,t\rangle_{2d}) + \wt(|\mu^2,d-1-t\rangle_{2d})  
+ \big(\Delta(t,2d)+\Delta(d-1-t,2d)-\nabla({s},{r})\big)\, \delta \\
& \, = \wt(|\mu^1,t\rangle_{2d}) + \wt(|\mu^2,d-1-t\rangle_{2d})  - t\,\delta/2 \\
& \, = \pi_{{s} + {r}} + 2\Lambda_1 - \big(w_{2d}(\mu^1)+w_{2d}(\mu^2)+ 
\Delta({s},1)+\Delta({r},1)+t/2\big)\delta. 
\end{aligned}$$
where we used the formulas in \S \ref{subsec:xgrading}. Therefore we must study the effect of
removing a $d$-cohook on the integer
$$w_{2d}(\mu^1)+w_{2d}(\mu^2)+ \Delta({s},1)+\Delta({r},1)+t/2.$$ 
Note that by Lemma \ref{lem:B}, removing a $d$-cohook has no effect on ${s+r}$. 

\smallskip

First, assume that there is an integer $x \in \beta_t(\mu^1)$ such that $x \notin \beta_{d-1-t}(\mu^2)$
and write
$$\beta_{t'}({\mu'}^1)=\beta_t(\mu^1)\setminus\{x\},\quad
\beta_{d-1-t'}({\mu'}^2)=\beta_{d-1-t}(\mu^2)\sqcup\{x\}.$$
Fix $p \in \{1,\ldots,2d\}$ and $z \in \beta_{s_p}(\mu^{1,p})$ such that
$$x = p-2d+2dz.$$
Fix $u_0$ such that $z = z_{u_0}$ where
$\beta_{s_p}(\mu^{1,p}) = \{z_1 > z_2 > \cdots \}$.
We have 
$$|\mu^{1,p}| = \sum_{u\geqslant 1} (z_u+u-1-s_p).$$
We also have
$$\begin{aligned}
|{\mu'}^{1,p}| - |\mu^{1,p}| &\, 
= \sum_{u\geqslant 1} (z'_u+u-1-s_p+1) - \sum_{u\geqslant 1} (z_u+u-1-s_p)  \\
& \, = \sum_{1 \leq u <u_0} (z_u+u-s_p) + \sum_{u\geqslant u_0} (z_{u+1}+u-s_p) -  \sum_{u\geqslant 1} (z_u+u-1-s_p) \\
&\, = s_p -z. 
\end{aligned}$$
We deduce that we have
$$w_{2d}({\mu'}^1) = w_{2d}(\mu^1)+s_p-z.$$ 
A similar computation yields 
$$w_{2d}(\mu^2) = w_{2d}({\mu'}^2)+r_p'-z = w_{2d}({\mu'}^2)
+r_p+1-z.$$ 
We deduce that
$$ w_{2d}({\mu'}^1) + w_{2d}({\mu'}^2) = w_{2d}({\mu}^1) + w_{2d}({\mu}^2) + s_p-r_p-1.$$
Since we have removed the element $x$ from $\beta_t(\mu^1)$, we have $t' =t-1$.

\smallskip

Next, for any integers $v,$ $e$ with $e \geqslant 1$, we have
\begin{equation*} \Delta(v,e)- \Delta(v-1,e) = \left\lfloor \frac{v-1}{e} \right\rfloor.\end{equation*}
We deduce that
$$\Delta({s}',1) -\Delta({s},1) = 1-s_p,\quad
\Delta({r}',1) -\Delta({r},1) = r_p.$$ In particular, the value of
$$ w_{2d}(\mu^1)+w_{2d}(\mu^2)+ \Delta({s},1)+\Delta({r},1) + t/2$$
decreases by $1/2$.

\smallskip

Now assume that there is an integer $y \in \beta_{d-1-t}(\mu^2)$ such that $y-2d \notin \beta_{t}(\mu^1)$
and let us write
$$\beta_{t'}({\mu'}^1)=\beta_t(\mu^1)\sqcup\{y-2d\},\quad
\beta_{d-1-t'}({\mu'}^2)=\beta_{d-1-t}(\mu^2)\setminus\{y\}.$$ 
Fix $p \in \{1,\ldots,2d\}$ and $z \in \beta_{r_p}(\mu^{2,p})$ such that
$$y = p-2d+2dz.$$
The same computation as above gives
$$w_{2d}({\mu'}^2) = w_{2d}(\mu^2)+r_p-z,\quad
w_{2d}(\mu^1) = w_{2d}({\mu'}^1)+(s_p+1)-(z-1),$$
from which we deduce that
$$ w_{2d}({\mu'}^1) + w_{2d}({\mu'}^2) = w_{2d}({\mu}^1) + w_{2d}({\mu}^2) + r_p-s_p-2.$$
Since $t' = t+1$, the value of
$$ w_{2d}(\mu^1)+w_{2d}(\mu^2)+ \Delta({s},1)+\Delta({r},1) + t/2$$
decreases by $1/2$. 
\end{proof}

\smallskip

\subsubsection{The $\frakg_{2d}$-representation on $\scrU_\k$ in the unitary case}
We now have all the tools to extend the action of $\frakg_{2d}'$ to $\frakg_{2d}$
on $\scrU_\k$. Recall that we assume that $f$ is even, so $f=2d$. Hence $I_{2d}$
is a cyclic quiver of size $2d$ and we have an isomorphism
$\widetilde{\fraks\frakl}_f\simto\frakg'_{2d}$  sending
$\alpha^\vee_k$ to $\alpha^\vee_{q^{k}}$.
Let $\frakg_{2d}$ be the Kac-Moody algebra associated with the quiver $I_{2d}$, i.e., we 
have $\widehat{\fraks\frakl}_f\simeq\frakg_{2d},$ and
set $\X_{2d} = \P_{2d}\, \oplus\, \bbZ \delta/2$.

\smallskip

Let $\sigma$ be the automorphism of the Lie algebra $\frakg_{2d}$ given by the quiver automorphism
of $I_{2d}$ such that $i\mapsto -i$. In particular, we have $\sigma(\alpha_i^\vee)=\alpha_{-i}^\vee$. 

\begin{definition}
We define the weight of a symbol by
$\wt (\Theta_t(\mu)) = \sigma_*^t(\wt (\mu,t))$, where $\wt(\mu,t)$ is as in \eqref{wt}.
\end{definition}

We can now prove the following.

\begin{proposition} Let $\mu$ be a bipartition and $t \in \mathbb{Z}$. Let $\Theta$, $\Theta'$ be symbols.
\hfill  \begin{itemize}[leftmargin=8mm]
\item[$\mathrm{(a)}$] $\wt (\Theta_t(\mu)) \equiv \mathrm{wt}(|\mu,Q_t\rangle_{2d})\ \text{mod}\ \bbZ\,\delta,$
\item[$\mathrm{(b)}$] if $\Theta'$ is
obtained from $\Theta$ by removing a $d$-hook, then $\wt (\Theta) = \sigma_*\!\wt(\Theta') - \delta/2,$
\item[$\mathrm{(c)}$]  $E_\Theta$, $E_{\Theta'}$ 
belong to the same $\ell$-block if and only if we have $\wt (\Theta) = \wt (\Theta')$. 
\end{itemize}
\end{proposition}

\begin{proof}
We have $Q_t = (-1)^tP_t$. Further, 
\begin{equation*}
\begin{split}
\wt(\mu,t)&\equiv \Lambda_{P_t}- \sum_i n_i(\mu,P_t)\, \alpha_i^\text{cl}\ \text{mod}\ \bbZ\,\delta, \\
 \mathrm{wt}(|\mu,Q_t\rangle_{2d})& \equiv \Lambda_{Q_t}-\sum_{i} n_i(\mu,Q_t)\,\alpha_i^\cl\ \text{mod}\ \bbZ\,\delta.
\end{split}
\end{equation*}
This implies part (a) because
$$n_{(-1)^t i}(\mu,P_t)=n_{i}(\mu,(-1)^t P_t)=n_{i}(\mu,Q_t).$$
Part (b) follows from Lemma \ref{lem:C}.

\smallskip

Now, we prove part (c). Recall that the representations $E_\Theta$ and $E_{\Theta'}$ 
are in the same $\ell$-block if and only if $\Theta$ and $\Theta'$ have the same
rank and the same $d$-cocore $\Theta_\circ$, see \S \ref{sec:blockB}.

\smallskip

Assume that this is the case. 
Since the number of $d$-cohooks one
needs to add to $\Theta_\circ$ to obtain $\Theta$ and $\Theta'$ depends only on the rank of the symbols,
we deduce from Lemma \ref{lem:C} that $\Theta$ and $\Theta'$ have the same weight.

\smallskip

Conversely, assume now that $\wt (\Theta) = \wt (\Theta')$. Then \eqref{eq:classical} implies that
\begin{align}\label{form1}\pi_{{s+r}} = \pi_{{s}'+{r}'}.\end{align}
By definition of $s,r,s'$ and $r'$ we have
$$\sum_ps_p=t,\quad\sum_pr_p=d-1-t,\quad\sum_ps'_p=t',\quad\sum_pr'_p=d-1-t',$$ 
hence
\begin{align}\label{form2}\sum_p (s_p+r_p) = \sum_p (s_p'+r_p') = d-1.\end{align}
From \eqref{form1} and \eqref{form2}, we deduce that the $2d$-tuples ${s+r}$ and ${s}'+{r}'$ are equal. 
Hence, by Lemma \ref{lem:B}, 
the symbols $\Theta$ and $\Theta'$ have the same $d$-cocore $\Theta_\circ$. Let $w$ and $w'$
be the number of $d$-cohooks we need to add to $\Theta_\circ$ to obtain $\Theta$ and
$\Theta'$ respectively. Since $\sigma_*(\delta) = \delta$, by part (b) we have
\begin{align*}
\wt (\Theta) &= (\sigma_*)^w\wt (\Theta_\circ) - w\,\delta /2,\\
 \wt (\Theta') &= (\sigma_*)^{w'}\wt (\Theta_\circ) - w'\,\delta /2.
 \end{align*}
Since $\sigma_*(\delta) = \delta$ and $\wt (\Theta) = \wt (\Theta')$, 
we deduce that $w = w'$. In particular, the symbols $\Theta$ and $\Theta'$ have the same rank,
which shows that the corresponding representations lie in the same $\ell$-block.
\end{proof}

We have proved the following refinement of Theorem \ref{thm:B}.

\begin{theorem} \label{thm:BBB}
Recall that $\ell$ and $q$ are odd, and $\ell \nmid q(q^2-1)$.
Assume that $f$ is even.
The representation datum on $\scrU_\k$
gives a $\frakg_{2d}$-representation yielding a $\frakg'_{2d}$-module isomorphism
$\bigoplus_{t\in\bbN}\bfF({Q_t})_{2d}\simto[\scrU_\k]$ such that
$|\mu, Q_t\rangle_{2d}\mapsto d_\scrU([E_{\Theta_t(\mu)}])$.
The weights of the $\frakg_{2d}$-module $[\scrU_\k]$ belong to the lattice $\X_{2d}$.
The weight of $d_\scrU([E_{\Theta_t(\mu)}])$ is equal to $\wt(\Theta_t(\mu))$.
The classes in $[\scrU_\k]$ of the simple unipotent modules
are weight vectors. Two simple unipotent modules belong to the same block
if and only if their classes have the same weight in $\X_{2d}$.
\qed
\end{theorem}

\smallskip

\subsubsection{Determination of the ramified Hecke algebras}
As in the case of finite unitary groups, see \S \ref{sec:crystal-unitary}, we can invoke Proposition \ref{prop:unicite}
to show that every ramified Hecke algebra associated
to a weakly cuspidal module is indeed a Hecke algebra, whose parameters can be
read off from the weight of the cuspidal module.

\begin{corollary}\label{cor:hecke-BC} 
Let $D$ be a weakly cuspidal $\k G_r$-module. Then 
\begin{itemize}[leftmargin=8mm]
\item[$\mathrm{(a)}$] the weight of the class $[D]$ with respect to the $\frakg'_{2d}$-action on $[\scrU_\k]$ is of the form 
$\Lambda_{Q}$ for some element
$Q = (Q_1,Q_2)$ in $I_{2d}^2$,
\item[$\mathrm{(b)}$] for each $m \geqslant 0$ and $n = r+2m$,
the map $\phi_{\k,m}\, :\, \bfH^{q^2}_{\k,m}\longrightarrow\End_{\k G_n}(F^m)^\op$ factors to an algebra isomorphism
  $\bfH^{Q;\,q^2}_{\k,m}\mathop{\longrightarrow}\limits^\sim \scrH(\k G_n,D).$
  \qed
\end{itemize}
\end{corollary}

\smallskip

In order to apply this corollary efficiently, we must determine explicitly the weakly cuspidal modules and their weights.
In the linear prime case, this is known by \cite{GrH}.
In the unitary prime case, by Theorem \ref{thm:BBB}, the weakly cuspidal modules form a basis of the subspace of highest weight vectors
of the representation of $\frakg'_{2d}$ on $[\scrU_\k]$, hence all ramified Hecke algebras can be explicitly computed.
If, moreover, we assume that the decomposition matrix is unitriangular with respect
to the $a$-function, see Conjecture \ref{conj:unitriangular} below, then the weakly cuspidal modules and their weights
can also be read off from the crystal 
isomorphism in Theorem \ref{thm:crystalBC} below. 

\smallskip

More precisely, any weakly cuspidal module is of the form $D = D_\Theta$
for some symbol $\Theta=\Theta_t(\mu)$. The parameter $Q$ can be explicitly computed from $\mu$ and $t$ 
according to the 
rule
\begin{align*}
\mathrm{wt}(D_\Theta)
&\equiv \Lambda_{Q_t}-\sum_{i} n_i(\mu,Q_t)\,\alpha_i^\cl\ \text{mod}\ \bbZ\,\delta,\\
&\equiv \Lambda_{Q_1}+\Lambda_{Q_2}\ \text{mod}\ \bbZ\,\delta.
\end{align*}
Since the isomorphism in Proposition \ref{prop:lreductionB} identifies the set 
$$\{[D_{\Theta_t(\mu)}]\,;\,\mu\in\scrP^2,\,D_{\Theta_t(\mu)}\ \text{is\ weakly\ cuspidal}\}$$
with a basis of the space of highest weight vectors
$$\{x\in\bfF(Q_t)_{2d}\,;\,E_i(x)=0,\,\forall i=1,2,\dots,2d\},$$
the weakly cuspidal modules correspond precisely to the elements in the crystal $B(s_t)_{2d}$ which are killed by
the operators $\widetilde E_1,\dots, \widetilde E_{2d}$, with the notation in \S \ref{sec:crystalBC} below.
These elements have been computed recently in \cite[thm.~5.9]{JL} in terms of bipartitions $\mu$ whose
abacus is totally $2d$-periodic.
N

\medskip

\subsection{Derived equivalences}\hfill\\

Blocks of $\scrU_\k$ are parametrized by pairs $(\Theta,w)$ where $\Theta$ is
a $d$-core (resp. a $d$-cocore) when
$f$ is odd (resp. $f$ is even) and $w$ is a non-negative integer called the
\emph{degree} of the block. The unipotent characters lying in $B_{\Theta,w}$ are 
parametrized by the symbols which are obtained from $\Theta$ by adding 
$w$ successive $d$-hooks (resp. $d$-cohooks) if $f$ is odd (resp. $f$ is even).

\smallskip

Our categorification theorems (see Theorem \ref{thm:BBB} and \ref{thm:BB}) show that blocks
of $\scrU_\k$ and $\scrU_{\k,t}$ correspond to weight spaces for the action of
the Kac-Moody algebra $\frakg_{2d}$. Using the work of Chuang-Rouquier \cite{CR}
(see also Theorem \ref{thm:reflection}) we obtain derived equivalences between blocks of $\scrU_\k$
which are in the same orbit under the action of the affine Weyl group $W_{2d}$.

\smallskip

If $\ell$ is a linear prime (i.e. $f$ is odd), the situation is similar to the case of
finite unitary groups. The $d$-cores, which are formed by pairs of
partitions which are both $d$-cores, form a single orbit under the action of $W_{2d}$. 
In particular, orbits of blocks of $\scrU_{\k,t}$ under $W_{2d}$ are parametrized
by their degree. Let $w \in \bbN$. Livesey showed in \cite[thm.~7.1]{Li12} how to construct
a \emph{good} block of $\scrU_{\k,t}$ of degree $w$ which is derived equivalent to its
Brauer correspondent. We deduce that Brou\'e's abelian defect group conjecture holds for groups of type $BC$ at odd linear primes.

\begin{theorem}\label{thm:derivedBC}
Recall that $\ell$ and $q$ are odd. Assume that the order $f$ of $q$ modulo $\ell$ 
is even. Let $B$ be a unipotent block of $G_n$ over $\k$ or $\scrO$, and $D$ be a defect group of $B$. If $D$ is abelian, then
$B$ is derived equivalent to its Brauer correspondent in $N_{G_n}(D)$. \qed
\end{theorem}

\medskip

\subsection{The crystal of $\scrU_\k$}\hfill\\

As in the case of finite unitary groups, see \S\ref{sec:isocrystals}, our goal is to
compare the crystal of the categorical representation on  $\scrU_\k$
with the crystal of the $\frakg_{2d}$-module $\bigoplus_t\bfF(Q_t)_{2d}$. 
In order to get an explicit crystal isomorphism, we will assume that the 
$\ell$-decomposition matrices of unipotent blocks are unitriangular with respect to
the $a$-function. This gives a parametrization of the simple objects of $\scrU_\k$
and a combinatorial way (yet conjectural) to compute the (weak) Harish-Chandra branching graph.

\smallskip

Throughout this section we will assume that $\ell$ is a \emph{unitary prime}, which
means that the order $f$ of $q$ modulo $\ell$ is even. We will write $f = 2d$.

\subsubsection{Ordering symbols} For $t\geqslant 0$, we define the charge
\begin{equation}
\label{eq:charge_typeBC}
s_t = \left\{\begin{array}{ll} (t,d-1-t) & \text{ if $t$ is even}, \\
(t-d,-1-t) & \text{ if $t$ is odd}. \end{array} \right.\end{equation}
The symbol associated with a bipartition $\mu$ 
is $\Theta = \Theta_t(\mu) = \{\beta_t(\mu^1),\beta_{-1-t}(\mu^2)\}$.
As in \cite[\S 5.5.10]{GJ}, we can consider the union of the two $\beta$-sets $\beta_t(\mu^1)$ and
$\beta_{-1-t}(\mu^2)$ (allowing multiplicities) and write their elements in non-increasing
order
$$ \kappa_{m_t}(\mu) := (\kappa_1 \geqslant \kappa_2 \geqslant \kappa_3 \geqslant \cdots ), $$
where the parameter $m_t$ is given by
$$m_t = (t,-1-t).$$
We can check that $\kappa_i = \lfloor i/2 \rfloor$ for $i$ large enough. Then, we
define the $a$-value of the symbol $\Theta$ by
$$ a(\Theta_t(\mu)) = a(\Theta) := \sum_{i = 1}^{\infty} (i-1)\left(\kappa_i - \left\lfloor\frac{i}{2}\right\rfloor\right).$$
This differs slightly from the $a$-value $a^{(m_t,1)}(\mu)$ defined in \cite{GJ}~:
from \cite[prop.~5.5.11]{GJ} we have $a(\Theta_t(\mu)) = 
a^{(m_t,1)}(\mu) + a(\Theta_t(\emptyset))$. 

\smallskip

On the set of symbols in the same Harish-Chandra series and of the same rank we define as in \cite[\S 5.7.5]{GJ} a partial order by
$$ \Theta_t(\mu)  \preccurlyeq \Theta_t(\lambda)  \iff \mu \ll_{m_t} \lambda \iff \lambda = \mu \text{ or } \kappa_{m_t}(\lambda) 
< \kappa_{m_t}(\mu).$$
We deduce from \cite[prop.~5.7.7]{GJ} that
\begin{equation}\label{eq:a-compatible}
\Theta_t(\mu)  \prec \Theta_t(\lambda) \ \Longrightarrow \ a(\Theta_t(\mu)) < a(\Theta_t(\lambda)).\end{equation}

\subsubsection{Parametrization of unipotent modules}
By Proposition \ref{prop:isodecmap} we know that when $\ell$ is odd, there are as many
unipotent $\k G_n$-modules as symbols of odd defect and rank $n$. However there is
no known natural parametrization of these modules in terms of symbols unless the decomposition
matrix has unitriangular shape. This property has been proved when $f$ is odd
\cite{GrH}, i.e., when $\ell$ is a linear prime, or when $f$ and $n$ are small. Geck conjectured in
\cite[conj.~ 3.4]{GH97} that it should always hold.

\begin{conjecture}[Geck]\label{conj:unitriangular}
Assume $\ell$ is odd. There exists a parametrization 
$$\Irr(\scrU_\k)  = \{[D_\Theta]\,|\,\lambda\in\scrS_\text{odd}\}$$
of the simple unipotent modules and an ordering $\leqslant$ on symbols such that 
the decomposition map satisfies the following unitriangularity property
$$d_{\scrU}([E_\Theta]) \in \ [D_\Theta] + 
\sum_{\Theta' < \Theta}  \bbZ[D_{\Theta'}].$$
Moreover, the ordering $\leqslant$ can be chosen to be compatible with the $a$-function
so that $\Theta' \leqslant \Theta$ implies $a(\Theta') \leqslant a(\Theta)$.
\end{conjecture}

This conjecture is expected to hold for any finite reductive group as long as $\ell$
is a very good prime. In that case the ordering $\leqslant$ should be induced from the partial order
on families of unipotent characters defined by Lusztig.

\smallskip

\subsubsection{Comparison of the crystals}\label{sec:crystalBC}
Recall that to any categorical representation one can associate a perfect basis, and
hence an abstract crystal, see Proposition \ref{prop:PBfromcategorification}.
In \S \ref{ssec:cataction2d} we constructed a categorical representation of $\frakg_{2d}'$
on the category $\scrU_\k$ of unipotent representations over $\k$. From this
categorical representation we get an abstract crystal
$B(\scrU_\k)=\big(\Irr(\scrU_\k),\widetilde E_i,\widetilde F_i\big)$.
As for finite unitary groups, this crystal is related to the (weak) Harish-Chandra 
series. More precisely, the (uncolored) crystal graph associated with $B(\scrU_\k)$
coincides with the weak Harish-Chandra branching graph, and its connected components with
the weak Harish-Chandra series. See \S\ref{ssec:crystals-and-series} for more details.

\smallskip

Recall from Theorem \ref{thm:B} that $[\scrU_\k]$ is isomorphic to a direct
sum of Fock spaces $\bfF(Q_t)_{2d}$ where $t$ runs over $\bbN$ and
$Q_t = (Q_1,Q_2) = ((-q)^t,(-q)^{-1-t})$.
Recall from \eqref{eq:charge_typeBC} 
that we have defined the charge $s_t = (s_1,s_2)$
by
$$s_t=\begin{cases}(t,d-1-t)&\ \text{if}\ t\ \text{is\ even},\\
 (t-d,-1-t)&\ \text{if}\ t\ \text{is\ odd}\end{cases}$$
and $m_t =(t,-1-t)$. Since $q^d = -1$ in $\k$, we have $Q_p=q^{2s_p}$ for each $p=1,2$.
In other words $s_t$ is a charge for $\bfF(Q_t)_{2d}$ with respect to $q$.
Let $B(s_t)_{2d}$ be the corresponding abstract crystal of $\bfF(Q_t)_{2d}$,
with the canonical labeling
$B(s_t)_{2d}=\{b(\mu,s_t)\,|\,\mu\in\scrP^2,\,t\in\bbN\}.$
Finally, we set $B_{2d}=\bigsqcup_{\,t\in\bbN}B(s_t)_{2d}$.

\smallskip 

Under the assumption that the decomposition matrix is unitriangular with respect
to the $a$-function, we can compare the crystal $B_{2d}$ with the abstract crystal $B(\scrU_\k)$
coming from the categorical action of $\widehat \fraksl_{2d}$ on $\scrU_\k$.

\begin{theorem}\label{thm:crystalBC}
Recall that $\ell$, $q$ are odd and $f$ is even. If Conjecture \ref{conj:unitriangular} holds, then
the map $b(\mu,s_{t})\mapsto [D_{\Theta_t(\mu)}]$ is a crystal isomorphism
$B_{2d} \simto B(\scrU_\k)$.
\end{theorem}

\begin{proof}
Let $t \in \bbZ$ and $\mu$ be a bipartition.
By \cite[thm.~6.6.12]{GJ} the Uglov basis satisfies
$$ b(\mu,s_{t}) \in |\mu,s_t\rangle + \sum_{\begin{subarray}{c} \lambda \neq \mu \\ \lambda \ll_{m_t} \mu \end{subarray}}  \bbZ | \lambda,s_t\rangle.$$
Let $\phi : |\mu,s_t\rangle \mapsto d_{\scrU}([E_{\Theta_t(\mu)}])$ be 
the morphism of $\frakg_{2d}'$-modules studied in \S \ref{ssec:cataction2d}.
By the previous formula and the definition of $\preccurlyeq$ we have
$$ \phi(b(\mu,s_{t})) \in d_{\scrU}([E_{\Theta_t(\mu)}]) + \sum_{\begin{subarray}{c} \Theta' \neq \Theta_t(\mu) \\ \Theta' \preccurlyeq \Theta_t(\mu) \end{subarray}}  \bbZ d_{\scrU}([E_{\Theta'}]).$$
Let $n$ be the rank of the symbol $\Theta_t(\mu)$. If Conjecture \ref{conj:unitriangular}
holds, then we can choose a total order $\leqslant$ on symbols of rank $n$ and odd defect,
compatible with the $a$-function and for which the decomposition map $d_\scrU$ is unitriangular. In particular, using \eqref{eq:a-compatible} we have
$$ \Theta' \neq \Theta_t(\mu)  \text{ and } \Theta' \preccurlyeq \Theta_t(\mu) \ \Longrightarrow
a(\Theta') < a(\Theta_t(\mu)) \ \Longrightarrow \Theta' \leqslant \Theta_t(\mu).$$
This, together with the previous formula shows that
$$\phi(b(\mu,s_{t})) \in [D_{\Theta_t(\mu)}] + \sum_{\begin{subarray}{c} \Theta' \neq \Theta_t(\mu) \\ \Theta' \leqslant \Theta_t(\mu) \end{subarray}}  \bbZ [D_{\Theta'}]$$
and we can invoke Proposition \ref{prop:perfect-bases} to conclude.
\end{proof}

\begin{remark} Note that the definition of $s_t$ makes sense for negative integers
$t$ as well. Then $b(\mu,s_t) \longmapsto b(\mu^\dag,s_{-1-t})$ induces an isomorphism
of crystals. This is a particular case of changing the charge $(a,b)$ to $(b-d,a+d)$.
\end{remark}

\end{document}